\documentclass{amsart}
\usepackage{cite}
\usepackage{amsmath, amsthm, amssymb,tikz}
\usepackage{accents}
\usepackage{csquotes}
\usepackage{tikz-cd}
\usepackage{young,ytableau,youngtab}
\usepackage{MnSymbol}
\usepackage[hidelinks]{hyperref}
\usetikzlibrary{arrows}
\usepackage{longtable}
\usepackage{mathtools}
\usepackage{bbm}
\numberwithin{equation}{section}
\newtheorem{lemma}{Lemma}[section]
\newtheorem{theorem}[lemma]{Theorem}
\newtheorem{conjecture}[lemma]{Conjecture}
\newtheorem{corollary}[lemma]{Corollary}
\newtheorem{proposition}[lemma]{Proposition}
\theoremstyle{definition}
\newtheorem{example}[lemma]{Example}
\newtheorem{remark}[lemma]{Remark}
\newtheorem{definition}[lemma]{Definition}
\newtheorem{defproposition}[lemma]{Definition-Proposition}

\def\tr{{\rm tr}}
\def\16{{\bf 16}}
\def\1{{\bf 1}}
\def\2{{\bf 2}}
\def\3{{\bf 3}}
\def\4{{\bf 4}}
\def\val{\mathrm{val}}
\def\bar{\overline}
\def\tilde{\widetilde}
\def\hat{\widehat}

\def\Sp{{\mathrm{Sp}}}
\def\SL{{\mathrm{SL}}}
\def\GL{{\mathrm{GL}}}

\def\Sym{{\mathrm{Sym}}}
\DeclareMathOperator{\Ad}{Ad}
\DeclareMathOperator{\Hilb}{Hilb}

\newcommand{\tH}{\widetilde{H}}
\renewcommand{\th}{\tilde{h}}
\DeclareMathOperator{\Res}{\mathrm{Res}}
\DeclareMathOperator{\coarm}{\mathrm{coarm}}

\DeclareMathOperator{\Ind}{\mathrm{Ind}}

\DeclareMathOperator{\id}{\mathrm{id}}

\DeclareMathOperator{\Link}{Link}
\DeclareMathOperator{\Lie}{Lie}

\DeclareMathOperator{\Fl}{Fl}

\DeclareMathOperator{\Jac}{Jac}

\DeclareMathOperator{\Pic}{Pic}

\DeclareMathOperator{\cox}{\textbf{cox}}

\DeclareMathOperator{\ch}{\mathrm{ch}}

\DeclareMathOperator{\gr}{gr}
\DeclareMathOperator{\HHH}{\mathrm{HHH}}
\DeclareMathOperator{\SH}{\mathbb{SH}}
\DeclareMathOperator{\eee}{\mathbf{e}}

\DeclareMathOperator{\MM}{\mathbf{M}}
\newcommand{\cG}{\mathcal{G}}
\newcommand{\wt}{\mathrm{wt}}
\newcommand{\frs}{\mathbf{f}_{\vec{p},\vec{q}}}
\newcommand{\fpq}{\frs}
\newcommand{\fgamma}{\mathbf{f}_{\gamma}}

\newcommand{\tW}{\widetilde{W}}

\newcommand{\cO}{\mathcal{O}}

\newcommand{\cH}{\mathcal{H}}
\newcommand{\cE}{\mathcal{E}}

\newcommand{\cN}{\mathcal{N}}

\newcommand{\cU}{\mathcal{U}}

\newcommand{\fs}{\mathfrak{s}}
\newcommand{\fgl}{\mathfrak{gl}}

\newcommand{\fm}{\mathfrak{m}}
\newcommand{\fn}{\mathfrak{n}}

\newcommand{\fv}{\mathfrak{v}}
\newcommand{\fl}{\mathfrak{l}}
\newcommand{\fp}{\mathfrak{p}}

\newcommand{\fg}{\mathfrak{g}}

\newcommand{\bP}{\mathbf{P}}
\newcommand{\bO}{\mathbb{O}}

\newcommand{\bPic}{\overline{\Pic}}

\newcommand{\bI}{\mathbf{I}}

\newcommand{\C}{{\mathbb C}}
\newcommand{\R}{{\mathbb R}}
\newcommand{\N}{{\mathbb N}}
\newcommand{\A}{\mathbb{A}}
\newcommand{\Z}{{\mathbb Z}}
\newcommand{\G}{\mathbb{G}}
\newcommand{\Q}{{\mathbb Q}}
\newcommand{\F}{{\mathbb F}}

\newcommand{\HH}{\mathbb{H}}

\newcommand{\SYT}{\mathrm{SYT}}
\newcommand{\fqgamma}{\hat{\mathbf{f}}_\gamma}
\newcommand{\fqrs}{\hat{\mathbf{f}}_{\vec{p},\vec{q}}}
\newcommand{\fqpq}{\fqrs}

\newcommand{\St}{\mathrm{St}}
\newcommand{\ad}{\mathrm{ad}}
\newcommand{\aff}{\mathrm{aff}}

\DeclareMathOperator{\cchar}{\text{char}}
\DeclareMathOperator{\touch}{touch}
\newcommand{\Dyck}{\mathbb{D}}
\newcommand{\nDyck}{\mathbb{D}^{\mathrm{nest}}}

\DeclareMathOperator{\area}{\mathrm{area}}

\DeclareMathOperator{\coarea}{\mathrm{coarea}}

\newcommand{\NEpath}[4]{
    \fill[white!25]  (#1) rectangle +(#2,#3);
    \fill[fill=white]
    (#1)
    \foreach \dir in {#4}{
        \ifnum\dir=0
        -- ++(1,0)
        \else
        -- ++(0,1)
        \fi
    } |- (#1);
    \draw[help lines] (#1) grid +(#2,#3);
    \draw[dashed] (#1) -- +(#2,#3);
    \coordinate (prev) at (#1);
    \foreach \dir in {#4}{
        \ifnum\dir=0
        \coordinate (dep) at (1,0);
        \else
        \coordinate (dep) at (0,1);
        \fi
        \draw[line width=2pt,-stealth] (prev) -- ++(dep) coordinate (prev);
    };
}

\def\Young#1{\vbox{\smallskip\offinterlineskip
    \halign{&\vbox{##}\kern-\Thickness\cr #1}}}

\newdimen\Squaresize \Squaresize=12pt
\newdimen\Thickness \Thickness=.1pt
\newdimen\Correction \Correction=7pt

\def\Vide#1{\hbox{
       \vbox to \Squaresize{\vss
          \hbox to \Squaresize{\hss#1 \hss}\vss}
    \hskip-\Correction}
   \kern-\Thickness}
\def\Carre#1{\hbox{\vrule width \Thickness
   \vbox to \Squaresize{\hrule height \Thickness\vss
      \hbox to \Squaresize{\hss$\scriptstyle#1$\hss}
   \vss\hrule height\Thickness}
   \unskip\vrule width \Thickness}
   \kern-\Thickness}

\begin{document}
\title{Shalika germs for tamely ramified elements in $GL_n$}
\author{Oscar Kivinen}
\address{Aalto University}
\author{Cheng-Chiang Tsai}
\address{Academia Sinica and NSYSU}


\begin{abstract}
We prove explicit combinatorial formulas for various germ expansions of orbital integrals of tamely ramified elements in $GL_n(F)$, where $F$ is a nonarchimedean local field. The relevant combinatorics arises from the theory of the elliptic Hall algebra and the representation-theoretic knot superpolynomials defined by Cherednik--Danilenko and Morton--Samuelson. 

As a byproduct, we give explicit formulas for the weight polynomials of affine Springer fibers in type A and standard orbital integrals of tamely ramified regular semisimple elements. Our formulas subsume most earlier facts about the structure of Shalika germs of $GL_n$ in the literature.

As a further corollary, we show that point-counts of compactified Jacobians of locally planar curves are given by non-negative integral polynomials. Our results also provide further evidence for the Oblomkov-Rasmussen-Shende conjecture relating compactified Jacobians and HOMFLY-PT invariants of algebraic knots.

\end{abstract}

\maketitle 
\tableofcontents
\section{Introduction}
The Shalika germs of a reductive group over a non-archimedean local field are among the most fundamental invariants in its harmonic analysis, controlling the behaviour of orbital integrals near the identity; outside a handful of cases, however, they have resisted explicit computation. In this paper we determine them for every tamely ramified regular semisimple element of $GL_n(F)$. 
Let $G/F$ be a reductive group defined over a non-archimedean local field $F$. Given a compactly supported locally constant function $f\in C^\infty_c(G(F))$, the Shalika germ expansion of orbital integrals \cite{DeBacker, Shalika} can be stated as follows.
\begin{proposition}\label{prop:Shalika}
Given any element $\gamma\in U\subseteq G(F)$, where $U$ is a neighborhood of the identity depending on $f$, the orbital integral $I_\gamma(f)$ satisfies the expansion
$$I_\gamma(f)=\sum_{\bO\in \mathcal{U}(F)/\sim}\Gamma_\bO(\gamma)I_\bO(f),$$
where $\cU(F)/\sim$ is the set of unipotent orbits, $I_\bO$ are the corresponding unipotent orbital integrals, and $\Gamma_\bO: U \to \C$ are unique functions called Shalika germs.
\end{proposition}

Computing the Shalika germs for a given group $G/F$ is an important and in general open problem. The germs carry a great deal of harmonic-analytic information: they yield formulas for regular semisimple orbital integrals (Section \ref{sec:orbitalintegrals}) and are essential to the determination of character values of admissible representations \cite{Murnaghan, Spice}. In a more geometric direction, for $G=GL_n$ and $\gamma$ regular semisimple and tamely ramified, the orbital integrals $I_\gamma$ of suitable test functions compute weight polynomials of affine Springer fibers in type $A$ and of the associated compactified Jacobians (Theorem \ref{thm:weightpolynomials}).
Our formulas for these orbital integrals as well as the associated germ expansions is governed by the representation theory of the elliptic Hall algebra, and the rational functions that occur are closely related to the representation-theoretic knot superpolynomials of Cherednik--Danilenko \cite{CD} and Morton--Samuelson \cite{MS}. The elliptic Hall algebra \cite{BS1} was introduced with automorphic forms over function fields in mind, but its role in producing the Shalika germs of $p$-adic orbital integrals over an arbitrary local field is far from obvious. Moreover, our methods suggest intimate connections with the Hilbert schemes of points on $\A^2$ and HOMFLY-type knot invariants (see \cite{GKR} for a survey), in particular yielding a conjectural geometric incarnation of the germs.

The main results are stated in more detail in Subsection \ref{subsec:mainresults} below. 
We first provide some historical context. 

In \cite{W2}, Waldspurger defines another set of germs, arising from expanding the distribution $I_\gamma(-)$ in terms of truncated characters of the representations $\Ind_{P}^G(\St)$, where $P$ is a parabolic subgroup and $\St$ is the Steinberg representation of its Levi factor. The truncation is an essential technical tool, originally due to Clozel \cite{C89}. We refer to the resulting germs as {\em Steinberg germs} (see Theorem \ref{thm:steinberg} and Section \ref{sec:combinatorics}).

The paper \cite{W2} also suggests a recursive algorithm to determine the Steinberg germs for tamely ramified elements: 
by comparing the germ expansions for carefully chosen test functions, one shows that the germs of an element $\gamma\in G(F)$ should be determined inductively by those of
smaller depth elements $\gamma'\in GL_{n'}(E)$, where $n'|n$ and $E/F$ is a tamely ramified extension. 
However, the resulting algorithm is cumbersome to implement, and even then, one must understand the change of basis from truncated Steinberg distributions to unipotent orbital integrals in order to relate the Steinberg germs to the original Shalika germs.

In Section \ref{sec:combinatorics}, we recast this recursion, the change of basis, and a number of other results about germ expansions on $G(F)$ using the elliptic Hall algebra. 
This greatly clarifies the resulting computations, giving formulas essentially computable by hand and revealing several new structural properties.
As Waldspurger anticipates in \cite{W2}:

\begin{displayquote}
{\em L'auteur est
convaincu qu'il existe une bonne combinatoire, moins naïve que celle utilisée ici, qui
devrait permettre de calculer les germes.} \hspace{5cm} - J.-L. Waldspurger \cite{W2}
\end{displayquote} 

\subsection{Main results}
\label{subsec:mainresults}
Let $G=GL_n$, $\fg=\Lie(G)=\fg\fl_n$. As above, $F$ is a non-archimedean local field, with $\cO$ its ring of integers, $\mathfrak{m}_F$ the maximal ideal, and $k=\cO/\mathfrak{m}_F$ the residue field. We also choose a uniformizer $\varpi\in F$ and endow $F$ and any of its finite extensions with the standard valuation. Let $\gamma\in\fg(\cO)$ be a regular semisimple element that lives in a maximal torus that splits over a tamely ramified extension. Let $F(\gamma)$ be the commutative algebra generated by $\gamma$ in the algebra of matrices.

Our main results may be summarized as follows:
\begin{enumerate}
\item We attach to $\gamma$ an \emph{orbital symmetric function} $\fgamma$ whose expansions in two standard bases of the ring of symmetric functions (depending on a parameter) $\Sym_q$, are simultaneously the Shalika and the Steinberg germ expansions of $\gamma$ (Definition \ref{def:ellipticmastersymfn}). The classical change of basis between these two families of germs becomes an elementary change of basis in $\Sym_q$. There is a third expansion featuring in the proofs, related to the Shalika one by a Fourier transform and truncation to compact elements, which is manifestly positive and has the clearest combinatorial significance. 
\item For inertially elliptic $\gamma$, the orbital symmetric function coincides up to a change of variables with the {\em annular symmetric function} (Section \ref{sec:knots}). The latter is constructed via an iterated application of "plethystic" operators along the slopes determined by the Newton pairs of $\gamma$, rather than by a recursion (Theorem \ref{thm:mastersymfn}). The general elliptic case admits a similar, albeit slightly more complicated formula (Theorem \ref{thm:mainformula}).
\item Orbital integrals of parahoric test functions are recovered from $\fgamma$ through a Hall pairing against homogeneous symmetric functions, and for inertially elliptic $\gamma$, they recover weight polynomials of affine Springer fibers (Theorem \ref{thm:weightpolynomials}). These have closed formulas in terms of combinatorial data we call {\em $k$-nested Dyck paths}, arising from the Newton pairs of $\gamma$.
\item The Shalika germs satisfy an explicit transition rule under the passage from $GL_n(F)$ to $GL_{n/e}(F')$ along a totally tamely ramified extension (Theorem \ref{thm:transitionforSha}). 
In one of the standard normalizations, this transition matrix is given by an explicit combinatorial formula using graphs (Theorem \ref{thm:graphformula}).
\item As consequences we obtain the positivity and integrality of compactified-Jacobian point counts (Corollary \ref{cor:nonnegative}), a canonical $t$-deformation of the germs, and information about the finite part of the affine Springer action on the cohomology of the Iwahori affine Springer fibers.
\end{enumerate}
We now describe these results in more detail.
As above, the starting point is the introduction of a symmetric function in infinitely many variables $\fgamma\in \Sym_q$,
called the orbital symmetric function. 

By definition (see Definition \ref{def:ellipticmastersymfn}), $\fgamma$ is attached to any tamely ramified regular semisimple element $\gamma\in \fg\fl_n(F)$. It is defined by the equations
\begin{equation}
\fgamma:=\sum_{\lambda\vdash n} \Gamma_{\lambda^t}(\gamma)\th_\lambda=\sum_{\lambda\vdash n}\Gamma_\lambda^{St}(\gamma)e_\lambda
\end{equation}
where $\Gamma_\lambda(\gamma)$ are the Shalika germs from Proposition \ref{prop:Shalika} and $\Gamma_\lambda^{St}(\gamma)$ are the Steinberg germs from Theorem \ref{thm:steinberg}. The index $\lambda$ tracks nilpotent orbits via their Jordan type and standard parabolics via their block-triangular structure. The expressions $\th_\lambda, e_\lambda$ are bases of $\Sym_q$ defined in Section \ref{sec:symfns}. We note that the second equality is consistent with the Shalika--to--Steinberg change of basis on the harmonic analysis side (Fourier transform composed with a sign twist and truncation), and should be seen as its combinatorial manifestation. 

\begin{remark}
As we will show, the germs appearing in the expressions above are rational functions in $q$ without poles at $q=|k|$, the cardinality of the residue field. We therefore adapt the systematic abuse of notation that for the harmonic analysis expressions, $q$ also denotes $|k|$ and that $\fgamma$ denotes the specialization of $\fgamma$ at $q=|k|$. Unless explicitly stated, it should be clear from the context which one is meant.
\end{remark}

\begin{remark}
We note that \cite{W1,W2} work with the group $GL_n(F)$ rather than the Lie algebra $\fg\fl_n(F)$; the connection is straightforward and is explained in Section \ref{sec:orbitalintegrals}. 
We also note that much of the harmonic analysis literature imposes further conditions on $F$, such as on the characteristic or residue characteristic. Our results apply unconditionally; we have tried our best to explain the connections to existing results on germ expansions in Section \ref{sec:orbitalintegrals}.
\end{remark}


For the sake of simplicity, we will now assume that $\gamma$ is inertially elliptic, meaning that $F(\gamma)/F$ is a totally ramified degree $n$ extension, and refer to Remark \ref{rmk:RedtoLeviIntro} for the general case.
For inertially elliptic $\gamma$, we show in Corollary \ref{cor:inertiallyelliptic} that $\fgamma$ essentially coincides with a representation-theoretically defined symmetric function which we denote $\fpq$ and call the {\em annular symmetric function}. We now explain the construction; details will be given in Section \ref{sec:eha}.

With the above assumptions, $\gamma$ admits a Newton--Puiseux type expansion in terms of a uniformizer $u$ of $F(\gamma)$ with $u^n=:\varpi\in F$:
\begin{equation}
\label{eq:puiseux}
\gamma=\sum_{r\in\frac{1}{n}\Z_{\ge 0}} a_r u^{nr},
\end{equation}
where $a_i\in\cO^{\times}$. As a concrete example, one may take $F=k((\varpi))$ and $u$ to be of the form
\begin{equation}
\label{eq:umatrix}
u=\begin{pmatrix}
0 & 0 & \cdots & 0 & \varpi\\
1 & 0 & \cdots & 0 & 0\\
0 & 1 & \ddots & 0 & 0\\
0 & 0 & \ddots & 0 & 0\\
0 & \cdots & 0 & 1 & 0
\end{pmatrix}.
\end{equation}
From Eq. \eqref{eq:puiseux}, we extract the finite set of {\em root valuations}
\[
RV(\gamma):=\left\{r\in\frac{1}{n}\Z_{\ge 0}\;|\;r\not\in\mathrm{span}_{\Z}\{1,\,r'\in\frac{1}{n}\Z_{\ge 0}\;|\;r'<r,\;a_{r'}\not=0\}\right\}
\]
and define
\begin{definition}
\label{def:puiseux} 
Write $RV(\gamma)=\{r_1, r_2,\cdots,r_k\}$, where $r_1>\cdots>r_k$. Let $r_k=\frac{m_k}{n_k}$ be the reduced expression, and inductively write $r_i=\frac{m_i}{n_in_{i+1}\cdots n_{k}}$, where $m_i\in\Z_{>0}$, $n_i\in\Z_{\ge 0}$ are coprime. We call $(m_1,n_1),\ldots,(m_k,n_k)$ the sequence of {\em Puiseux pairs} of $\gamma$.

Given the Puiseux pairs, we define another sequence $(p_1,q_1),\ldots, (p_k,q_k)$ from the previous one by setting $p_i=n_i$, $q_i=m_i-m_{i+1}n_i$. Note that we always have $\gcd(p_i,q_i)=1$. We call $(p_1,q_1),\ldots,(p_k,q_k)$ the sequence of {\em Newton pairs} of $\gamma$. 
\end{definition}
Given any sequence of pairs of integers $(a_1,b_1),\ldots, (a_k,b_k)$, we will denote it by $(\vec{a},\vec{b})$. In subscripts, we will also drop the parentheses and simply write $\vec{a},\vec{b}$.
As explained in Section \ref{sec:eha}, for any pair of nonnegative coprime integers $a,b\in \Z_{\geq 0}$, one may define a homomorphism $\varphi_{a/b}: \Sym_{q}\to \Sym_q$ by letting
\begin{equation}\label{eq:phi_m/n}
\varphi_{a/b}(e_\ell)=E_{a,b,\ell}:=\sum_{\pi\in \Dyck_{\ell a,\ell b}}q^{\area(\pi)}e_\pi
\end{equation}
and extending by multiplication. Here $\Dyck_{\ell a,\ell b}$ is the set of Dyck paths in a $(\ell b\times \ell a)$ rectangle and $\area(\pi)$ is the area statistic (Definition \ref{def:dyck}). The $e_\ell$ and $e_\pi$ are elementary symmetric functions; note that that $\varphi_{a/b}(e_\ell)$ has degree $\ell b$. We refer to Proposition \ref{prop:Emnk} for the details.
\begin{definition}
We now define the annular symmetric function $\fpq$. 
\label{def:introannularsym}
Given a sequence of pairs of coprime integers $(\vec{p},\vec{q})=(p_1,q_1),\ldots,(p_k,q_k)$, we define the annular symmetric function $$\frs:=\varphi_{q_k/p_k}(\cdots \varphi_{q_1/p_1}(e_1))\cdots)$$
\end{definition}
\begin{remark}\label{rmk:slopeoperators}
The homomorphisms $\varphi_{a/b}$ are the $t=1$ shadows of the action of the slope subalgebras of the elliptic Hall algebra (Section \ref{sec:eha}) on the Fock space: for each slope the algebra contains a Heisenberg subalgebra, whose nonnegative half isomorphic to a ring of symmetric functions, and the change of slope is implemented by the $\widehat{\SL_2(\Z)}$-symmetry on the EHA. In the limit at $t=1$, the action on Fock space degenerates to give the degree-changing operators $\varphi_{a/b}$.
\end{remark}

\begin{remark}
\label{rmk:puiseux}
In the above discussion, apart from the discussion on germs, one may let $F$ be any complete DVF. When $F=\C((\varpi))$, the Newton--Puiseux sequence 
determines the topological type of the plane curve singularity $\{\det(xI-\gamma)=0\}\subset \C^2$, i.e. the isotopy class of its link in $S^3$. Conversely, if we impose $q_k\geq p_k$, the topological type uniquely determines the Puiseux pairs. The terminology {\em annular symmetric function} reflects the fact that after imposing this inequality, not only the link but the isotopy type of an annular braid closure inside the solid torus is determined by $\gamma$. Relatedly, Definition \ref{def:introannularsym} is closely related to an intermediate symmetric function appearing in the definition of the representation-theoretic {\em knot superpolynomials} from \cite{CD,MS} and was motivated by it.
We refer to Section \ref{sec:knots} for further discussion and references on the topological perspective.

Given a more general DVF and $\gamma\in\fg$ as above, the Puiseux pairs, the topological type, the semigroup, and the resolution graph of the spectral curve of $\gamma$, whenever they are defined, all determine each other and can be thought of as defining equisingularity classes of $\gamma$.
\end{remark}
Our first main result identifies the orbital symmetric function of $\gamma$ with the annular symmetric function of its Newton sequence:
\begin{theorem}
\label{thm:mastersymfn}
Let $F$ be a non-archimedean local field as before, so that $|k|=q<\infty$. Let $\gamma\in\fg(\cO)$ be inertially elliptic as before, and denote by $(\vec{p},\vec{q})$ be the sequence of Newton pairs of $\gamma$.
Then
\begin{equation}
\label{eq:totram}
\fgamma=q^{\Xi(\gamma)}\omega\frs|_{q\mapsto q^{-1}}
\end{equation}
where $\omega$ is the involution on $\Sym_{q}$ swapping $e_{\lambda}$ and $h_{\lambda}$ (see \S\ref{subsec:Sym}), $\frs|_{q\mapsto q^{-1}}$ denotes the image under $q\mapsto q^{-1}$, and $\Xi(\gamma)\in\Z$ is the dimension of the affine Springer fiber as in Definition \ref{def:RV}. In particular, the Shalika germs of $\gamma$ only depend on its Newton pairs. 
\end{theorem}

\begin{remark}\label{rmk:RedtoLeviIntro}
For more general tamely ramified elliptic elements, one needs to modify the statement of Theorem \ref{thm:mastersymfn} slightly, to account for intermediate unramified extensions. We will restate and generalize Theorem \ref{thm:mastersymfn} to this setting in Theorem \ref{thm:mainformula}.

Finally, if $\gamma$ is not elliptic but still tamely ramified and regular semisimple, it will belong to a Levi subgroup $L$ in which it is elliptic. 
Since $G=GL_n$, we have that $L\cong \prod_{i=1}^\ell GL_{\lambda_i}$ for some integer partition $\lambda\vdash n$, so that 
$\gamma=\gamma_1\oplus\cdots \oplus \gamma_\ell$ for some elliptic $\gamma_i\in GL_{\lambda_i}$. In Corollary \ref{cor:mastersymfnmanycomponents}, we show that $\fgamma$ is proportional to $\prod_i \mathbf{f}_{\gamma_i}$, the product of the orbital symmetric functions of the factors. When $F=k((\varpi))$ and $F(\gamma)$ is a product of totally tamely ramified extensions, this implies that $\fgamma$ is uniquely determined by the Puiseux series of the branches of the spectral curve, which correspond to the different blocks of $\gamma$. 
\end{remark}

We now explain how Theorem \ref{thm:mastersymfn} can be applied to give combinatorial formulas for orbital integrals, which is our second main result. As is by now well-known, in many cases orbital integrals can be reformulated in geometric terms as point-counts on affine Springer fibers. As we explain in Section \ref{sec:applications}, our combinatorial formulas have interesting algebro-geometric content as well.

For any partition $\lambda\vdash n$, let $\bP_{\lambda}\subset G(\cO)$ be the standard parahoric subgroup and $\1_\lambda\in C^{\infty}(\fg(\cO))$ the characteristic function of $\Lie\bP_{\lambda}$ divided by its measure.
As initially observed by Waldspurger in \cite{W1,W2}, it is natural to collect different germ expansions of $\gamma$ into a symmetric function. We claim that the orbital symmetric function $\fgamma$ is {\em the} canonical avatar of germ expansions, seen from various perspectives. On the one hand, $\fgamma$ encodes orbital integrals via the Hall pairing of symmetric functions, as we will show in Theorem \ref{thm:orbintformula}:

\begin{theorem}
\label{thm:weightpolynomials}
Let $\gamma\in \fg$ be a tamely ramified regular semisimple element and $\lambda\vdash n$ be a partition. 
Then 
$$ I_\gamma(\1_\lambda)=\langle \fgamma, h_\lambda\rangle$$ where we pair using the Hall inner product (Definition \ref{def:innerproducts}). 

When $\gamma$ is additionally inertially elliptic, we also have $I_{\gamma}(\1_\lambda)=|\Sp_\gamma^\lambda(k)/\Lambda_\gamma(k)|$,  where $\Sp_\gamma^\lambda$ is the affine Springer fiber of $\gamma$ in the partial affine flag variety $G(F)/\bP_\lambda$. See Sections \ref{sec:orbitalintegrals} and \ref{sec:applications} for the precise definitions and for our choice of normalization for the Haar measures used.
\end{theorem}
This yields the following combinatorial formula in the inertially elliptic case, which we prove as Corollary \ref{cor:orbitalnesteddyck}. Combining Theorem \ref{thm:mastersymfn}, Theorem \ref{thm:weightpolynomials} and the nested Dyck path expansion of the annular symmetric function (Proposition \ref{prop:annularformula}), the orbital integral of every standard parahoric test function is a manifestly non-negative sum over nested Dyck paths:
\begin{equation}
\label{eq:orbitalnesteddyckintro}
I_\gamma(\1_\lambda)=q^{\Xi(\gamma)}\sum_{\Pi\in\nDyck_{(\vec p,\vec q)}}q^{-\area(\Pi)}\,\bigl|\MM(\lambda,\lambda(\Pi))\bigr|.
\end{equation}
Here $\nDyck_{(\vec p,\vec q)}$ is the set of $k$-nested Dyck paths determined by the Newton pairs of $\gamma$ (Definition \ref{def:nesteddyck}), $\area(\Pi)$ is the total area of the tower, $\lambda(\Pi)\vdash n$ is the partition recorded by its top-level runs, and $|\MM(\lambda,\lambda(\Pi))|$ is the number of non-negative integer matrices with row sums $\lambda$ and column sums $\lambda(\Pi)$ (Definition \ref{def:rowcolmatrices}), which coincides with the Hall pairing $\langle h_\lambda,h_{\lambda(\Pi)}\rangle$. Note that in the spherical case $\lambda=(n)$ has $|\MM((n),\lambda(\Pi))|=1$ for every $\Pi$ and recovers the (reversed) area generating function of Corollary \ref{cor:sphericalorbital}.



Back on the side of invariant distributions, the Shalika germs satisfy a well-known homogeneity property (Proposition \ref{prop:homogeneity}), which for  $\fg=\fg\fl_n$ reads
$$\Gamma_{\lambda}(\varpi\gamma)=|\varpi|^{-n(\lambda^t)}\Gamma_\lambda(\gamma),$$
where $n(\lambda^t):=\sum_i(i-1)\lambda^t_i$. 
This has a clear parallel on the symmetric function side, where the $\nabla$-operator of Bergeron and Garsia (Definition \ref{def:nablaoperator}) scales each modified Macdonald polynomial $\tH_\lambda$ by $q^{n(\lambda)}t^{n(\lambda^t)}$. In particular, at $t=1$, it scales the corresponding specialization $\tH_\lambda[X;q,1]=:\th_\lambda$, by $q^{n(\lambda)}$. We show (Proposition \ref{prop:homog-orbital}) that these two facts are directly related through the equation
\begin{equation}
\label{eq:homogintro}
\nabla|_{t=1}\, \fgamma=\mathbf{f}_{\varpi\gamma}.
\end{equation}
In other words, applying $\nabla$ to the orbital symmetric function corresponds to multiplying $\gamma$ by a uniformizer. 
\begin{remark}
Again, there is a knot-theoretic manifestation of Eq. \eqref{eq:homogintro}. As explained in Proposition \ref{prop:homog-annular}, over $F=\C((\varpi))$, multiplication of $\gamma$ by $\varpi$ corresponds to adding a full twist to the link of the singularity. This is consistent with the expectations from conjectures in \cite{CD,ORS}.
\end{remark}

Finally, in Theorem \ref{thm:transitionmatrix}, we translate the application of the slope $m/n$-plethysm to the transition matrix of Shalika germs, under the passage from $GL_n(F)$ to $GL_{n/e}(F')$ for a totally tamely ramified extension $F'/F$ of degree $e$. 
We state a vague version of the result here for convenience, and refer to Section \ref{sec:w2} for the precise statement. 
\begin{theorem}\label{thm:transitionforSha}
There is an explicit transition matrix $M^{d/e}=(M^{d/e}_{\lambda,\lambda'})$, indexed by partitions $\lambda\vdash n$ and $\lambda'\vdash n/e$, that expresses the Shalika germs of $\gamma$ in terms of those of $\gamma'$, i.e.\ $\Gamma_{\lambda^t}(\gamma)=\sum_{\lambda'}M^{d/e}_{\lambda,\lambda'}\,\Gamma_{{\lambda'}^t}(\gamma')$. Here $\gamma'\in GL_{n/e}(F')$ is the next element in the recursion, $F'/F$ is totally tamely ramified of degree $e$ (or unramified when $d=0$), and $d$ is an integer determined by $\gamma$. In the ordinary normalization the entries are given by an explicit composition sum weighted by the coefficients $\wt_{d/e}(\alpha)$ of Eq.~\eqref{eq:compweightfull} (Theorem \ref{thm:transitionmatrix}); in the degenerate-Whittaker and renormalized normalizations of Remark \ref{rmk:shelstad}, they admit an equivalent, manifestly combinatorial description as a weighted sum over certain graphs (Theorem \ref{thm:graphformula}).
\end{theorem}

Everything stated so far concerns the specialization $t=1$. In fact, on the side of symmetric functions the entire construction has a natural lift in to general $t$:

The orbital symmetric function $\fpq$ is the $t=1$ limit of a $(q,t)$-symmetric function $\fqrs$, the \emph{refined annular symmetric function} (Definition \ref{def:fullmastersymfn}); accordingly the germs admit a canonical $t$-deformation (Theorem \ref{thm:tdeformation}, with the deformed germs defined in Eq. \eqref{eq:deffgammahat}). We expect the parameter $t$ to carry genuine meaning on both sides. Geometrically, $t$ records the perverse filtration on the cohomology of the compactified Jacobian of the spectral curve $\{\det(xI-\gamma)=0\}$, in the sense of \cite{ORS}, and setting $t=1$ forgets the filtration and returns the point count, that is, the orbital integral. Arithmetically, we expect the deformation to be a manifestation of parahoric versions of the local factors of \emph{orbital $L$-functions} of $GL_n$ emphasized by Arthur \cite{Arthur}. Over a function field the local factor of such an $L$-function at $\gamma$ is the Dedekind zeta function of the order $\cO_F[\gamma]$. Its 
value at $t=1$ is the spherical orbital integral $I_\gamma(\1_{(n)})$ of Corollary \ref{cor:sphericalorbital}. The parahoric versions of 
these $L$-functions appear to be new, and just as the parahoric orbital integrals at $t=1$ are all recovered from the single symmetric function $\fpq$ by pairing against the $h_\lambda$ (Theorem \ref{thm:weightpolynomials}), we expect all their $t$-deformations to be controlled by the single refined symmetric function $\fqrs$. From this vantage point the canonical deformation of a Shalika germ is not an artifact of the symmetric-function formalism, but we leave the detailed exploration for future work.
We also mention that the occurrence of these refined functions in the geometry of Hilbert schemes of points on $\A^2$ further suggests a categorification of the deformed Shalika and Steinberg germs, which we also hope to pursue in the future.

%

To conclude this introduction, we mention some further directions and limitations of our methods. In Section \ref{sec:cptjac}, we give direct applications to the theory of plane curve singularities. Namely, we show the non-negativity and integrality of point-counts of compactified Jacobians of integral planar curves (Corollary \ref{cor:nonnegative}). In Section \ref{sec:boundsOI}, we also discuss applications to bounding orbital integrals on $GL_n$, with applications to automorphic $L$-functions. 

We close by delineating the scope of our methods and the natural directions beyond it. The construction is specific to $GL_n$ and its inner forms, as it relies on the elliptic Hall algebra. For a general reductive group, the most accessible case seems to be that of $\gamma$ lying in a maximal unramified torus, which corresponds to the slope-zero part of our formalism. It seems likely that similar combinatorial results may be found for the germs, even without the additional organizational tool the EHA provides. Further, we expect the (annular) Hochschild homology of Rouquier complexes of Soergel bimodules to play the role of the annular symmetric function $\fpq$ in that setting. The restriction to tame ramification is also essential for us throughout: the basic recursive step underlying $\fgamma$ breaks down for wildly ramified extensions, and it would be very interesting to know whether a finite algorithm computes the Shalika germs of an arbitrarily ramified $\gamma$ in terms of a discrete set of invariants, such as the Artin conductor.

\subsection{Outline of the paper}
The paper is organized as follows. Section \ref{sec:orbitalintegrals} reviews orbital integrals and the various germ expansions, and Section \ref{sec:symfns} collects the background on symmetric functions. Section \ref{sec:eha} introduces the elliptic Hall algebra, constructs the annular symmetric functions, and explains their relation to HOMFLY-type knot invariants and the representation-theoretic superpolynomials. Section \ref{sec:combinatorics} is the technical heart of the paper, where the closed formulas for the orbital symmetric function and the Shalika germs are proved. A reader interested only in the consequences for orbital integrals may read it selectively and pass to the applications. Section \ref{sec:examples} is devoted to worked examples, and Section \ref{sec:applications} develops the applications to plane curve singularities, bounds on orbital integrals, and orbital $L$-functions.

\subsection*{Acknowledgments}
The authors thank Francois Bergeron, Pierre-Henri Chaudouard, Stephen DeBacker, Eugene Gorsky, Thomas Hales, Bertrand Lemaire, Yen-Chi Roger Lin, Anton Mellit, Fiona Murnaghan, Andrei Neguț, Alexei Oblomkov, Koji Shimizu, Yan Soibelman, Loren Spice, Minh-Tam Trinh, Jean-Loup Waldspurger, Zhiwei Yun, and Wei Zhang for interesting conversations.
\section{Orbital integrals and germ expansions}
\label{sec:orbitalintegrals}

In this section, we set up the harmonic analysis background needed for the rest of the paper. We proceed as follows. After fixing normalizations of measures, we recall the Shalika germ expansion (Subsection \ref{subsec:shalikasubsec}) and introduce the Steinberg germs (Subsection \ref{subsec:steinbergsubsec}). In Subsection \ref{sec:germsAHA} we give an alternative, characteristic-independent construction of the relevant distributions using the cocenter of the affine Hecke algebra. The main result used later is Theorem \ref{thm:pseudo-steinberg}, which constructs a replacement for the Steinberg distributions in arbitrary characteristic. In Subsection \ref{subsec:paraind} we recall the parabolic induction/restriction formalism that reduces the computation of germs for non-elliptic elements to the elliptic case. 

Let $F$ be a non-archimedean local field with residue field $k$ and $\cO$ its ring of integers. Below, the group $G$ will be $G=GL_n$, or a product $G=\prod GL_{n_i}$ for which one can reduce the discussion to the former case. 
Let $C_c^\infty(G(F))$ be the ring of compactly supported, locally constant complex valued functions on $G(F)$.
For any $\gamma\in G(F)$ and $f\in C_c^{\infty}(G(F))$, the orbital integral is formally defined as
\begin{equation}
I_\gamma(f):=\int_{g\in C_{G(F)}(\gamma)\backslash G(F)}f(g^{-1}\gamma g)dg
\end{equation}

To make sense of this expression, one needs to explain the measures used.
For general $\gamma\in G(F)$ or $\gamma\in\fg(F)$, the conjugacy class of $\gamma$ is locally closed and the centralizer unimodular. Therefore, the orbit of $\gamma$ under the adjoint action admits an invariant measure. Another way to see this is to identify the orbit with the orbit inside $\fg\fl_n(F)^*$ via the natural equivariant embedding of varieties $GL_n\hookrightarrow \fg\fl_n$ and the Killing form. The coadjoint orbit admits a natural algebraic symplectic structure, whose top wedge is an invariant volume form defined up to a scalar.

We will not attempt to fix the normalization for arbitrary $\gamma$, but do so for semisimple and unipotent elements, since these two cases will be used throughout the paper.
\begin{definition}
\label{def:haarmeasures}
For $\gamma\in G(F)$ semisimple, the measure $dg$ is normalized as follows: On $G(F)$ we have the up to a scalar unique Haar measure, which we will normalize so that $G(\cO)$ has measure $1$. In this case, $C_{G(F)}(\gamma)$ is a product of general linear groups over extensions $F'$ of $F$, and we use the same normalization. The measure $dg$ in $I_\gamma(f)$ is the quotient of these two Haar measures, i.e. the unique $G(F)$-invariant measure on $C_{G(F)}(\gamma)\backslash G(F)$ such that for any measurable $\phi$, $$\int_{G(F)}\phi(g)dg=\int_{C_{G(F)}(\gamma)\backslash G(F)}\int_{C_{G(F)}(\gamma)}\phi(zg)dz\,dg$$ The same is done in the Lie algebra case.
\end{definition}

When $\gamma\in G(F)$ is unipotent, we will normalize the measure as follows, following \cite{W1, H74}. Let $\lambda\vdash n$ and consider the Richardson parabolic $P=P(\lambda^t)$ for the unipotent orbit associated to $\lambda$. Let $N$ be its unipotent radical and $M$ its Levi factor. Let $I_{\lambda}^G$ be the orbital integral on the unipotent orbit with Jordan type $\lambda$. For $f\in C_c^\infty(G(F))$ let $f^P\in C_c^\infty(M(F))$ be defined by
\begin{equation}
\label{eq:constantterm}
f^P(m)=\delta_P(m)^{1/2}\int_{G(\cO)}\int_{N(F)} f(k^{-1}mnk)dndk
\end{equation}
where $\delta_P$ is the modular function for $P$.
By \cite[Proposition 5]{H74}, the linear forms given by the unipotent orbital integrals $I^G_\lambda(-)$ for $\lambda\vdash n$ are proportional to $f\mapsto f^P(1)$. We normalize the measure on the unipotent orbit so that $I^G_\lambda(f)=f^P(1)$, where on $N(F)$ we take the Haar measure such that $N(\cO)$ has measure $1$.
This discussion applies verbatim with $G$ replaced by $\fg$, etc., and the resulting formula is
\begin{equation}
\label{eq:constanttermLie}
f^{\fp}(X)=\int_{G(\cO)}\int_{\fn(F)}f(\mathrm{ad}(k^{-1})(X+Y))dYdk
\end{equation}
and 
we normalize the nilpotent orbital integral $I^{\fg}_\lambda(f)=f^{\fp}(0)$, again with the Haar measure on $\fn(F)$ such that $\fn(\cO)$ has measure $1$. We will drop the superscripts $G, \fg$ when understood from the context. This is even more merited in view of 
\begin{lemma}
    \label{lem:equalityofintegrals}
    Let $f\in C^\infty(\fg(\cO)/\Lie(\bI))$ be the characteristic function of a standard parahoric.
    Then the restriction of $f$ to $G(\cO)\subset \fg(\cO)$ equals the characteristic function of the corresponding standard parabolic subgroup, and $$f^P(1)=f^{\fp}(0).$$ 
\end{lemma}

See Lemma \ref{lem:definitionsagree} where the content here is discussed in the more general setting in \S\ref{subsec:paraind} in terms of parabolic induction/restriction.

\subsection{Shalika germs}\label{subsec:shalikasubsec}

For any subset $\Sigma\subset G(F)$ we denote by $J(\Sigma)$ the space of invariant distributions on $G(F)$ supported on elements of the form $g^{-1}\gamma g$ with $g\in G(F)$, $\gamma\in\Sigma$. The famous Howe conjecture states that 

\begin{theorem}[\cite{H74}, \cite{HC}, \cite{C89}, \cite{BM00}] For any compact subset $\Sigma\subset G(F)$ and an open subgroup $K \subset G(F)$, the restriction of $J(\Sigma)$ to $C_c(G(F)/K)$ is finite-dimensional.
\end{theorem}

The same is true when we replace $\Sigma$ by a compact subset in $\fg(F)$, $K$ by an open sub-$\cO$-module in $\fg(F)$ and $C_c(G(F)/\Sigma)$ by $C_c(\fg(F)/K)$. In this article we will make use of precise versions of the above finiteness.
A particularly important one is the following theorem, proved by Hales \cite[Thm. 1]{Hales} for the span of regular semisimple orbital integrals, and which follows in general from Proposition \ref{prop:uniquedual} (see also Court\'{e}s \cite[Thm 1.10]{Cou06}):

\begin{theorem}
\label{thm:uni} 
Recall that $G=GL_n$. 
Let $\cU\subset G(F)$ be the ($F$-points of the) unipotent variety and $\bI\subset G(F)$ be an Iwahori subgroup. Then the restriction of $J(G(\cO))$ to $C_c(G(F)/\bI)$ is equal to that of $J(\cU)$ to $C_c(G(F)/\bI)$. Both restrictions have a basis given by unipotent orbital integrals.
\end{theorem}

\begin{remark}
\label{rmk:characteristic}
Hales \cite{Hales} works over characteristic zero $F$, as does Waldspurger \cite{W1, W2}. However, the most essential ingredient for Hales is the original Shalika expansion \cite{Shalika}, which also works in positive characteristic assuming finiteness of unipotent orbits and convergence of the unipotent orbital integrals. This was proved in \cite{McNinch} when the characteristic is {\em good} for $G$. In particular, for $GL_n$ it holds in arbitrary characteristic. Moreover, the more general Proposition \ref{prop:uniquedual} works in arbitrary characteristic, so Theorem \ref{thm:uni} is true in this generality.
\end{remark}

\begin{remark} As far as invariant distributions are concerned, any test function in $C_c(G(F)/K)$ can be averaged by $K$-conjugation into $C_c(K\backslash G(F)/K)$. Similarly, in Theorem \ref{thm:uni} we can replace $C_c(G(F)/\bI)$ by $C_c(\bI\backslash G(F)/\bI)$.
\end{remark}

Theorem \ref{thm:uni} is a variant of the Shalika germ expansion, and it can be reinterpreted as follows.
\begin{corollary}
\label{cor:shalika} 
For any $\gamma\in G(\cO)$, there exist constants $\Gamma_{\lambda}(\gamma)$ where $\lambda$ runs over unipotent orbits of $G(F)$ (i.e. partitions of $n$) such that for any $f\in C_c(G(F)/\bI)$, we have
\begin{equation}\label{eq:germ-expansion}
I_\gamma(f)=\sum_{\lambda}\Gamma_{\lambda}(\gamma)I_\lambda(f).
\end{equation}
\end{corollary}

The Lie algebra variant is as follows. 
\begin{theorem}\label{thm:shalika-for-Lie-algebra} 
 [Thm. 2.1.5., \cite{DeBacker}] Let $\fg$ be a reductive Lie algebra. Let $\cN\subset\fg(F)$ be the ($F$-points of the) nilpotent cone and $\Lie \bI\subset \fg(F)$ be an Iwahori subalgebra. Then the restriction of $J(\fg(\cO))$ to $C_c(\fg(F)/\Lie \bI)$ is equal to that of $J(\cN)$ to $C_c(\fg(F)/\Lie \bI)$. Both restrictions have a basis given by nilpotent orbital integrals.
\end{theorem}

As explained in \cite{DeBacker}, this statement is true for an arbitrary connected reductive group $G$ provided that $\operatorname{char} k\gg\operatorname{rank} G$.
The sharp bound needed when $G=GL_n$ is $\operatorname{char} k>2n$. Whenever the theorem works, it asserts the existence of unique functions $\Gamma_{\lambda}(\gamma)$ such that (\ref{eq:germ-expansion}) holds for $f\in C_c(\fg(F)/\Lie \bI)$ where $I_{\lambda}(f)$ is the integral on nilpotent orbits this time; by abuse of language we will denote the germs again by $\Gamma_{\lambda}(\gamma)$. This does no harm, thanks to the following proposition:

\begin{proposition}
\label{prop:germcomparison}
Let $\gamma\in\fg(\cO)$. Fix a nilpotent orbit $\lambda$ in $\fg(F)$ -- it corresponds to a unipotent orbit in $G(F)$ under $x\mapsto 1+x$.
Under this matching, the following are equal:
\begin{enumerate}
\item The Lie algebra Shalika germ $\Gamma_{\lambda}(\gamma)$ from Theorem \ref{thm:shalika-for-Lie-algebra}.
\item The Lie algebra Shalika germ $\Gamma_{\lambda}(c+\gamma)$ from Theorem \ref{thm:shalika-for-Lie-algebra}, for any $c\in\cO$.
\item The Lie group Shalika germ $\Gamma_{\lambda}(c+\gamma)$ from Theorem \ref{thm:uni} and Corollary \ref{cor:shalika}, for any $c\in\cO$ such that $c+\gamma\in G(\cO)$.
\end{enumerate}
\end{proposition}
\begin{proof}
(1) equals (2) since the central translation does not affect orbital integrals.

Just like the unipotent orbital integrals on $G(F)$ may be computed using Eq. \eqref{eq:constantterm}, so may the nilpotent ones on $\fg(F)$ using Eq. \eqref{eq:constanttermLie}. By Lemma \ref{lem:equalityofintegrals} the  restriction to $G(\cO)$ of the characteristic function of any standard parahoric subalgebra of $\fg(F)$ is the the characteristic function of the corresponding standard parahoric subgroup and the nilpotent orbital integral equals the unipotent orbital integral of the restriction.
By \cite[Corollaire 4.4.]{W1}, the nilpotent and unipotent orbital integrals are determined as distributions on Iwahori-(bi-)invariant functions by their values on $\1_{\Lie(\bP_\lambda)}\in C_c(\fg(\cO)/\Lie \bI)$ where $\lambda\vdash n$. Therefore (2) and (3) have the same Shalika germs, either for the group or the Lie algebra. Therefore (2) equals (3).\end{proof}


\begin{remark} Theorem \ref{thm:shalika-for-Lie-algebra} is expected to hold even if $\cchar k\leq 2n$. For example, in \cite[\S5]{Lemaire} it is shown that (\ref{eq:germ-expansion}) holds for any fixed $f$ and for fixed $\gamma$ which is {\em quasi-regular} (see the abstract to \cite{Lemaire}, and note in particular that when $\cchar k>n$ quasi-regularity is equivalent to regular semisimplicity) and in a small enough neighborhood $V_f\ni \gamma$. Consequently, (\ref{eq:germ-expansion}) holds for quasi-regular elements in a small enough neighborhood and for all functions in $C_c(\fg(\cO)/\Lie \bI)$, since the latter is finite-dimensional. 

This is what we used in the proof of Proposition \ref{prop:germcomparison}, and thus our Shalika germ results for tamely ramified regular semisimple $\gamma$ will also work for Lemaire's germs. Even more generally, in \cite[\S5.2, pp. 505]{Lemaire} Lemaire defines normalized Shalika germs $\tilde{b}_i$ defined on the set of all quasi-regular elements using homogeneity. Proposition \ref{prop:homog-orbital} together with propositions \ref{prop:homogeneity} and \ref{prop:germcomparison} show that such $\tilde{b}_i$ also agree with the Lie algebra/group Shalika germs and in particular the above results, as well as the computations in Section \ref{sec:combinatorics} extend to this case as well given $f\in C_c(\fg(\cO)/\Lie \bI)$.
\end{remark}


As mentioned in the introduction, the Shalika germs enjoy a homogeneity property, which we now restate for future reference. Consider a nilpotent orbit $\mathbb{O}\subset\fg(F)$.
\begin{proposition}
\label{prop:homogeneity}
For $G=GL_n$, we have 
$$\Gamma_{\mathbb{O}}(\varpi\gamma)=|\varpi|^{-\frac{1}{2}\dim \mathbb{O}}\Gamma_\mathbb{O}(\gamma).$$
\end{proposition}

\begin{remark} For general reductive group when $\cchar k$ is very good we have $\Gamma_{\mathbb{O}}(\varpi^2\gamma)=|\varpi|^{-\dim \mathbb{O}}\Gamma_{\varpi^{-2}\mathbb{O}}(\gamma)=|\varpi|^{-\dim\mathbb{O}}\Gamma_\mathbb{O}(\gamma)$ since $\mathbb{O}$ and $\varpi^{-2}\mathbb{O}$ are the same orbit in $\fg(F)$. However $\mathbb{O}$ and $\varpi^{-1}\mathbb{O}$ are typically different orbits and the identity in Proposition \ref{prop:homogeneity} does not hold in general.
\end{remark}
\begin{remark}
In terms of partitions, if $\mathbb{O}\leftrightarrow \lambda\vdash n$, 
$$\frac{1}{2}\dim\mathbb{O}=\sum_{i=1}^{\ell(\lambda^t)}(i-1)\lambda_i^t=:n(\lambda^t)$$
\end{remark}

\subsection{Steinberg germs}\label{subsec:steinbergsubsec}

We now define the Steinberg germs, which are one of the main ingredients in our recursive computation of $\fgamma$.
They arise through the expansion of $I_\gamma$ via truncated Steinberg distributions. To explain the truncation, we first recall various notions of compact elements.

\begin{definition}
\label{def:compactelts}
Let us say an element $g\in G(F)=GL_n(F)$ is {\em compact mod center} if all its eigenvalues have the same valuation; we will suppress the ``mod center'' and just call them {\em compact} when no confusion should arise. Let $G(F)_c\subset G(F)$ be the subset of all compact elements. 

Let us also call $\gamma\in \fg(F)$ compact if it belongs to some parahoric subalgebra, or equivalently that it is conjugate to an element in $\fg(\cO)$.
\end{definition}
\begin{remark}
    Note that in the group case, not all compact-modulo-center elements are literally central translations of compact elements. Instead, they become compact in the usual sense under the map to $PGL_n$.
\end{remark}
We will also need the following alternative characterization of the compact modulo center elements in $GL_n$:
\begin{lemma}
\label{lem:compactelts}
The set $G(F)_c$ coincides with the union of conjugates of normalizers of standard parahorics. In particular, for $GL_n$ we have 
$$G(F)_c=\Ad(G(F))\left(\bigcup_{\substack{e|n,\alpha\vdash n/e,\\ d\in \Z}} u^{nd/e} \bP_{\alpha^e} \right)$$ where $u$ is the matrix from the introduction and $\bP_{\alpha^e}$ is a standard parahoric subgroup. 
\end{lemma}
\begin{remark}
    Note that $u^{nd/e}, d\in \Z$, normalizes $\bP_{\alpha^e}$. Compare also to the beginning of \cite[Section 4]{W2}.
\end{remark}

Unsurprisingly, the Steinberg germs inherit their name from the Steinberg representation:
\begin{definition}
Let $\St_n$ be the Steinberg representation of $G(F)=GL_n(F)$. More generally, for $\lambda\vdash n$ (Definition \ref{def:partition}) let $P(\lambda)\subset GL_n(F)$ be the corresponding parabolic and $L(\lambda)$ its Levi subgroup. Let $\St_{\lambda}$ be the parabolic induction\footnote{Usually parabolic inductions are normalized by the modulus character. In our case, the modulus character is trivial on compact (mod center) elements because it is a homomorphism to $(\R^+,\times)$ that is trivial on center. Since we will immediately restrict to compact elements, the normalization will have no effect on what follows.} 
of the Steinberg representation of $L(\lambda)$ to $G(F)$
\end{definition}
We will identify $\St_{\lambda}$ with its character, an invariant distribution on $G(F)$. 
To define the germs, we denote by $\St_{\lambda,c}$ the restriction of $\St_{\lambda}$ to $G(F)_c$, i.e. the distribution which is truncated to be $0$ outside $G(F)_c$. We also denote by $J_{\St,c}\subset J(G(F)_c)$ the subspace spanned by $\St_{\lambda,c}$. 
Write $C_c(\bI\backslash G(F)/\bI)=\bigoplus_{k\in\Z} C_c(\bI\backslash G(F)^{\operatorname{val}=k}/\bI)$ where $G(F)^{\operatorname{val}=k}$ is the subset of elements whose determinant has valuation $k$. 
In \cite[Proposition 2.4.]{W1} and \cite[Proposition III 4.]{W2}, when $\cchar F=0$, Waldspurger proved:

\begin{theorem}\label{thm:steinberg}
\label{Wal} For $k\in\Z$ the restriction of $\St_{\lambda,c}$ to $C_c(\bI\backslash G(F)^{\operatorname{val}=k}/\bI)$ is non-zero iff $\lambda$ is divisible by $n/\gcd(k,n)$. The restriction of $J_{\St,c}$ to $C_c(\bI\backslash G(F)^{\operatorname{val}=k}/\bI)$ has a basis given by the restrictions of these $\St_{\lambda,c}$.
\end{theorem}

\begin{theorem}\label{thm:steinberg2}
For any regular semisimple $\gamma\in G(F)^{\operatorname{val}=k}$ that is compact mod center, there exist unique constants $\Gamma_\lambda^{St}(\gamma)$ indexed by $\lambda\vdash\gcd(k,n)$, so that for any $f\in C_c(\bI\backslash G(F)^{\operatorname{val}=k}/\bI)$ we have 
\begin{equation}\label{eq:St-germ}
I_\gamma(f)=\sum_{\lambda}\Gamma^{St}_\lambda(\gamma)\St_{n'\lambda,c}(f)
\end{equation}
where $n':=n/\gcd(k,n)$. We shall call the constants $\Gamma^{St}_\lambda(\gamma)$ the {\em Steinberg germs of }$\gamma.$
\end{theorem}

\begin{remark}
A few remarks are in order. 
The most important one is that in \cite{W1,W2} there is an assumption on $\cchar F=0$. However, it is easy to see that the proof of \cite[Proposition 2.4.]{W1} only uses characteristic-independent facts about the representation theory of $G(F)$. In the next section, Section \ref{sec:germsAHA}, we will in particular construct unique distributions $\St_{\lambda,c}$
satisfying Theorems \ref{thm:steinberg}, \ref{thm:steinberg2} and \cite[V 11, V12]{W2}. {\em If} one were able to carry out Clozel's work in \cite{C89} in positive characteristic, then these $\St_{\lambda,c}$ could safely be identified with the truncated characters of parabolic inductions of Steinberg representations. While it is somewhat awkward we cannot do this right now, it will not affect the computation of the Shalika germs themselves.

Further, we note that the "$St$"-superscript stands for Steinberg, and should not be confused with the notion of "stability" in the automorphic forms literature.

\end{remark}

\subsection{Germs via rigid cocenters of affine Hecke algebras}
\label{sec:germsAHA}

In this subsection, we construct the Steinberg distributions from the previous subsection in a characteristic-independent way, using the cocenter of the (extended) affine Hecke algebra $$\cH:=C_c(\bI\backslash G(F)/\bI)$$ following \cite{CiHe, HeCocenter}. 

Note that as an abstract algebra, $\cH$ only depends on the cardinality of the residue field of $F$, through for example the well-known description by generators and relations.  In particular, the results in this subsection apply in all characteristics. 

Consider the space $J(G(F)_c)$ of invariant distributions supported on the compact mod center elements of $G(F)$. There is a natural map $$J(G(F)_c)\rightarrow\cH^*$$ to the linear dual of the AHA given by evaluating the distributions on the functions. 
The $G(F)$-invariance of the distributions amounts to this map factoring through the {\em cocenter}
\begin{equation}\label{eq:embed}
J(G(F)_c)\rightarrow(\cH/[\cH,\cH])^*\rightarrow\cH^*
\end{equation}
We will soon see that the first map in Eq. \eqref{eq:embed} can be further shown to factor through the dual of the {\em rigid} cocenter introduced in \cite{CiHe}.
First, note that $\cH$ and $\tr(\cH):=\cH/[\cH,\cH]$ are graded by the valuation of the determinant. In particular, we have the graded algebra decomposition $$\cH=\bigoplus_{k\in\Z} C_c(\bI\backslash G(F)^{\operatorname{val}=k}/\bI)$$ and similarly for the cocenter. Let us denote $\cH^{\val=k}$, $\tr(\cH)^{\val=k}:=\tr(\cH^{\val=k})$ the corresponding subspaces. 

This is further refined by the {\em Newton decomposition} of the group $G(F)$ as well as $\cH$ from \cite{HeCocenter}. From the Cartan decomposition, we have $\bI\backslash G(F)/\bI\cong\tilde{W}$, where $\tilde{W}\cong\Z^n\rtimes S_n$ is the extended affine Weyl group associated to $G=GL_n$. We also write $W^{\mathrm{fin}}:=S_n$ in the above. We will identify $X_*(T)\cong \Z^n$ in a standard way.
We have a decomposition $\tW\cong \Omega\ltimes W^{\aff}$ where $W^{\aff}$ is the affine Weyl group and $\Omega\cong \Z$.
\begin{definition}
    The {\em Kottwitz map} is the projection $\kappa:\tilde{W}\twoheadrightarrow\Omega$.
    The {\em Newton map} $\nu': \tW\to \frac{1}{n}\Z^{n}$ is defined as follows. If $w^k.x=\lambda+x$ we let 
    $\nu'(w)=\lambda/k$. By sorting, we may take this to be the unique dominant element in the $W^{\mathrm{fin}}-$orbit of $\nu'(w)$, which gives another map $\nu^+:\tW\to (\frac{1}{n}\Z^n)_{\ge 0}$.

    Together, we get a map $$\pi:=(\kappa,\nu^+): \tW\to \Omega\times (\frac{1}{n}\Z^n)_{\ge 0}=:\aleph$$
    Note that $\pi(w)=\pi(w')$ whenever $w$ and $w'$ are conjugate.
\end{definition}
\begin{definition}
    For $\nu\in\aleph$, the {\em Newton stratum} of $G(F)$ is $$G(\nu):=\bigcup_{w\in \tW \text{ minimal length, } \pi(w)=\nu} \Ad(G(F))\bI w\bI$$
\end{definition}
\begin{theorem}[Theorem A, \cite{HeCocenter}]
\label{thm:groupnewton}
    We have the {\em Newton decomposition of $G(F)$} is $$G(F)=\bigsqcup_{\nu\in \aleph} G(\nu)$$ i.e. the group decomposes as a disjoint union of Newton strata.
\end{theorem}
For the algebraic/combinatorial point of view, we have the following theorem of He and Nie \cite[Thm. 6.7]{HN14} and He \cite[Thm. 11]{HeCocenter}:
\begin{theorem}
\label{thm:HN}
    We have a Newton decomposition for the cocenter $\tr(\cH)$ of the Iwahori-Hecke algebra $\cH$
    \begin{equation}\label{eq:Newton-trace}
    \tr(\cH)=\bigoplus_{\nu\in\aleph}\tr(\cH)_\nu
    \end{equation}
    where $\tr(\cH)_\nu$ is spanned by the images of the Iwahori-Matsumoto generators $T_w$ with $\pi(w)=\nu\in \aleph$ and $w$ is of minimal length. Moreover, the image of $T_w$ for minimal length $w$ depends only on its conjugacy class. These $T_w$, one for each conjugacy class with $\pi(w)=\nu$, form a basis of $\tr(\cH)_\nu$.
\end{theorem}

We are ready to define the {\em rigid cocenter} following \cite{HeCocenter} and its relation to $J(G(F)_c)$. Let us call $\nu^+\in(\frac{1}{n}\Z^n)_{\ge 0}$ central if it lives in the diagonal $\frac{1}{n}\Z_{\ge 0}$.

\begin{definition}
\label{def:rigtr}
    The {\em rigid cocenter} of $\cH$ is $$\tr(\cH)^{rig}:=\bigoplus_{\substack{\nu=(\kappa,\nu^+)\in \aleph, \\\nu^+ \text{ central}}} \tr(\cH)_\nu$$
\end{definition}
By \cite[Prop. 21]{HeCocenter}, $\tr(\cH)^{rig}$ is exactly the image of the subspace of $\bI$-bi-invariant $C^\infty$-functions on $G(F)$ represented by functions supported on the compact-mod-center elements. More precisely, we have
\begin{proposition}
\label{prop:cptnewton}
    The set $G(F)_c\subset G(F)$ of compact-mod-center elements in $G(F)$ is exactly 
    $$G(F)_c=\bigsqcup_{\substack{\nu=(\kappa,\nu^+)\in \aleph, \\\nu^+ \text{ central}}} G(F)_\nu$$
\end{proposition}

\begin{corollary}
\label{cor:factors} Identify $\tr(\cH)^{rig}$ as a direct summand of $\tr(\cH)$ (as vector spaces) using \eqref{eq:Newton-trace}. Then the map \eqref{eq:embed} factors as $$J(G(F)_c)\to (\tr(\cH)^{rig})^*\hookrightarrow (\tr(\cH))^*$$
In particular, the image of a distribution in $J(G(F)_c)$ under the map \eqref{eq:embed} is determined by its image in $(\tr(\cH)^{rig})^*$.
\end{corollary}

Combining Theorem \ref{thm:HN} and Corollary \ref{cor:factors}, we get
\begin{proposition}
\label{prop:uniquedual}
The image of any $D\in J(G(F)_c)$ under \eqref{eq:embed}
is determined by $D([\bI w\bI])$ where $w$ runs over a set of minimal length representatives for conjugacy classes in $\tilde{W}$ s.t. $\nu^+(w)$ is central.
\end{proposition}

We remark that by definition $\nu^+(w)$ is central iff $\nu^+(w^m)=m\nu^+(w)$ is central. Note further that $\tr(\cH)^{rig}$ is still graded by $k\in \Omega$ given by $\kappa:\tilde{W}\rightarrow\Omega$. In fact, composing with the Cartan decomposition $G(F)=\sqcup_{w\in\tilde{W}} IwI$ we have a map $G(F)\twoheadrightarrow\Omega$, which is a homomorphism and, abusing notation slightly, the so-called Kottwitz map $\kappa$ on $G(F)$. This map is just
$$\kappa:=\val \circ \det: G(F)\to \Z$$
by identifying $\Omega\cong\Z$. We will be interested in understanding the restrictions of the unipotent orbital integrals $I_\lambda(-)$ and the truncated Steinberg characters $\St_{\lambda,c}$ to $\tr(\cH^{rig})$. Further, we want to understand the truncations of the latter for a fixed $k\in \Omega$.

Imitating \cite[IV 1.]{W2} we define
\begin{definition}
\label{def:w2functions}
Suppose $\lambda\in P(n/e)$.
Let $f^{d,e}_\alpha\in\cH^{\val=\frac{nd}{e}}$ be the characteristic function
of $u^{nd/e}\bP_{\alpha^e}$ where $\bP_{\alpha^e}$ is the standard parahoric associated to $\alpha^e$.
\end{definition}

Here
$$u=\begin{pmatrix}
0 & 0 & \cdots & 0 & \varpi\\
1 & 0 & \cdots & 0 & 0\\
0 & 1 & \ddots & 0 & 0\\
0 & 0 & \ddots & 0 & 0\\
0 & \cdots & 0 & 1 & 0
\end{pmatrix}$$ is the matrix in Eq. \eqref{eq:umatrix} except that we can take $\varpi$ to be any uniformizer for $F$. We have that $u$ normalizes $\bI$, and $u^{n/e}$ normalizes $\bP_{\alpha^e}$.

Note that by Lemma \ref{lem:compactelts}, $f^{d,e}_\alpha$ is supported on the compact-mod-center elements and therefore its image in the cocenter lies in the rigid part. Let us now fix $k\in\Z\cong\Omega$ and only look at the span $[\bI w\bI]$ for $w\in W^{\aff}\times\{k\}$ that are compact, i.e. restrict to $\tr(\cH)^{rig,\val=k}$. 
Let $e=n/\gcd(n,k)$ and $d=ke/n$, so that $k=nd/e$. Note that as $\Omega\subset\tW$ we can also view $k$ as a length zero element in $\tilde{W}$. From Proposition \ref{prop:cptnewton} it is then not hard to see the following:
\begin{proposition}\label{prop:HeNiebasis}
    There is a basis of $\tr(\cH)^{rig,\val=k}$ given by (images in the cocenter of) characteristic functions $[\bI w_i\bI]$ where $w_i=(\sigma_i, k)\in W^{\aff}\rtimes \Omega$ with $\sigma_i\in S_{n/e}$ a set of minimal length representatives for conjugacy classes in $S_{n/e}$. Here we embed (as sets) $$S_{n/e}\hookrightarrow S_n\hookrightarrow W^{\aff}\times \{k\}\subset \tilde{W}$$ where the first inclusion is given by permutation of the first $n/e$ elements. In particular, this space has dimension the number of partitions of $n/e$.
\end{proposition}
\begin{proof}
    By construction, each $w_i$ is compact, of minimal length, and the elements are in distinct conjugacy classes. 
    For each $k \in\Omega$ we would like to know the number of compact conjugacy classes in $W^{\aff}\times\{k\}$. If it is the number of partitions of $n/e$ we are done. 
    
    Note that the map $G(F)\twoheadrightarrow \Omega$ is given by $g\mapsto\val(\det(g))\in\Z\cong\Omega$, and there is a section $\tW \hookrightarrow GL_n(F)$ with image generated by permutation matrices and diagonal matrices with diagonal entries in $\varpi^\Z$, where $\varpi\in F$ is a fixed uniformizer. It's easy to see that an element $w\in\tW$ has $\nu^+(w)$ central iff its image in $GL_n(F)$ under this section is compact (e.g., one can verify both properties by replacing $w$ with some power of $w$ that lives in the lattice part $\Z^n$), which is the case iff all eigenvalues for a $g\in GL_n(F)$ have the same valuation. 
    
    When $\val(\det(g))=k$, that $g$ is compact iff all of its eigenvalues have valuation $k/n=d/e$, where $d/e$ denotes the reduced expression. For compact $g$ in the image of the section $\tW \subset GL_n(F)$, we need each cycle of the permutation to have length divisible by $e=n/\gcd(k,n)$. Conversely, for every partition of $n$ for which all parts are divisible by $e$, we have a unique compact conjugacy class in $W^{\aff}\times\{k\}$ mapping to $k$. Hence the space has dimension equal to the number of partitions of $n/e$.
\end{proof}
\begin{example}
If $n=4$, $e=2$, $d=1$, $k=2$, then $S_{n/e}=S_2$ which is abelian.
The elements $(1,2), (s,2)\in W^{\aff}\rtimes \Omega$ send $(a,b,c,d)\in \R^4$ to $(c+1,d+1,a,b)$ and $(d+1,c+1,a,b)$ respectively.
As elements of $\tW=S_n\ltimes \Z^n$
we have $w_1=((13)(24),(0,0,1,1))$ and $w_2=((1324),(0,0,1,1))$.
\end{example}
\begin{corollary}
    The images of the functions $f_\alpha^{d,e}$ from Proposition \ref{prop:wpair2} for $\alpha\vdash n/e$ also give a basis of $\tr(\cH)^{rig,\val=k}$ indexed by partitions of $n/e$.
\end{corollary}
\begin{proof}
The linear independence is clear by imitating e.g. \cite[Corollaire 4.4.]{W1} again. By Proposition \ref{prop:HeNiebasis} the dimension is the number of partitions of $n/e$.
\end{proof}
By the previous Corollary, for any $w\in W^{\aff}\times\{k\}\subset\tilde{W}$ we can write
\[
[\bI w\bI]=\sum_{\alpha\vdash\frac{n}{e}}c(w,\alpha)f_{\alpha}^{d,e}
\]
for some constants $c(w,\alpha)$, and for chosen $w_i$ as above, the matrix with entries $c(w_i,\alpha)$ is a change-of-basis matrix. By the well-known generators-and-relations description of the Iwahori-Hecke algebra $\cH$, where $[\bI w\bI]$ corresponds to the ''standard basis" $T_w$, together with the Cartan decomposition of $G(F)$, we see that $c(w,\alpha)$ are rational functions in $q$ that depend only on $n$, $k$ and $\alpha$, but {\em not} on the local field $F$. We forego the explicit computation of these rational functions, although it should be an interesting exercise.

By Lemma \ref{lem:compactelts} combined with Definition \ref{def:w2functions} and Proposition \ref{prop:cptnewton} we also have that the image of a distribution $D\in J(G(F)_c)$ in $(\tr(\cH)^{rig,\val=k})^*$ is also determined by $D(f^{d,e}_\alpha)$ where $d, e, k\in \Z$ are as before, and $\alpha\vdash n/e$ is a partition. 

Consider now the unipotent orbital integrals $I_\lambda$ and let $k=0$. Let $f_\alpha=f^{0,1}_\alpha$ be the characteristic function of the standard parahoric $\bP_\alpha$. Pairing with $f_{\alpha}$ for varying $\alpha, \lambda$, Proposition \ref{prop:wpair1} (which is basically computed in \cite[Proposition 4.2]{W1} whose proof is characteristic-independent) shows that we get an invertible matrix, in particular the $I_\lambda$ give a basis of $J(G(\cO))=J(G(F)_c)^{\val=0}$. Combined with Eq. \eqref{eq:embed} and Proposition \ref{prop:HeNiebasis}, we get a new proof of Theorem \ref{thm:uni}. We also have

\begin{theorem}
\label{thm:pseudo-steinberg}
Fix $n$ as before and let $d, e\in \Z$, with $e|n$ be arbitrary. There exist unique elements $\tilde{\St}_{\lambda,c,k}\in (\tr(\cH)^{rig,\val=k})^*$ whose pairing with the family of test functions $f_{\alpha}^{d,e}$ is given by the right-hand-side of Proposition \ref{prop:wpair2},
that is 
$$
\tilde{\St}_{\lambda,c,k}(f_{\alpha}^{d,e})=(-1)^{nd-k}\langle e_\alpha, h_\lambda\rangle
$$
When $\cchar F=0$ the elements $\tilde{\St}_{\lambda,c}=\sum_k \tilde{\St}_{\lambda,c,k}$ coincide with the image of the truncated Steinberg characters $\St_{\lambda,c}$ under \eqref{eq:embed}. Moreover, these $\tilde{\St}_{\lambda,c,k}$ for $\lambda\vdash n/e$ form a basis of $(\tr(\cH)^{rig,\val=k})^*$.
\end{theorem}
\begin{proof}
Theorem \ref{thm:steinberg} combined with Proposition \ref{prop:wpair2} shows the first assertion. The second follows from the invertibility of the square matrix given by the pairings $\tilde{\St}_{\lambda,c,k}(f_\alpha^{d,e})=(-1)^{nd-k}\langle e_\alpha,h_\lambda\rangle$ for varying $\alpha,\lambda$.
\end{proof}
Combining the results for all $k\in\Z$, we get an element $\tilde{\St}_{\lambda,c}\in(\tr(\cH)^{rig})^*$ that serves as the image of this version of $\St_{\lambda,c}$ under \eqref{eq:embed} satisfying Theorems \ref{thm:steinberg} and \ref{thm:steinberg2}. The resulting recursive computation of Shalika germs, using the essentially combinatorial proof of \cite[Lemme V 12]{W2} with $\tilde{\St}_{\lambda,c}$ in place of the truncated Steinberg characters, is carried out in Section \ref{sec:combinatorics}.




\subsection{Parabolic induction}\label{subsec:paraind}
We now recall the parabolic restriction formalism, which reduces the computation of orbital integrals and germs for non-elliptic elements to the elliptic case. Note that the constant term map $f\mapsto f^P$ used to normalize unipotent orbital integrals at the beginning of this section is a special case of the parabolic restriction defined below. The main application is Corollary \ref{cor:mastersymfnmanycomponents}, which expresses the orbital symmetric function of a non-elliptic element in terms of those of its elliptic blocks.

\begin{definition}\label{def:parabolic-restriction-group} Suppose $M\subset P=MN\subset G$ are compatible Levi subgroup and parabolic subgroup defined over $\cO$. For $f\in C_c^{\infty}(G(F))$, we define its parabolic restriction (also called parabolic descent, or constant term) $\Res_M^G(f)\in C_c^{\infty}(M(F))$ as in (\ref{eq:constantterm}):
\[
\Res_M^G(f)(m):=\int_{G(\cO)}\int_{N(F)}f(gmng^{-1})dndg
\]
where the measure is normalized so that $G(\cO)$ and $N(\cO)$ have measure $1$.
\end{definition}

\begin{definition}\label{def:parabolic-restriction-algebra} Let $M,N,P$ be as above and $\fm:=\Lie M$. For $f\in C_c^{\infty}(\fg(F))$, we define its parabolic restriction $\Res_{\fm}^{\fg}(f)\in C_c^{\infty}(\fm(F))$ as in (\ref{eq:constanttermLie}):
\[
\Res_{\fm}^{\fg}(f)(X):=\int_{G(\cO)}\int_{\Lie N(F)}f(\Ad(g)(X+Y))dYdg
\]
where the measure is normalized so that $G(\cO)$ and $\Lie N(\cO)$ have measure $1$.
\end{definition}

Recall that $G(F)_c^{\val=0}\subset G(F)$ is the subset of elements whose eigenvalues all have valuation $0$. As $G=GL_n$, we realize $G(F)_c^{\val=0}$ also as a subset of $\fg(F)$.

\begin{lemma}\label{lem:definitionsagree} If $f\in C_c^{\infty}(G(F)_c^{\val=0})$, then Definition \ref{def:parabolic-restriction-group} and \ref{def:parabolic-restriction-algebra} agree; we have $\Res_M^G(f)=\Res_{\fm}^{\fg}(f)$. 
\end{lemma}

\begin{proof} In both definitions, the resulting $\Res_M^G(f)$ and $\Res_{\fm}^{\fg}(f)$ is evidently supported on $M(F)_c^{\val=0}$. Here one may view $M$ as a product of general linear groups and $M(F)_c^{\val=0}$ is again the subset of elements whose all eigenvalues have valuation $0$. For $m\in M(F)_c^{\val=0}$, $\Ad(m):N(F)\rightarrow N(F)$ preserves the Haar measure on $N(F)$. This shows
\[
\int_{N(F)}f(gmng^{-1})dn=\int_{\fn(F)}f(g(m+n)g^{-1})dn
\]
and thus the two definitions agree.
\end{proof}

\begin{example} One has obviously that $\Res_M^G(\1_{G(\cO)})=\1_{M(\cO)}$ from Definition \ref{def:parabolic-restriction-group}, and hence also $\Res_{\fm}^{\fg}(\1_{G(\cO)})=\1_{M(\cO)}$.
\end{example}

\begin{proposition}
\label{prop:reductiontolevi}
Suppose $\gamma\in M(F)$ is $G$-regular, meaning $\gamma$ is regular when viewed as an element of $G$. Then 
$$I_\gamma^G(f)=\left|\det(\Ad(\gamma)_{\fg/\fm}-\id_{\fg/\fm})\right|^{-1/2}\cdot I_\gamma^M(\Res_M^G(f))$$
\end{proposition}

\begin{proposition}
\label{prop:reductiontolevi-alg}
Suppose $\gamma\in \fm(F)$ is $G$-regular, meaning $\gamma$ is regular when viewed as an element of $\fg$. Then 
$$I_\gamma^G(f)=\left|\det(\ad(\gamma)_{\fg/\fm})\right|^{-1/2}\cdot I_\gamma^M(\Res_{\fm}^{\fg}(f))$$
\end{proposition}

\begin{definition} We define $\Ind_{\fm}^{\fg}:C_c^{\infty}(\fm(F))^*\rightarrow C_c^{\infty}(\fg(F))^*$ to be the adjoint of $\Res_{\fm}^{\fg}$. It is called {\em parabolic induction}. Same for $\Ind_M^G:C_c^{\infty}(M(F))^*\rightarrow C_c^{\infty}(G(F))^*$
\end{definition}

In particular, Proposition \ref{prop:reductiontolevi-alg} effectively says that we have the equality of invariant distributions $$\left|\det(\ad(\gamma)_{\fg/\fm})\right|^{-1/2}\cdot\Ind_M^G I_\gamma^M(-)=I_\gamma^G(-)$$
More generally, an important property is that $\Res_{\fm}^{\fg}$ and $\Res_{M}^{G}$ send $G(F)$-coinvariant to $M(F)$-coinvariant and equivalently $\Ind_{\fm}^{\fg}$ and $\Ind_M^G$ send $M(F)$-invariant to $G(F)$-invariant \cite[Lemma 13.1]{Kottwitz}.

\begin{proposition}\label{induction-unipotent} Suppose $\mathbb{O}$ is a nilpotent orbit of $\fm(F)$ and $\tilde{\mathbb{O}}$ the induced orbit in the sense of Lusztig-Spaltenstein, i.e. $\tilde{\mathbb{O}}$ contains an open dense subset of $\mathbb{O}+\Lie N(F)$. Then $\Ind_{\fm}^{\fg}I_{\mathbb{O}}^M=I_{\tilde{\mathbb{O}}}^G$.
\end{proposition}

\begin{remark}
In terms of partitions (in the sense of Definition \ref{def:partition}), if $M=GL_{n_1}\times\cdots\times GL_{n_r}$, and $\mathbb{O}$ is a unipotent orbit corresponding to a sequence of partitions $$\lambda^{(1)}\vdash n_1,\;\ldots,\; \lambda^{(r)}\vdash n_r$$ the induced orbit is $$\tilde{\mathbb{O}}\leftrightarrow (\lambda_1^{(1)}+\cdots+\lambda^{(r)}_1,\ldots,\lambda_k^{(1)}+\cdots+\lambda_k^{(r)})$$ where $k$ is the length of the longest $\lambda^{(i)}$. For example, when $M=T$, the zero orbit in $T$ induces to the principal one in $GL_n$.

\end{remark}

\begin{corollary}
\label{cor:inducedorbits}
For $\gamma\in\fm(F)$ {\em which is $G$-regular}, we have
\[
\Gamma_{\tilde{\mathbb{O}}}^G(\gamma)=\left|\det(\ad(\gamma)|_{\fg/\fm})\right|^{-1/2}\sum_{\substack{\mathbb{O}\subset \fm(F) \\ \Ind_{\fm}^{\fg}\mathbb{O}=\tilde{\mathbb{O}}}}\Gamma_{\mathbb{O}}^M(\gamma),
\]
where $\Gamma_{\tilde{\mathbb{O}}}^G(\gamma)$ and $\Gamma_{\mathbb{O}}^M(\gamma)$ are defined as in Theorem \ref{thm:shalika-for-Lie-algebra} for $G$ and $M$, respectively, and the sum is over all nilpotent orbits $\mathbb{O}$ of $M$ inducing to $\tilde{\mathbb{O}}$, and zero if there are none.
\end{corollary}

\section{Symmetric functions and combinatorics}
\label{sec:symfns}
In this section, we define the combinatorial objects we use and review the theory of symmetric functions, in order to formulate our results. Apart from Propositions \ref{prop:wpair1}, \ref{prop:wpair2}, this section can be skipped on a first reading and referred back to as necessary.
\subsection{Combinatorics}
We begin with two combinatorial definitions.
\begin{definition}
\label{def:partition}
A {\em partition} of an integer $n>0$, written $\lambda\vdash n$ or $\lambda\in P(n)$ is a nonincreasing sequence of positive integers  
$$\lambda_1\geq \ldots\geq \lambda_k>0,\;\;\sum_i\lambda_i=n$$ and a {\em composition} of $n$, written $\alpha \vDash n$ is an ordered collection $(\alpha_1,\ldots,\alpha_k)$ of positive integers such that $\sum_i \alpha_i=n$. In both cases, we write $\ell(\lambda)=\ell(\alpha)=k$ for the length of the composition or partition and denote by $\lambda^t, \alpha^t$ the conjugate partition (resp. composition). Note that the parts of a composition $\alpha\vDash n$ can be sorted to give a partition $\text{sort}(\alpha)\vdash n$.
\end{definition}
Recall that partitions $\lambda\vdash n$ index conjugacy classes in the symmetric group $S_n$. Similarly, to any composition $\alpha \vDash n$ we can associate the {\em Young subgroup} $S_{\alpha_1}\times\cdots \times S_{\alpha_k}\subseteq S_n$, whose conjugacy class only depends on $\text{sort}(\alpha)$, the partition obtained by sorting $\alpha$.
We will draw the Young/Ferrers diagrams of partitions in French notation. We think of them as lying in $\Z^2_{\geq 0}$ with the first box always at $(0,0)$.
For a box $\square\in \lambda$ with coordinates $(i,j)$ we denote 
\begin{equation}
    \label{eq:armslegs}
    a(\square)=\lambda_i-j-1, \; l(\square)=\lambda^t_j-i-1, \; a'(\square)=i, \; l'(\square)=j
\end{equation}
the arm-, leg-, coarm-, and coleg-lengths of the box. The {\em $q,t$-content} of a box is defined to be $q^{a'(\square)}t^{l'(\square)}$. 
Finally, we have
\begin{definition}
A {\em standard Young tableau} $T$ is a filling of the Ferrers diagram of $\lambda\vdash n$ with the letters $1,\ldots,n$ such that the letters increase in columns and rows. The set of all such tableaux for fixed $\lambda$ is denoted $\text{SYT}(\lambda)$.
\end{definition}
Given a standard Young tableau $T\in \text{SYT}(\lambda)$, we denote the box $\square\in \lambda$ labeled $i$ by $\square_i$. 
We also denote by $z_i$ the $q,t-$content of the box $\square_i$. 

The following definition will become useful in the discussion of positivity properties in Sections \ref{sec:combinatorics} and \ref{sec:applications}.
\begin{definition}
\label{def:rowcolmatrices}
For two compositions $\lambda, \mu$ define $$\MM(\lambda,\mu)$$ to be the set of nonnegative integer matrices (of size $\ell(\lambda)\times \ell(\mu)$) whose rows sum to $\lambda$ and columns sum to $\mu$.
\end{definition}

We will also need the following Lemma in Sections \ref{sec:eha} and \ref{sec:combinatorics}.
\begin{lemma}\label{lem:comptab}
To each composition $\alpha \vDash n$ is associated a unique Young tableau $T(\alpha)$ defined as follows. To each $\alpha_i, 1\leq i\leq k$ we assign the sequence of numbers $\sum_{j=1}^{i-1}\alpha_j+1,\sum_{j=1}^{i-1}\alpha_j+2,\ldots \sum_{j=1}^{i}\alpha_j$ and form a tableau by taking one-row diagrams with these fillings, and then dropping them on top of each other, with the rule that gravity brings boxes as low as possible. In particular, the tableau decomposes as a sequence of horizontal $\alpha_i$-strips.
\end{lemma}
\begin{example}
To the compositions $4=2+2$, $4=1+2+1$ and $4=1+3$ we assign the tableaux 
$$
\Young{
\Carre{3}&\Carre{4}\cr
\Carre{1}&\Carre{2}\cr
}\;\;\;\;\;
\Young{
\Carre{4}&\cr
\Carre{2}&\cr
\Carre{1}&\Carre{3}\cr
} \;\;\;\;\; 
\Young{
\Carre{2}&&\cr
\Carre{1}&\Carre{3}&\Carre{4}\cr
}
$$
\end{example}
\begin{remark}
An alternative description is as follows: $T(\alpha)$ is the unique SYT of shape $\text{sort}(\alpha)$ with descent set $\{\alpha_1,\alpha_1+\alpha_2,\ldots,\alpha_1+\cdots+\alpha_{k-1}\}$. 
\end{remark}

The final combinatorial object we will need, especially throughout Sections \ref{sec:knots} and \ref{sec:combinatorics}, is the notion of (rational) Dyck paths.
\begin{definition}
\label{def:dyck}
Let $m,n,k\geq 1, (m,n)=1$. Then the set
$$\Dyck_{km,kn}$$ will be the collection of lattice paths only with North and East steps inside a $kn\times km$-rectangle in $\Z^2_{\geq 0}$, fitting under the diagonal line of slope $m/n$.

The {\em area} of a Dyck path $\pi\in \Dyck_{km,kn}$ is defined to be the number of full squares between the path and the diagonal. Similarly, the {\em coarea} of $\pi$ is defined as $(m-1)(n-1)/2-\area(\pi)$, which is the number of squares below the path. 
\end{definition}

\subsection{Symmetric functions}\label{subsec:Sym}
Let $\Sym_{q,t}$ be the ring of symmetric functions over $\Q(q,t)$ in the alphabet $\{X_1,\ldots,X_n,\ldots\}$ and denote the five usual bases of monomial, homogeneous, elementary, Schur, and power sum symmetric functions by
$$\{m_\lambda\},\{h_\lambda\},\{e_\lambda\},\{s_\lambda\},\{p_\lambda\}$$
Here $\lambda$ runs over all integer partitions in the sense of Definition \ref{def:partition}. 

In addition to the above bases, the symmetric functions of Macdonald will be important for us. For the reader who wants more details, the theory is very well covered in many sources. See for example \cite[Section 3]{HaimanCDM}.

Recall that the {\em modified Macdonald polynomials} $\tH_\lambda[X;q,t], \lambda\vdash n$ are the unique symmetric functions with the properties
\begin{eqnarray}
\tH_\mu[X(1-q);q,t]&\in&\Q(q,t)\{s_\lambda|\;\lambda\geq \mu\}\\
\tH_\mu[X(1-t);q,t]&\in&\Q(q,t)\{s_\lambda|\;\lambda\geq \mu^t\}\\
\langle \tH_\mu[X;q,t],s_{(n)}\rangle&=&1
\end{eqnarray}
Here the last pairing is the Hall inner product, defined in Definition \ref{def:innerproducts}.

We do not require much of the advanced theory of Macdonald polynomials, but let us note down the following definition as well as some specializations.
\begin{definition}
\label{def:nablaoperator} 
The operator $\nabla$ of Bergeron and Garsia scales by definition each $\tH_\lambda$ by $q^{n(\lambda)}t^{n(\lambda^t)}$ where $n(\lambda)=\sum_{i=1}^{\ell(\lambda)}(i-1)\lambda_i$.
\end{definition}
From \cite[Proposition 3.5.8.]{HaimanCDM} we have
\begin{lemma}[The limit of $\tH_\lambda$ as $q\to 1$]	
\label{lem:plethhomog}
The modified Macdonald symmetric function $\tH_\lambda$ at $q=1$ is given by
\begin{equation}
\label{eq:limitmacdonald}\tH_\lambda[X;1,t]=(1-t)^{|\lambda|}[\lambda^t]_t!h_{\lambda^t}[X/(1-t)]=:\th_{\lambda^t}[X;t]
\end{equation}
in other words a plethystically transformed homogeneous symmetric function, up to normalization. 
\end{lemma}
From the $q,t$-symmetry $\tH_\lambda[X;t,q]=\tH_{\lambda^t}[X;q,t]$ we immediately have
\begin{corollary}[The limit as $t\to 1$]\label{cor:tHlambda}
$$\tH_\lambda[X;q,1]=\th_\lambda[X;q]$$
\end{corollary}
In addition to the $q,t$-symmetry, we have the symmetry under inverting $q$ and $t$ \cite[Proposition 3.5.12.]{HaimanCDM}
\begin{equation}
\label{eq:invertqt}
t^{-n(\mu)}q^{-n(\mu^t)}\omega \tH_\mu[X;q,t]=\tH_\mu[X;q^{-1},t^{-1}]
\end{equation}
so that in particular $\lim_{t\to 1}\tH_\mu[X;q^{-1},t^{-1}]=q^{n(\mu^t)}\omega\th_\mu[X;q^{-1}]$. Note that this implies 
\begin{equation}
\label{eq:invertq}
\omega\th_\mu[X;q]=q^{n(\mu^t)}\th_\mu[X;q^{-1}]
\end{equation}

We will denote $\th_\lambda:=\tH_\lambda[X;q,1]$ and call these the {\em specialized Macdonald symmetric functions} or the {\em plethystically transformed homogeneous symmetric functions}.  Later on, we will also need the prefactor
\begin{equation}
\label{eq:ccoeff}
c_\lambda(q):=(1-q)^{|\lambda|}[\lambda]_q!=\prod_{i=1}^{\ell(\lambda)}\prod_{j=1}^{\lambda_i}(1-q^j)
\end{equation}
in Eq. \eqref{eq:limitmacdonald}. Note that 
$$\th_\lambda=c_\lambda^{Wal}(q)e_\lambda\left[\frac{X}{q-1}\right]=(-1)^n c_\lambda^{Wal}\omega e_\lambda\left[\frac{X}{1-q}\right],$$ where 
\begin{equation}
\label{eq:ccoeffwal}
c_\lambda^{Wal}:=(-1)^nc_\lambda=\prod_{i=1}^{\ell(\lambda)}\prod_{j=1}^{\lambda_i}(q^j-1)
\end{equation}
is exactly the prefactor defined on \cite[pp. 201]{W1}.

Note that compared to the Macdonald polynomials, $\th_\lambda$ are much simpler in behaviour. For instance, they are multiplicative: $$\th_\lambda\th_\mu=\th_{\lambda+\mu}$$ and one can deduce combinatorial expansions for them in terms of the other standard bases via known relations between $h_\lambda$ and these bases. This will turn out to be important in the proof of Theorem \ref{thm:transitionmatrix}. Further, the $\nabla$-operator from Definition \ref{def:nablaoperator} becomes a ring homomorphism on symmetric functions in this limit.

We will also need a few different inner products on the ring of symmetric functions, the interplay of whom turns out to play a key role. We remark that by ''inner product" we simply mean a symmetric bilinear form
valued in $\Q(q,t)$.
\begin{definition}
\label{def:innerproducts}
\begin{enumerate}
\item The {\em Hall inner product} is the inner product on $\Sym_{q,t}$ defined by
$$\langle s_\lambda,s_\mu\rangle=\delta_{\lambda\mu}$$
\item The {\em $q$-inner product} is $$\langle f,g\rangle_q:=\langle f,g\left[\frac{X}{(1-q)}\right]\rangle$$
\item The {\em $q,t$-inner product} is $$\langle f,g\rangle_{q,t}:=\langle f,g\left[\frac{(1-q)}{(1-t)}X\right]\rangle$$
\item The {\em geometric inner product} is $$(f,g)=-q^{\deg f}\left\langle (\nabla^{-1}(f))[X(1-t^{-1})],g[X(1-t^{-1})] \right\rangle_{q,t^{-1}}$$
\end{enumerate}
\end{definition}
\begin{definition}[Frobenius characteristic]\label{def:frobenius}
The {\em Frobenius characteristic map} is the isomorphism $$\ch: \bigoplus_{n\geq 0}\text{Rep}(S_n)\xrightarrow{\sim} \Sym_\Z$$ sending the irreducible representation $V_\lambda$ of $S_n$ to the Schur function $s_\lambda$. It is an isometry with respect to the natural inner product on characters and the Hall inner product.
\end{definition}
%
\begin{remark}\label{rmk:glambda}
The geometric inner product is the one naturally arising from the geometry of Hilbert schemes of points. It is characterized as the unique inner product satisfying $(\tH_\lambda,\tH_\mu)=\delta_{\lambda\mu}\,g_\lambda$, where
\begin{equation}
    \label{eq:glambda}
    g_\lambda=\prod_{\square\in\lambda}\bigl(1-q^{a(\square)}t^{-l(\square)-1}\bigr)\bigl(1-q^{-a(\square)-1}t^{l(\square)}\bigr),
\end{equation}
 In particular the order of vanishing of $g_\lambda$ at $t=1$ is $\operatorname{ord}_{t=1}g_\lambda=\#\{\square\in\lambda: a(\square)=0\}=\ell(\lambda)$: the first factor $1-q^{a}t^{-l-1}$ vanishes at $t=1$ exactly when $a(\square)=0$ (one box per row, so $\ell(\lambda)$ of them), while the second never vanishes there.
\end{remark}

\subsection{Combinatorial evaluation of some distributions}
The first appearance of symmetric functions in connection to harmonic analysis on $G(F)$ are the following formulas due to Howe and Waldspurger, which give a useful way of recasting values of truncated Steinberg distributions and nilpotent orbital integrals in terms of Hall pairings.

\begin{proposition}
\label{prop:wpair1}
Let $\1_\lambda$ be the characteristic function of the standard parahoric subgroup $\bP_\lambda\subset GL_n(F)$
corresponding to the partition $\lambda$, {\em divided by the measure of $\bP_\lambda$}. Let $I_{\mu^t}(-)$ be the orbital integral over the unipotent class of type $\mu^t$.
Then we have 
\begin{equation}
    \label{eq:wpair1}
    I_{\mu^t}(\1_\lambda)=\langle h_\lambda,\th_\mu\rangle
\end{equation}
where $\langle -, -\rangle$ is usual  Hall inner product from Definition \ref{def:innerproducts}. 

\end{proposition}
\begin{proof}
    Let $\psi_\lambda$ be the characteristic function of $\bP_\lambda$, i.e. $\psi_\lambda=c_\lambda(q)c_n(q)^{-1} \1_\lambda$\footnote{Warning: Waldspurger uses the notation $\varphi_\lambda$ for our $\psi_\lambda$. For him, $\psi_\lambda$ denotes a different function. We reserve the notation $\varphi$ for the plethysms introduced in the next sections.}. By \cite[Prop. 4.2.]{W1} we get that in the normalizations we have chosen, 
    $$I_{\mu^t}(\psi_\lambda)=
    c_\lambda^{Wal}(q)c_{\mu}^{Wal}(q)c_n^{Wal}(q)^{-1}\langle x_\lambda, y_\mu\rangle_{Zel}
    $$ where $x_\lambda, y_\mu$ in the notations of \cite{W1, Zelevinsky} correspond to our $h_\lambda$ and $e_\mu$ by \cite[10.2.]{Zelevinsky}.

    Since $e_\mu\left[\frac{X}{q-1}\right]=(-1)^n h_\mu\left[\frac{X}{1-q}\right]$ we get $$
    I_{\mu^t}(\psi_\lambda)=(-1)^{4n} c_\lambda c_\mu c_n^{-1}\langle h_\lambda, h_\mu\left[\frac{X}{1-q}\right]\rangle
    $$ 
Absorbing $c_\mu$ to $\th_\mu$ and taking into account the difference between $\psi_\lambda$ and $\1_\lambda$ we get 
$$I_{\mu^t}(\1_\lambda)=c_\lambda^{-1} c_nI_{\mu^t}(\psi_\lambda)=\langle h_\lambda, \th_\mu\rangle$$
\end{proof}

A similar result holds for the distributions $\tilde{\St}_{\lambda,c,k}$ from Theorem \ref{thm:pseudo-steinberg} (or equivalently, the truncated Steinberg characters $\St_{\lambda,c,k}$ from Theorems \ref{thm:steinberg}, \ref{thm:steinberg2} in characteristic zero).
\begin{proposition}[Lemme V 11. \cite{W2}]
\label{prop:wpair2}

Suppose $e\geq 1, (d,e)=1$ and $\mu\in P(n/e)$. Denote also $k=nd/e$, and
let $f^{d,e}_\mu$ be the characteristic function
of $u^{nd/e}\bP_{\mu^e}$ where $\bP_{\mu^e}$ is the parahoric associated to $\mu^e$, again divided by the measure of the subset as in the previous Proposition. Note
Here $u$ is the matrix in Eq. \eqref{eq:umatrix}, or any suitable lift of the generator of $\Omega\subset \tW$ to the group $G(F)$. Note that $f^{0,e}_{\mu}=\1_{\mu}$.

Then 
$$\tilde{\St}_{e\lambda, c, k}(f^{d,e}_{\mu})=(-1)^{nd-nd/e}\langle h_\mu, e_\lambda\rangle$$ where the pairing is the Hall inner product and $\tilde{\St}_{e\lambda,c,k}(-)$ is the distribution from Theorem \ref{thm:pseudo-steinberg}.

In particular, when $d=0$, we have $\tilde{\St}_{e\lambda,c,k}(f^{0,e}_{\mu})=\langle h_\mu, e_\lambda\rangle$.
\end{proposition}
\begin{proof}
    In characteristic zero, $\tilde{\St}_{e\lambda,c,k}=\St_{e\lambda,c,k}$ and this is a direct translation of \cite[V. II]{W2}, with the identifications $x_\mu=h_\mu, y_\lambda=e_\lambda$ as above. The general case then follows since the formula is purely combinatorial and independent of the characteristic of $F$.
\end{proof}

\section{The elliptic Hall algebra and annular knot invariants}
\label{sec:eha}
In this section, we define the elliptic Hall algebra and recall its action on the ring of symmetric functions. We then use this action to define the refined annular symmetric functions $\fqrs$. The $t\to 1$ degeneration of this action, which is the key to our main results, is then studied and used to define the annular symmetric functions $\fpq$.

\subsection{The elliptic Hall algebra}\label{subsec:eha}

Throughout this section, we work with three parameters $q_1, q_2, q_3$ subject to $q_1 q_2 q_3 = 1$. The dictionary between these and the Macdonald parameters $(q,t)$ used elsewhere in the paper is
\begin{equation}\label{eq:dictionary}
q_1 = q, \qquad q_3 = t^{-1}, \qquad q_2 = (q_1 q_3)^{-1} = q^{-1} t.
\end{equation}

\begin{definition}\label{def:eha}
The {\em elliptic Hall algebra} is the $\Q(q_1, q_2)$-algebra $\cE = \cE_{q_1, q_2}=\cE_{q_1,q_2,q_3}$, where $q_3 := (q_1 q_2)^{-1}$, generated by elements $P_{m,n}$, $(m,n) \in \Z^2 \setminus (0,0)$ subject to
$$[P_{m_1,n_1}, P_{m_2,n_2}] = 0$$ if $(m_1,n_1)$, $(m_2,n_2)$ lie on the same line through the origin, and 
$$[P_{m_1,n_1}, P_{m_2,n_2}] = \frac{\theta_{m_1+m_2, n_1+n_2}}{\alpha_1}$$ if $(m_1,n_1), (m_2,n_2), (m_1+m_2, n_1+n_2)$ form a quasi-empty triangle. Here, for $(m,n)$ coprime and $k>0$,
$$\exp\left(\sum_{k=1}^\infty P_{km, kn} \alpha_k x^k\right) = \sum_{\ell=1}^\infty \theta_{\ell m, \ell n} x^\ell, \qquad \alpha_k = \frac{(q_1^k - 1)(q_2^k - 1)(q_3^k - 1)}{k}.$$
\end{definition}

\begin{proposition}[Triangular decomposition]
Let $\cE^>$, $\cE^<$, $\cE^0$ be the subalgebras generated by $\{P_{m,1}\}_{m \in \Z}$, $\{P_{m,-1}\}_{m \in \Z}$, and $\{P_{\pm k, 0}\}_{k > 0}$, respectively. The multiplication map gives a $\Q(q_1, q_2)$-linear isomorphism $\cE^< \otimes \cE^0 \otimes \cE^> \to \cE$.
\end{proposition}

For the rest of this paper, we restrict our attention to the nonnegative part $\cE^\geq \subset \cE$, which is by definition generated by $P_{m,1}, m\in \Z$ and $P_{k,0}, k > 0$. 
Over the field $\Q(q_1,q_2)$, the positive part $\cE^>$ is isomorphic to the shuffle algebra of \cite{FHHSY} by \cite{NegutShuffle}; we use this identification freely.

\begin{remark}[The integral form]\label{rmk:integralform}
One cannot set $q_3\to 1$ directly in Definition \ref{def:eha}, since the structure constant $\alpha_k$ contains the factor $(q_3^k-1)/k$, which vanishes there. This is remedied by finding an integral form $\cE^>_\Z\subset\cE^>$. Such a subalgebra was explicitly described by \cite{NegutIF}. Explicitly, $\cE^>_\Z$ is a $\Z[q_1^\pm,q_2^\pm]$-subalgebra generated by rescaled elements $\bar P_{km,kn}=\zeta_{km,kn}P_{km,kn}$, where the factor $$\zeta_{km,kn}:=\frac{q_2^k-1}{(q_1^k-1)(q_3^k-1)}(q_1-1)^{kn}(q_3-1)^{kn}$$ is chosen to clear the denominators that $P_{km,kn}$ carries in the shuffle presentation of \cite[(45)]{NegutIF}. We do {\em not} use this construction. The $t=1$ action built below is a representation-level phenomenon, whose regularity at $t=1$ comes from the geometric normalization of the operators $P_{km,kn}$ in the modified Macdonald basis. The details are explained in Theorem \ref{thm:t1specialization}.
\end{remark}

We will also need the following subalgebras of the above.
\begin{definition}
\label{def:slopesubalg}
Let $m,n\in \Z^2_{\geq 0}$, $(m,n)=1$. The {\em slope $\frac{m}{n}$-subalgebra} of $\cE^{\geq}$ is the subalgebra $\cE^{m/n}$ generated by $P_{km,kn}, k\geq 0$.
\end{definition}
\begin{theorem}[\cite{NegutShuffle}]
\label{thm:slopesubalg}
Let $\Sym_{q,t}$ be the algebra of symmetric functions over $\Q(q,t)$ as introduced in Section \ref{sec:symfns}. There is an algebra isomorphism
$$\widehat{\varphi}_{m/n}: \Sym_{q,t}\to \cE^{m/n}$$ sending $p_k\mapsto P_{km,kn}$.
\end{theorem}

\subsection{The action on symmetric functions}\label{subsec:fock}
We now study the action of $\cE^\geq$ on symmetric functions. 
The basic fact is that by \cite[Theorem 0.2]{SVKTheory} and \cite[Theorem 3.5]{FT}, we have
\begin{theorem}
\label{thm:fockrep}
There is an action of $\cE_{q_1,q_2,q_3}$ on $\Sym_{q,t}$ by shuffle algebra operations, uniquely determined by the requirements that the operators $P_{m,0}$ for $m \geq 1$ act as multiplication by the power sum symmetric functions $p_{m}$, that $P_{0,1}$ acts as a Macdonald eigenoperator.
\end{theorem}
\begin{remark}\label{rmk:asymmetry}
The defining relations of $\cE$ are symmetric in the $q_i$, but the construction of the representation in Theorem \ref{thm:fockrep} is not. This is the source of the asymmetric role played by $t$ relative to $q$ in the rest of the paper. 
\end{remark}

As in Section \ref{sec:symfns}, the space $\Sym_{q,t}$ has a $\Q(q,t)$-basis $\tH_\lambda$, and it is enough to describe the matrix coefficients of $P_{km,kn}$ in this basis to completely determine this representation. By results in \cite{NegutKTheory} and \cite{FT}, these matrix coefficients are known.

Following conventions in the physics literature around the EHA, we denote the Macdonald basis by the bra-ket notation $\tH_\lambda=|\lambda\rangle$, together with the inner product $\langle \lambda|\mu\rangle=\delta_{\lambda\mu}g_\lambda$, which corresponds to the geometric inner product $(-,-)$ on symmetric functions as in Definition \ref{def:innerproducts}. With this notation at hand, the matrix coefficients can be expressed as follows. Throughout, $(m,n) \in \Z_{\geq 0}^2$ are coprime and $k \geq 1$.


The action of Theorem \ref{thm:fockrep} is the action of \cite{SVKTheory, FT} on $K_T(\Hilb) \cong \Sym_{q,t}$: under this isomorphism the torus-fixed-point classes $I_\lambda$ correspond to the modified Macdonald polynomials $\tH_\lambda$ \cite[Thm.~4.2]{GN}, and $P_{m,0}$ acts as multiplication by $p_m$ \cite[Rmk.~4.3]{GN}. It is in this fixed-point basis $\{\tH_\lambda\}$ that the matrix coefficients of $P_{km,kn}$ take the closed form of Theorem \ref{thm:matrixcoeff}.

\begin{theorem}[\protect{\cite[Eq.~(37)]{GN}, \cite{NegutKTheory}}]
\label{thm:matrixcoeff}
We have\footnote{This is \cite[Eq.~(37)]{GN} in our conventions: their operator $\mathcal{P}_{kn,km}$ is our $P_{km,kn}$ (see above), and we have switched $\mu$ and $\lambda$.}
$$
\langle \lambda|P_{km,kn}|\mu\rangle=\frac {\xi^{kn}}{[k]} \cdot \frac{g_\mu}{g_\lambda} 
\sum^{\SYT}_{\lambda = \mu + \square_1 + \ldots + \square_{kn}} \quad \left[\sum_{j=0}^{k-1}  (qt)^{j} \frac {z_{n(k-1)+1}z_{n(k-2)+1}\cdots z_{n(k-j)+1}}{z_{n(k-1)}z_{n(k-2)}\cdots z_{n(k-j)}} \right] \cdot 
$$
$$\cdot \frac {\prod_{i=1}^{kn} z_i^{S'_{m/n}(i)} (qt z_i - 1 ) }{\left(1 - qt\frac {z_2}{z_1}\right)\cdots \left(1 - qt\frac {z_{kn}}{z_{kn-1}} \right)} \prod_{1\leq i<j \leq kn} \omega'~^{-1} \left(\frac {z_j}{z_i} \right)  \prod_{1\leq i \leq kn}^{\square \in \mu} \omega'~^{-1} \left(\frac {z(\square)}{z_i} \right) $$
\end{theorem}
\noindent
where $$\omega'\left(x\right)=\frac{(x-1)(x-qt)}{(x-q)(x-t)},\; \xi=\frac{(q-1)(t-1)}{qt(qt-1)},\; [k]=\frac{(q^k-1)(t^k-1)}{(q-1)(t-1)}$$ and $$S_{m/n}(i)':=\lfloor \frac{im}{n}\rfloor -\lfloor \frac{(i-1)m}{n}\rfloor.$$

Here the outer sum is over standard Young tableaux on the skew shape $\lambda/\mu$, that is, chains $\mu=\lambda^{(0)}\subset\lambda^{(1)}\subset\cdots\subset\lambda^{(kn)}=\lambda$ adding one box $\square_i$ at each step; $z_i:=z(\square_i)=q^{a'(\square_i)}t^{l'(\square_i)}$ as before, and $z(\square)$ for $\square\in\mu$ is the content of a fixed box of $\mu$.

We now show that the action of Theorem \ref{thm:fockrep} has a well-defined specialization at $t=1$, on $\Sym_q$ in the basis $\{\th_\lambda := \tH_\lambda[X;q,1]\}$. This is a property of the representation on $\Sym_{q,t}$: it follows from the regularity at $t=1$ of the matrix coefficients of Theorem \ref{thm:matrixcoeff}, and uses no integral form of the algebra (cf.\ Remark \ref{rmk:integralform}).

\begin{theorem}\label{thm:t1specialization}
Let $R = \C(q)[t]_{(t-1)}$ and $V_R := \bigoplus_\lambda R\,\tH_\lambda \subset \Sym_{q,t}$. The matrix coefficients $\langle \lambda \mid P_{km,kn} \mid \mu \rangle$ of Theorem \ref{thm:matrixcoeff}, in the modified Macdonald basis, are regular at $t=1$ (equivalently $q_3=1$) for all $\lambda, \mu, k, m, n$. Consequently each operator $P_{km,kn}$ preserves the lattice $V_R$, and the geometric action of Theorem \ref{thm:fockrep} specializes at $t=1$ to a well-defined $\C(q)$-linear action on
$$V_R \otimes_R \C(q) \cong \Sym_q$$
with the specialized operators denoted $P_{km,kn}|_{t=1}$. 
\end{theorem}
The proof is given in Subsection \ref{subsec:weightformulas} below, after the necessary order of vanishing analysis has been carried out in Proposition \ref{prop:nonvanishingcoeffs} and Lemma \ref{lem:wt-nonvanishing}.

By Theorem \ref{thm:t1specialization}, the "geometric" action of $\cE^\geq$ on $\Sym_{q,t}$ specializes at $t=1$ to a well-defined action on $\Sym_q$ in the basis $\{\th_\lambda\}$. We now identify this $t=1$ action explicitly.

\begin{proposition}\label{prop:multops}
The $t=1$ specialization of the action (Theorem \ref{thm:t1specialization}) is by multiplication operators: each operator $P_{km,kn}|_{t=1}$, or equivalenly each $H_{m,n}|_{t=1}$, 
where $$1+\sum_{k=1}^\infty\frac{H_{km,kn}}{x^k}=\exp\left(\sum_{k=1}^\infty\frac{P_{km,kn}}{kx^k}\right),$$
acts on $f \in \Sym_q$ as multiplication by an explicit symmetric function.
\end{proposition}
\begin{proof}
As explained in \cite{BS1,NegutShuffle}, $\Q(q,t)$, $\cE^\geq$ is generated by the elements $H_{m,n}$, $(m,n)\in\Z\times\N$, together with the $P_{k,0}$, $k>0$. The $P_{k,0}$ act as multiplication by $p_k$ by construction, so it remains to show that each $H_{m,n}|_{t=1}$ acts by multiplication. By \cite[Theorem 2.15]{NegutOps}, the operator $H_{m,n}$ acts on $f\in\Sym_{q,t}$ by the contour integral
$$H_{m,n}\cdot f[X]=\oint\cdots\oint\frac{\prod_{i=1}^n z_i^{S_{m/n}'(i)}}{\prod_{i=1}^{n-1}(1-qt\frac{z_{i+1}}{z_i})\prod_{i<j}\omega'(z_j/z_i)}\cdot$$
$$\cdot\bigwedge\nolimits^\bullet\left(-\frac{X}{z_1}\right)\cdots \bigwedge\nolimits^\bullet \left(-\frac{X}{z_n}\right)\cdot
f\left[X-(1-q)(1-t)\sum_{i=1}^n z_i\right]\prod_{a=1}^n\frac{dz_a}{2\pi iz_a}$$
where $\bigwedge^\bullet(-\frac{X}{z})=\sum_{k=0}^\infty \frac{h_k[X]}{z^k}$ and the contours are concentric circles $|z_n|<\cdots<|z_1|$ separating $0$ from $\infty$ (see \cite{NegutOps, NegutKTheory} for details). At $t=1$, the plethystic shift $f[X-(1-q)(1-t)\sum_i z_i]$ becomes the identity, so $f$ factors out of the integral, and the remaining integral over $z_1,\ldots,z_n$ produces a symmetric function of $X$ alone. Therefore $H_{m,n}|_{t=1}$ acts as multiplication by this symmetric function.
\end{proof}
\begin{remark}
This Proposition is conjectured in \cite{BergeronNotes,BergeronSurvey}. 
\end{remark}

By Proposition \ref{prop:multops}, the slope $m/n$-plethysm $\widehat{\varphi}_{m/n}$ from Theorem \ref{thm:slopesubalg} applied to elements of $\Sym_q$ produces explicit symmetric functions at $t=1$. The case of elementary symmetric functions is particularly important for us, and completely determines the behavior of $\widehat{\varphi}_{m/n}$ at $t=1$.
\begin{proposition}
\label{prop:Emnk}
Suppose that $m,n>0$ and $\gcd(m,n)=1$. At $t=1$ the operator $\widehat{\varphi}_{m/n}(e_k)|_{t=1}$ becomes a multiplication operator by the degree $kn$ symmetric function:
$$
E_{m,n,k}:=\sum_{\pi\in \Dyck_{kn,km}}q^{\area(\pi)}e_{\pi}.
$$
Here $\pi$ is a rational Dyck path in a $(kn\times km)$ rectangle below the line of slope $m/n$, $\area(\pi)$ is the area between $\pi$ and the diagonal, and $e_\pi:=\prod_{i} e_{h_i(\pi)}$ where $h_i(\pi)$ are the horizontal run lengths of $\pi$ (i.e., the lengths of the maximal East runs between consecutive North steps, which sum to $kn$).
\end{proposition}
\begin{proof}
Given Proposition \ref{prop:multops}, this is \cite[Eq. (4.5.4)]{BergeronSurvey} (see also \cite{BGLX}). 
\end{proof}
This motivates the following.
\begin{definition}
\label{def:mnpleth}
The {\em slope $m/n$-plethysm in the $t=1$ limit} is the homomorphism
$\varphi_{m/n}: \Sym_{q}\to \Sym_{q}$
 defined by $\varphi_{m/n}(e_k)=E_{m,n,k}$. 
\end{definition}
As it will be always clear from the notation and context, we will refer to this limiting version $\varphi_{m/n}$ as the slope $m/n$-plethysm as well.
For any $\lambda\vdash n$, we will also use the notation
    $$E_{m,n,\lambda}:=\varphi_{m/n}(e_\lambda)=\prod_{j=1}^{\ell(\lambda)} E_{m,n,\lambda_j}.$$


\subsection{Superpolynomials and the annular symmetric function}
\label{sec:knots}
We now explain how the $\cE^{\geq 0}$-action on $\Sym_{q,t}$ can be used to produce invariants of algebraic knots.
We give the definition first and then discuss its relation to existing constructions.

Recall from Definition \ref{def:puiseux} and Remark \ref{rmk:puiseux} that to any irreducible germ of a plane curve $\{f=0\}\subset \C^2$  we may associate both a sequence of Newton pairs and the isotopy class of a link $L=\Link_0(f)\subset S^3$. 
It is well known that these data are essentially equivalent. We now explain the correspondence, and refer to \cite{EN} for a more thorough account.

\begin{figure}
\begin{center}
    \includegraphics[scale=0.3]{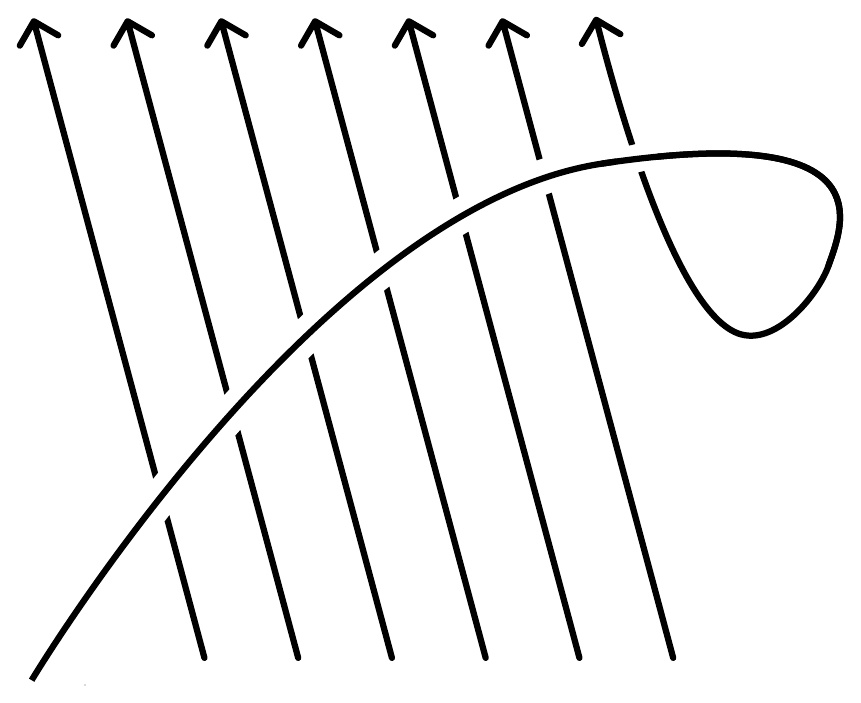}
    \caption{The Coxeter braid $\cox_7$.}
    \label{fig:coxeter}
\end{center}
\end{figure}

When $f$ is quasi-homogeneous and there is a single Newton pair $(p,q)$, we have that $L$ equals the torus knot $T(p,q)$. It is the braid closure in $S^3$ of $\cox_p^q$ , where $\cox_p$ is the Coxeter braid in Fig. \ref{fig:coxeter}). In fact, after fixing a coordinate system for $\C^2$, the sphere $S^3$ decomposes into two solid tori, with the knot $T(p,q)$ lying in a single one. We denote this knot in the solid torus by $T^q_p$. In what follows, it is best thought of as the annular closure of the diagram of $\cox_p^q$ shown in Fig. \ref{fig:coxeter} (in the blackboard framing).

When $f$ has several Newton pairs $(\vec{p},\vec{q})$, we define 
\begin{equation}
\label{eq:satellite}
T(\vec{p},\vec{q}):=T_{p_k}^{q_k}*(T_{p_{k-1}}^{q_{k-1}}*(\cdots *(T_{p_1}^{q_1})\cdots),
\end{equation}
where the right-hand side is constructed as follows.
\begin{definition}
For arbitrary knots $L_1, L_2$ in the solid torus, or more generally elements in the skein algebra of the annulus, we define the {\em satellite} of $L_1$ by $L_2$, denoted $L_1*L_2$ by thickening $L_1$ to an annulus and placing the diagram of $L_2$ inside this annulus. 
\end{definition}
\begin{remark}
\label{rmk:iteratedknots}
The sequence
$(\vec{p},\vec{q})$ is denoted $(\vec{\mathbf{r}},\vec{\mathbf{s}})$ in \cite{CD}.  
In fact, it can be any sequence of (not necessarily positive) coprime integers, in which generality we obtain {\em iterated torus knots}.
However, the Newton pairs of a plane curve germ are always positive and always have $p_k>1$.
\end{remark}

The upshot of the above discussion is that by choosing the convention $q_1>p_1$ for the first Newton pair (essentially, picking which solid torus $L$ lies inside), the isotopy type of the link $L$ equals $T(\vec{p},\vec{q})$, which is uniquely determined by $(\vec{p}, \vec{q})$ and vice versa. In light of the discussion in the introduction, it is also worth noting that this means one can canonically attach a link $L:=T(\vec{p},\vec{q})$ to any totally tamely ramified element $\gamma\in \fg\fl_n(F)$, where $F$ is local non-archimedean and $(\vec{p},\vec{q})$ is the Newton sequence of $\gamma$.

For general iterated torus knots as in Remark \ref{rmk:iteratedknots}, Morton--Samuelson \cite{MS} gave a definition of knot invariants called {\em superpolynomials} $\mathbf{P}_L(q,t,a)$, using the $\cE^\geq$-action on $\Sym_{q,t}$. We review their definition in Definition \ref{def:superpolynomial}. The construction is closely related to an earlier DAHA-based approach by Cherednik--Danilenko \cite{CD}, and appropriately interpreted, the two definitions coincide. Finally, almost simultaneously with these works, Gorsky--Neguț \cite{GN} gave an explicit combinatorial formula for $\mathbf{P}_L(q,t,a)$ in the case of torus knots. In fact, it is a special case of  the computation of matrix coefficients of $P_{m,n}$ in the Macdonald basis from the previous subsection.

\begin{remark}
The name ``superpolynomial'' was introduced by Dunfield--Gukov--Rasmussen in \cite{DGR}, where the existence of a knot polynomial satisfying certain properties was hypothesized. In particular, it was suggested that there should be a triply graded knot homology theory whose Poincar\'e polynomial satisfies said properties. Around the same time, Khovanov and Rozansky \cite{Khovanov} defined such a theory with the expected properties. Letting
\begin{equation}
\label{eq:superpoly}
P_L(q,t,a):=\sum_{i,j,k}q^it^ja^k\dim \HHH^{i,j,k}(L)
\end{equation}
be the Poincar\'e polynomial of this {\em triply graded Khovanov-Rozansky} (or HOMFLY) homology $\HHH(L)$ gives an alternative definition of the knot superpolynomial. The representation-theoretic superpolynomials $\mathbf{P}_L(q,t,a)$ from Definition \ref{def:superpolynomial} and the Khovanov-Rozansky superpolynomials of $P_L(q,t,a)$ in Eq. \eqref{eq:superpoly} are known to agree up to a change of variables for torus knots \cite{Mel,EH,HM} and conjectured to agree in general, but this has not been proved at the time of writing.

We also remark that the homology theory $\HHH$, and in particular the rational function $P_L(q,t,a)$, has been compared with Hilbert schemes of points on plane curve singularities and their compactified Jacobians in \cite{ORS}. On the other hand, we compare $\mathbf{P}_L(q,t,a)$ to the compactified Jacobian point-count in Proposition \ref{prop:superpolycptjac}, proving an expectation of Cherednik--Danilenko \cite{CD}. This gives significant evidence for the main conjecture of \cite{ORS}. We refer the reader to \cite{GKR} for more discussion about the construction of $\HHH$ and the precise relationship to plane curve singularities. 
\end{remark}




Most importantly for us, each of the approaches in \cite{GN,CD,MS} contains an intermediate symmetric function $\fqrs$, which is only an invariant of the underlying annular link. We call this symmetric function the {\em refined annular symmetric function}; the precise definition will be given below. In broad strokes, the superpolynomial itself is related to spherical orbital integrals and knots in $S^3$, whereas the (refined) annular symmetric function knows about germ expansions and knots in the solid torus as well.  

As above, for a sequence of pairs of coprime integers $(p_1,q_1),\ldots, (p_k,q_k)$ we have an iterated torus knot $T(\vec{p},\vec{q})=T_{q_k}^{p_k}*\cdots*T_{q_1}^{p_1}$ as in \eqref{eq:satellite}. By Theorem \ref{thm:slopesubalg}, we have algebra homomorphisms $\widehat{\varphi}_{q_i/p_i}: \Sym_{q,t} \to \cE^{q_i/p_i}$ for $i=1,\ldots,k$ sending $p_\ell \mapsto P_{\ell q_i, \ell p_i}$. 
The following is essentially the same as \cite[Definition 7.15]{MS}.

\begin{definition}
\label{def:fullmastersymfn}
The {\em refined annular symmetric function} associated to $(\vec{p},\vec{q})$ is 
\begin{equation}
\label{eq:fullydeformedmsf}
\fqrs=\widehat{\varphi}_{q_k/p_k}(\cdots(\widehat{\varphi}_{q_{2}/p_{2}}(P_{q_1,p_1}\cdot 1)\cdot 1)\cdots)\cdot 1
\end{equation}
\end{definition}
Note that $\fqrs$ is of degree $n=p_1\cdots p_k$.
A recursive description, useful for computations, is given by the following.
\begin{lemma}\label{lem:fqrs-recursive}
Set $f_{(p_1,q_1)} := P_{q_1,p_1}\cdot 1$. For $j=2,\ldots,k$, define $f_{(p_1,q_1),\ldots,(p_j,q_j)}$ by expanding $f_{(p_1,q_1),\ldots,(p_{j-1},q_{j-1})}$ in the power sum basis $\{p_\ell\}$, replacing each $p_\ell$ by the operator $P_{\ell q_j, \ell p_j}$, and acting on $1\in\Sym_{q,t}$. Then $\fqrs = f_{(p_1,q_1),\ldots,(p_k,q_k)}$.
\end{lemma}
\begin{proof}
This is immediate from Theorem \ref{thm:slopesubalg}, which sends each power sum $p_\ell$ to $P_{\ell q_j, \ell p_j}$ under $\widehat\varphi_{q_j/p_j}$.
\end{proof}
\begin{remark}
In Definition \ref{def:fullmastersymfn}, we can replace the vacuum $1\in\Sym_{q,t}$ in \eqref{eq:fullydeformedmsf} by a modified Macdonald polynomial $\tH_\lambda$. This can be done with varying $\lambda$ at each step, to obtain "colored" versions of the refined annular symmetric functions. 


\end{remark}

We note the following homogeneity property of $\fqrs$ for further reference.
\begin{proposition}
\label{prop:homog-annular}
Let $\nabla$ denote the Bergeron--Garsia $\nabla$-operator (Definition \ref{def:nablaoperator}), which scales each $\tH_\lambda$ by $q^{n(\lambda)}t^{n(\lambda^t)}$, where $n(\lambda)=\sum_{i=1}^{\ell(\lambda)}(i-1)\lambda_i$. Then for the refined annular symmetric function $\fqrs$,
$$\nabla \cdot \fqrs = \widehat{\mathbf{f}}_{(p_1, q_1), \ldots, (p_{k-1}, q_{k-1}), (p_k, q_k+p_k)},$$
\end{proposition}
\begin{proof}
Write $\fqrs=\cO\cdot 1$ with $\cO=\widehat{\varphi}_{q_k/p_k}(\cdots)$ the iterated operator. Since $\nabla\cdot 1=1$, we have $\nabla\cdot\fqrs=(\nabla\cO\nabla^{-1})\cdot 1$. Conjugation by $\nabla$ shifts the outermost operators $P_{\ell q_k,\ell p_k}$ (which realize $\widehat{\varphi}_{q_k/p_k}$) to $P_{\ell(q_k+p_k),\ell p_k}$, using the well-known identity $\nabla P_{q_k,p_k}\nabla^{-1}=P_{q_k+p_k,p_k}$ in the representation on $\Sym_{q,t}$; equivalently the outermost slope is shifted $q_k/p_k\mapsto (q_k+p_k)/p_k$.
\end{proof}

To pass from $\fqrs$ to the superpolynomial, one proceeds as follows.
\begin{definition}\label{def:evvec}
The {\em evaluation vector} from \cite[Eq. (39)]{GN} is
\begin{equation}
\label{eq:evvec}
\fv(a):=\sum_{\mu\vdash n}\frac{\tH_\mu}{g_\mu}\prod_{\square\in \mu}(1-aq^{a'(\square)}t^{l'(\square)})
\end{equation}
where 
$$g_\mu:=\prod_{\square\in\mu}(1-q^{a(\square)}t^{-l(\square)-1})\prod_{\square\in\mu}(1-q^{-a(\square)-1}t^{l(\square)})$$
and $a'(\square),l'(\square)$ denote the coarm and coleg of a box in the Ferrers diagram.
\end{definition}
\begin{remark}
The factor $g_\mu$ is the product of the weights of the $\G_m^2$-representation $\Lambda^\bullet T_\mu^\vee\Hilb^n(\A^2)$, where the variable $a$ encodes the exterior degree.
\end{remark}
\begin{definition}
\label{def:superpolynomial}
Let $L=T(\vec{p},\vec{q})$ be an iterated torus link. The {\em superpolynomial} of $L$ is 
$$\mathbf{P}_L(a,q,t):=(\fqrs,\mathfrak{v}(a))$$
where $\mathfrak{v}(a)$ is defined in \eqref{eq:evvec} and $(-,-)$ is the geometric inner product in Definition \ref{def:innerproducts}.
\end{definition}

The refined annular symmetric function and the superpolynomial depend on three variables $q,t,a$.
 At $a=0$ the evaluation vector simplifies to 
$$\fv(0)=\sum_{\mu\vdash n}\frac{\tH_\mu}{g_\mu}.$$
Thus, for an arbitrary $f=\sum_\mu c_\mu \tH_\mu$, we have
\begin{equation}
\label{eq:evpair}
\bigl(f, \fv(0)\bigr)=\sum_\mu c_\mu,
\end{equation}
where $(-,-)$ is the geometric inner product. On the other hand, from the fact that $\langle h_\mu,h_n\rangle=1$ for all $\mu\vdash n$ we see that Eq. \eqref{eq:evpair} limits as $t\to 1$ to the expression 
\begin{equation}
\label{eq:limevpair}
\bigl(f, \fv(0)\bigr)\xrightarrow{t\to 1} \langle f, h_n\rangle.
\end{equation}
We immediately have
\begin{proposition}
The superpolynomial at $a=0$ and $t=1$ equals
$$\mathbf{P}_L(0,q,1)=\langle \frs, h_n\rangle.$$
\end{proposition}

As in Theorem \ref{thm:mastersymfn}, orbital symmetric functions, which carry information about germ expansions,
are related to the $t=1$ specialization of the refined annular symmetric function $\fqrs$. 
Motivated by this, we define

\begin{definition}\label{def:degmaster}
The {\em annular symmetric function} $\fpq \in \Sym_q$ associated to a sequence of Newton pairs $(\vec p, \vec q) = ((p_1, q_1), \ldots, (p_k, q_k))$ is
$$\fpq := \varphi_{q_k/p_k}\bigl(\cdots \varphi_{q_2/p_2}\bigl(\varphi_{q_1/p_1}(e_1)\bigr) \cdots \bigr),$$
where $\varphi_{q/p}$ is as in Definition \ref{def:mnpleth}.
\end{definition}
Defining $\fpq$ this way avoids the question of whether the refined annular symmetric function $\fqrs$ has a well-defined limit at $t=1$. The next theorem states that this limit does in fact exist, and equals $\fpq$.
\begin{theorem}\label{thm:degfpqformula}
The limit $\lim_{t \to 1} \fqrs$ exists in $\Sym_q$, and equals $\fpq$ as defined in Definition \ref{def:degmaster}.
\end{theorem}
\begin{proof}
Existence of the limit follows from Theorem \ref{thm:t1specialization}, applied iteratively: each application of $\widehat\varphi_{q_i/p_i}$ in \eqref{eq:fullydeformedmsf} corresponds to acting by $P_{kq_i, kp_i}$ on the previous step, and the matrix coefficients of these operators are regular at $t=1$. Once existence is established, the equality $\lim_{t \to 1} \fqrs = \fpq$ follows from Proposition \ref{prop:Emnk}: at $t=1$, each operator $P_{kq_i, kp_i}$ acts by multiplication, and the slope $m/n$-plethysms $\widehat\varphi_{q_i/p_i}|_{t=1}$ become $\varphi_{q_i/p_i}$.
\end{proof}
The explicit computation of these $\fpq$ is key to our computation of Shalika and Steinberg germ expansions later on.
\begin{remark}
\label{rmk:msfdef}
We set $t=1$ rather than $q=1$ to pass from $\fqrs$ to $\fpq$. From the point of view of Macdonald theory, it would possibly be more natural to set $q=1$ and then rename $t$ to $q$. Our chosen convention introduces no harm, since the $(q,t)$-formulas only differ by a transposition upon switching $q$ and $t$.

We also caution the reader that the letter $q$ is used both for the Newton pair data $(\vec{p},\vec{q})$ and for the variable in $\Sym_q$ which is further specialized to $|\F_q|$ in the harmonic analysis setting. The different roles of $q$ should always be clear from context: when $q$ appears as a subscript in $\varphi_{q_i/p_i}$ or as part of a pair $(p_i,q_i)$, it refers to the Newton pair data; otherwise it is the symmetric function variable.
\end{remark}


\subsubsection*{The Cherednik--Danilenko construction}

We now briefly recall the connection between $\cE$ and the double affine Hecke algebra (DAHA) $\HH_n$, following \cite{SVCherednik, CD}. This can be used to compare our definition of the superpolynomial with that of Cherednik--Danilenko \cite{CD}, although the results here are not used elsewhere in the paper.

The DAHA $\HH_n$ is generated over $\Q(q,t)$ by $X_i^\pm, Y_i^\pm, T_i$ subject to relations we refer the reader to \cite[Section 1.3.]{CD} for. The {\em spherical DAHA} is $\SH_n:=\eee \HH_n \eee$ where $\eee$ is the symmetrizing idempotent for the finite Hecke algebra.
\begin{definition}
\label{def:polyrep}
The {\em polynomial representation} of $\SH_n$ is $\eee \Q(q,t)[X_1^\pm,\ldots,X_n^\pm]$, which restricts to an action of $\SH_n^+$ on symmetric polynomials $\Q(q,t)[X_1,\ldots,X_n]^{S_n}$.
\end{definition}
The key structural facts about the DAHA and its polynomial representation are as follows: 
\begin{proposition}
\label{prop:dahastruct}
\begin{enumerate}
\item There is an action of the braid group $B_3=\widehat{SL_2(\Z)}$ on $\HH_n$ by algebra automorphisms \cite[Section 1.3.]{CD}, which allows one to define elements $P^{(n)}_{a,b}\in \SH_n$ for $(a,b)\in\Z^2\setminus 0$ via 
$\SL_2(\Z)$-conjugation of $P^{(n)}_{0,k}=\eee(\sum_i Y_i^k)\eee$; and 
\item there is a surjection $\cE\twoheadrightarrow \SH_n$ sending $P_{a,b}\mapsto P^{(n)}_{a,b}$ \cite[Theorem 4.6.]{SVCherednik}.
\item The action of $\cE$ on $\Sym_{q,t}$ and the action of $\SH_n^+$ on $\Q(q,t)[X_1,\ldots,X_n]^{S_n}$ intertwine the surjections $\Sym_{q,t}\twoheadrightarrow \Q(q,t)[X_1,\ldots,X_n]^{S_n}$ and $\cE^>\twoheadrightarrow \SH_n^+$ \cite[Section 5.1]{SVCherednik}.
\end{enumerate}
\end{proposition}

To define the DAHA--superpolynomial, one proceeds as before, but by replacing the slope $m/n$--plethysms by the $\SL_2(\Z)$-action on DAHA, and the evaluation vector by the following evaluation homomorphism.
\begin{definition}[\cite{CD}]\label{def:dahaeval}
The {\em evaluation homomorphism} (or {\em coinvariant}) on the DAHA is the algebra map $\{-\}_{ev}: \HH_n\to \Q(q,t)$ defined by
$$X_a\mapsto q^{-(\rho,a)}, \qquad Y_b\mapsto q^{-(\rho,b)}, \qquad T_i\mapsto t.$$
\end{definition}

\begin{definition}[\protect{\cite[Eq. (4.26)]{CD}}]\label{def:dahajones}
The {\em DAHA-Jones polynomial} of Cherednik--Danilenko is
\begin{equation}\label{eq:dahajones}
JD^{(n)}_{\vec{p},\vec{q}}(q,t):=\bigl\{\gamma_{q_1/p_1}\bigl(\cdots(\gamma_{q_{k-1}/p_{k-1}}(P_{q_k,p_k}^{(n)}\cdot 1)\cdot 1)\cdots\bigr)\cdot 1\bigr\}_{ev}
\end{equation}
\end{definition}

Our earlier definition of the superpolynomial is the $n\to \infty$ limit of the DAHA--Jones polynomial in the following precise sense.
\begin{theorem}[Cherednik's stabilization conjecture]\label{thm:cherednik-stab}
We have
$$JD^{(n)}_{\vec{p},\vec{q}}(q,t)=\mathbf{P}_L(t^n,q,t).$$
Here $L=T(\vec{p},\vec{q})$ as above.
\end{theorem}
\begin{proof}
For torus knots, this is proved in \protect\cite[Section 3.4]{GN}. Since the proof only uses the properties in Proposition \ref{prop:dahastruct}, it is easy to see that the proof carries over to general $L$.
\end{proof}
\begin{remark}\label{rmk:homfly}
The construction of \cite{CD} was in part also a motivation for Morton--Samuelson's definition of the superpolynomial, which takes on a more manifestly topological form.
More precisely, the HOMFLY--PT skein algebra $H$ of the two-torus is isomorphic to a specialization of the EHA $\cE$ by \cite[Theorem 3]{MS} (in their notation, $H \cong \cE_{q,q}$, where $q$ is one of two parameters of the elliptic Hall algebra in Burban--Schiffmann's normalization \cite{BS1}). Under this identification, Morton--Samuelson's polynomial $J^\cE(K, \lambda; q, t, u)$ from \cite[Definition 7.15]{MS}, which is up to a change of variables our $\mathbf{P}_K(q,t,a)$, specializes to the $\lambda$-colored HOMFLY--PT polynomial of $K$ and (at $u = t^N$) to the DAHA--Jones polynomial $JD^{(n)}_{\vec{p},\vec{q}}(q,t)$ from above \cite[Theorem 7.4]{MS}. 
\end{remark}


\subsection{Expansions of annular symmetric functions}
\label{subsec:weightformulas}

In this subsection, we first write down an explicit formula for the annular symmetric function $\fpq$ as a sum over ``$k$-nested Dyck paths'' (Definition \ref{def:nesteddyck} and Proposition \ref{prop:annularformula} below). We then study the limit of Theorem \ref{thm:matrixcoeff} as $t\to 1$ and prove Theorem \ref{thm:t1specialization}. The resulting $\th_\lambda$-expansion of $\fpq$ will be crucial in the explicit computation of Shalika germs in Section \ref{sec:transitionmatrices}. 
The reader interested only in results about orbital integrals and not the germs themselves can proceed directly from Proposition \ref{prop:annularformula} to the next section.

We begin with the combinatorial objects indexing the $e_\lambda$-expansion of $\fpq$. Recall the slope $m/n$-plethysm $\varphi_{m/n}$ from Definition \ref{def:mnpleth}, which on $e_k$ is the area generating function $E_{m,n,k}=\sum_{\pi\in\Dyck_{kn,km}}q^{\area(\pi)}e_\pi$ of Proposition \ref{prop:Emnk}, and which is extended to a $\Q(q)$-algebra homomorphism of $\Sym_q$. Since $\fpq=\varphi_{q_k/p_k}(\cdots\varphi_{q_1/p_1}(e_1)\cdots)$ is an iterated such plethysm, unwinding it produces a tower of Dyck paths in which the horizontal runs of each path prescribe the boxes of the paths one level deeper.

\begin{definition}
\label{def:nesteddyck}
Fix Newton pairs $(\vec p,\vec q)=((p_1,q_1),\ldots,(p_k,q_k))$. A {\em $k$-nested Dyck path} $\Pi$ for $(\vec p,\vec q)$ is a collection of rational Dyck paths arranged in $k$ levels, defined recursively as follows. There is a single {\em incoming run} of length $1$ at level $1$. At level $j$ (for $1\le j\le k$), each incoming run of length $h$ is assigned a Dyck path
$$\pi\in\Dyck_{hq_j,\,hp_j}$$
of slope $q_j/p_j$ in the box of width $hp_j$ and height $hq_j$; the horizontal run lengths $h_1(\pi),\ldots,h_{r}(\pi)$ of $\pi$, which sum to $hp_j$, are the incoming runs at level $j+1$. The level-$k$ runs of $\Pi$ thus form a composition of $p_1\cdots p_k=n$, whose associated partition we denote $\lambda(\Pi)\vdash n$. The {\em area} of $\Pi$ is the sum of the areas of all of its constituent paths,
$$\area(\Pi):=\sum_{\pi\in\Pi}\area(\pi).$$
We write $\nDyck_{(\vec p,\vec q)}$ for the set of all $k$-nested Dyck paths for $(\vec p,\vec q)$.
\end{definition}

For $k=1$ a $1$-nested Dyck path is just an ordinary Dyck path $\pi\in\Dyck_{q_1,p_1}$ and $\lambda(\Pi)$ records its horizontal runs, recovering $\varphi_{q_1/p_1}(e_1)$.

\begin{proposition}
\label{prop:annularformula}
The annular symmetric function of Definition \ref{def:degmaster} is the area generating function of $k$-nested Dyck paths:
\begin{equation}
\label{eq:annularformula}
\fpq=\sum_{\Pi\in\nDyck_{(\vec p,\vec q)}}q^{\area(\Pi)}\,e_{\lambda(\Pi)}.
\end{equation}
\end{proposition}
\begin{proof}
By Proposition \ref{prop:Emnk} and Definition \ref{def:mnpleth}, each slope $q_j/p_j$-plethysm $\varphi_{q_j/p_j}$ is the $\Q(q)$-algebra homomorphism of $\Sym_q$ with $\varphi_{q_j/p_j}(e_\nu)=\prod_i E_{q_j,p_j,\nu_i}$. Expanding $$\fpq=\varphi_{q_k/p_k}(\cdots\varphi_{q_1/p_1}(e_1)\cdots)$$ one level at a time, applying $\varphi_{q_j/p_j}$ to an incoming run of length $h$ replaces it by the sum $\sum_{\pi\in\Dyck_{hq_j,hp_j}}q^{\area(\pi)}e_\pi$ over choices of a Dyck path in the box of width $hp_j$, and the horizontal runs of $\pi$ become the incoming runs at the next level. Iterating, the terms are indexed exactly by the $k$-nested Dyck paths $\Pi\in\nDyck_{(\vec p,\vec q)}$, each contributing $e_{\lambda(\Pi)}$. As every $\varphi_{q_j/p_j}$ fixes $q$, the area exponents add across the levels, giving the total weight $q^{\area(\Pi)}=q^{\sum_{\pi\in\Pi}\area(\pi)}$.
\end{proof}

In the torus knot case, we also have the following formula at general $t$:
Let $T\in \SYT(n)$ be a standard Young tableau on $n$ letters. Recall that $T$ has $n$ boxes labeled $1, \ldots, n$ with $(q,t)$-contents $z_1, \ldots, z_n$.
Recall also from Theorem \ref{thm:matrixcoeff} that $\omega'(x)=\frac{(x-1)(x-qt)}{(x-q)(x-t)}$, $S_i':=S_{q/p}'(i):=\lfloor iq/p \rfloor - \lfloor (i-1)q/p \rfloor$, and $\xi=\frac{(q-1)(t-1)}{qt(qt-1)}$. 
The following is essentially proved in \cite[Eqs.~(49)--(50)]{GN}.
\begin{proposition}
\label{prop:fmn-explicit}
The refined annular symmetric function of the torus knot $T(n,m)$ is given by
\begin{equation}
\label{eq:gorsky-negut}
    \hat{\mathbf{f}}_{m,n}=
    \sum_{\lambda\vdash n} \xi^n \frac{\tH_\lambda}{g_\lambda}
    \sum_{\SYT(\lambda)}\frac{\prod^{n}_{i=1}z_i^{S_{m/n}'(i)}(qtz_i-1)}{\left(1-qt\frac{z_2}{z_1}\right)\cdots\left(1-qt\frac{z_n}{z_{n-1}}\right)}\prod_{i<j}\omega'\left(\frac{z_j}{z_i}\right)^{-1}.
\end{equation}
\end{proposition}
\begin{proof}
Apply Theorem \ref{thm:matrixcoeff} with $k=1$, $\mu = \emptyset$ the empty partition. Then $[1]_q = 1$, $g_\emptyset = 1$, the bracket sum $\sum_{j=0}^{0}(qt)^j[\cdots]$ collapses to $1$, and the product $\prod^{\square \in \mu} \omega'^{-1}$ is empty. So $$\langle \lambda | P_{m,n} | \emptyset\rangle = \xi^n/g_\lambda \cdot \sum_{\SYT(\lambda)} \prod_{i=1}^n z_i^{S'_{m/n}(i)}(qtz_i - 1)/\prod (1-qtz_{i+1}/z_i) \cdot \prod_{i<j}\omega'^{-1}(z_j/z_i)$$ and $\hat{\mathbf{f}}_{m,n} = P_{m,n}\cdot 1 = \sum_\lambda \langle \lambda | P_{m,n} | \emptyset\rangle \tH_\lambda/g_\lambda$ gives the stated formula.
\end{proof}
The case of iterated cables requires care. As explained in Lemma \ref{lem:fqrs-recursive}, one expands each $\tH_\lambda$ in the power sum symmetric functions, replaces each $p_k$ by the operators $P_{kq,kp}$, uses Theorem \ref{thm:matrixcoeff} again, and proceeds. 
It is not obvious whether this gives rise to any simple closed combinatorial formula at general $t$.
Luckily, the expressions simplify considerably in the $t=1$ limit. 
For ease of exposition, we will focus on individual coefficients in the above $\tH_\lambda$-expansion, indexed by standard Young tableaux.

\begin{definition}
We call the coefficient of $\tH_\lambda$ in Eq. \eqref{eq:gorsky-negut} for a fixed $T\in \SYT(\lambda)$ the $(q,t)$-{\em weight} of the SYT $T$. We will denote it by $\hat\wt(T)_{m/n}$. We denote by $\wt(T)_{m/n} := \hat\wt(T)_{m/n}\big|_{t=1}$ the $t=1$ specialization (when it is regular).
\end{definition}

\begin{lemma}
By comparison to Eq. \eqref{eq:gorsky-negut}, a convenient formula for the $(q,t)$-weight is given by
\begin{equation}\label{eq:fullweight}
\hat\wt(T)_{m/n}=\frac{\prod_{i=1}^n z_i^{S_{m/n}(i)-1}}{\prod_{i=2}^{n}(1-\frac{1}{z_i})(1-qt\frac{z_{i-1}}{z_{i}})}\prod_{i<j}\omega\left(\frac{z_i}{z_j}\right)
\end{equation} 
where $$\omega(x)=
\frac{(1-x)(1-qtx)}{(1-qx)(1-tx)}
$$ and $$S_{m/n}(i)=\lceil\tfrac{im}{n}\rceil-\lceil\tfrac{(i-1)m}{n}\rceil.$$
\end{lemma}
\begin{proof}
First, note that $S_{m/n}(n-i)=S_{m/n}'(i)$ and $\omega'(x^{-1})=\omega(x)$. In particular, reversing the labeling on the $z_i$ in Eq. \eqref{eq:gorsky-negut} we see 
$$\hat\wt(T)_{m/n}=\frac{\xi^n}{g_\lambda}\frac{\prod_i z_i^{S_{m/n}(i)}\prod_{i}(qtz_i-1)}{\prod_{i=2}^n(1-qt\frac{z_{i-1}}{z_i})}\prod_{i>j}\omega\left(\frac{z_i}{z_j}\right)^{-1}.$$ 
On the other hand, one can check that 
$$g_\lambda= z_1\cdots z_n \frac{\xi^n \prod_{i=1}^n(1-z_i^{-1})(1-qtz_i)}{\prod_{i<j} \omega\left(\frac{z_i}{z_j}\right)\omega\left(\frac{z_j}{z_i}\right)}.$$
Plugging the latter into the former, we are done.
\end{proof}
Let us now study the limit as $t\to 1$.
\begin{proposition}
\label{prop:nonvanishingcoeffs}
$~$\\
\begin{enumerate}
\item Let $T\in \SYT(\lambda).$ Then the order of vanishing of the $(q,t)$-weight $\hat\wt(T)_{m/n}$ at $t=1$ equals
\begin{equation}
\label{eq:vanishing}
v(T)=|\lambda|-\ell(\lambda)-\pi(T)
\end{equation}
where $\pi(T)$ is the number of pairs of consecutive boxes $\square_i,\square_{i+1}$ in $T$ s.t. they lie in consecutive columns. Note that this number is always $\geq 0$ and independent of $m/n$.
\item Suppose that $v(T)=0$, so that the $t=1$ weight $\wt(T)_{m/n}$ does {\em not} vanish. Then it is equal to 
$$
\wt(T)_{m/n} = \frac{\prod_{i=1}^n z_i^{S_{m/n}(i)-1}}{\prod_{i=2}^{n}(1-z_i^{-1})(1-q\frac{z_{i-1}}{z_{i}})}.
$$
Here $z_i$ are $(q,t)$-contents of boxes in $T$ {\em now specialized at $t=1$} and as in \cite{GN}, we simply ignore $\ell(\lambda)-1+\pi(T)$ zero factors in the denominator.
\end{enumerate}
\end{proposition}
\begin{proof}
From the formula for the weight in Eq. \eqref{eq:fullweight} independence of the subscript $m/n$ is clear. Looking at the denominator, we have the claimed factor and the factors
$$
\omega\left(\frac{z_i}{z_j}\right)=\frac{(1-\frac{z_i}{z_j})(1-qt\frac{z_i}{z_j})}{(1-q\frac{z_i}{z_j})(1-t\frac{z_i}{z_j})}, \; i<j
$$
At a first glance it looks like this factor is always just $1$ at $t=1$. However, this only holds
if we have the wheel conditions $\frac{z_i}{z_j}\notin \{\frac{1}{q},\frac{1}{t},\frac{1}{qt}\}$. 
For example, if $\frac{z_i}{z_j}=1/q$, then we get $$\omega\left(\frac{z_i}{z_j}\right)=\frac{(1-1/q)(1-t)}{(1-t/q)}$$ which gives a zero at $t=1$ or order 1. Similarly, for $\frac{z_i}{z_j}=1/t$ we get a zero of order 1 from $\frac{(1-1/t)(1-q)}{(1-q/t)}$, and at $\frac{z_i}{z_j}=1/qt$ we get a pole of order 1 from 
$\frac{(1-1/qt)}{(1-1/t)(1-1/q)}$. 
Since each box which has a box to the right of it contributes a zero, each box with a box above it contributes a zero, and a box with a box diagonally above it contributes a pole, we see 
that this gives a total number of
$|\lambda|-1$ zeroes.

The factors $(1-z_i^{-1})$ in the denominator vanish at $t=1$ iff $z_i=t^k$ for some $k\geq 1$, i.e. $\square_i$ is at the beginning of a row and we ignore the first factor. Therefore, they contribute a pole of order $\ell(\lambda)-1$.

Finally, the factors $(1-qtz_{i-1}/z_{i})$ contribute a pole of order $\pi(T)$ at $t=1$, since this factor vanishes at $t=1$ iff
$z_i=qz_{i-1}$. The result follows.
\end{proof}
Recall from Lemma \ref{lem:comptab} that to each composition $\alpha \vDash n$ there is associated a unique Young tableau.
\begin{lemma}\label{lem:wt-nonvanishing}
The $(q,t)$-weight $\hat\wt(T)_{m/n}$ does not vanish at $t=1$ iff $T$ comes from a composition as in Lemma \ref{lem:comptab}.
\end{lemma}
\begin{proof}
We need to show that only the tableaux coming from compositions have $v(T):=|\lambda|-\ell(\lambda)-\pi(T)=0$, cf. Eq. \eqref{eq:vanishing}. 
Note that $v(T)=0$ iff $\pi(T)=|\lambda|-\ell(\lambda)$. Additionally, the latter is an upper bound (i.e. the condition is satisfied for every box except the ends of the rows, which is obviously the maximum number of boxes), so we are looking to maximize the number of consecutive pairs of boxes in consecutive columns.

For a tableau of shape $\lambda$, coming from a composition $\alpha=\alpha_1+\cdots+\alpha_r$ we have exactly
$$\pi(T)=(\alpha_1-1)+\cdots+(\alpha_r-1)=|\alpha|-\ell(\alpha)$$ by construction. Conversely, if $\pi(T)=|\lambda|-\ell(\lambda)$, the horizontal strip coming from the top boxes in the diagram must have consecutive labels. Stripping it away gives $\alpha_r$, and we continue inductively to build $\alpha_1+\cdots+\alpha_r$.
\end{proof}
We can now prove Theorem \ref{thm:t1specialization}.
\begin{proof}[Proof of Theorem \ref{thm:t1specialization}]
By Theorem \ref{thm:matrixcoeff}, the coefficient of 
$\tH_\lambda$ in $P_{km,kn}\,\tH_\mu$ is $\langle\lambda|P_{km,kn}|\mu\rangle/g_\lambda$, 
where
$$\langle \lambda \mid P_{km,kn} \mid \mu \rangle = \frac{\xi^{kn}}{[k]}\,g_\mu\sum_{T} B_k(T)\, W_T,$$
the sum being over standard Young tableaux $T$ on the skew shape $\lambda/\mu$ --- equivalently, chains 
$\mu = \lambda^{(0)} \subset \cdots \subset \lambda^{(kn)} = \lambda$ adding one box $\square_i$ (of content $z_i$) at each step. Here
$$B_k(T)=\sum_{j=0}^{k-1}(qt)^j\prod_{a=1}^{j}\frac{z_{n(k-a)+1}}{z_{n(k-a)}}$$
is the bracket factor, a sum of $k$ Laurent monomials in the contents, and
$$W_T=\frac{\prod_{i=1}^{kn}z_i^{S'_{m/n}(i)}(qt z_i-1)}{\prod_{i=2}^{kn}\bigl(1-qt\,z_i/z_{i-1}\bigr)}\ \prod_{1\le i<j\le kn}\omega'^{-1}\!\Bigl(\frac{z_j}{z_i}\Bigr)\ \prod_{i=1}^{kn}\prod_{\square\in\mu}\omega'^{-1}\!\Bigl(\frac{z(\square)}{z_i}\Bigr)$$
collects the content- and $\omega'$-factors (read, as noted after Theorem \ref{thm:matrixcoeff}, in the reduced form of \eqref{eq:fullweight} that clears the identically-degenerate factors). We must show that $\langle\lambda|P_{km,kn}|\mu\rangle/g_\lambda$ is regular at $t=1$.

Consider first the coprime case $k=1$, $\mu = \emptyset$, which is Eq.~\eqref{eq:gorsky-negut}. Here $[1]=1$, $B_1(T)=1$, the $\mu$-factor is empty, and the coefficient of $\tH_\lambda$ is $\sum_{T\in\SYT(\lambda)}\hat\wt(T)_{m/n}$ with $\hat\wt(T)_{m/n} = \tfrac{\xi^n}{g_\lambda}W_T$. The prefactor $\xi^n$ vanishes to order $n$ at $t=1$, while $W_T/g_\lambda$ acquires poles at $t=1$ whenever a wheel condition $z_i/z_j \in \{q^{-1}, t^{-1}, (qt)^{-1}\}$ holds in the factors $\omega'(z_j/z_i)$ or in the denominator of Eq.~\eqref{eq:fullweight}. By Proposition \ref{prop:nonvanishingcoeffs} these exactly balance: $\hat\wt(T)_{m/n}$ is regular of order $v(T) = |\lambda| - \ell(\lambda) - \pi(T) \geq 0$.

For the general case $k\ge1$ and arbitrary $\mu$ we account for each factor's order at $t=1$ separately. The boxes $\square\in\mu$ will be referred to as "frozen", as their $(q,t)$-contents $z(\square)=q^{a'(\square)}t^{l'(\square)}$ are the same in every tableau $T$ of the sum. Thus $\mu$ enters the sum only through the constant $g_\mu$ and through the mixed factors $\prod_i\prod_{\square\in\mu}\omega'^{-1}(z(\square)/z_i)$.

The factor $\xi^{kn}$ vanishes to order $kn$, and by Remark \ref{rmk:glambda} $\operatorname{ord}_{t=1}g_\lambda=\ell(\lambda)$. Also, $[k]=\frac{(q^k-1)(t^k-1)}{(q-1)(t-1)}\to k\frac{q^k-1}{q-1}\ne0$, so $1/[k]$ is regular and nonzero at $t=1$. The bracket factor $B_k(T)$ is a sum of $k$ Laurent monomials in the contents so regular and nonzero.

It remains to bound the order of $W_T$. 
The genuinely new factors there compared to the $k=1$ case are the $\omega'^{-1}(z(\square)/z_i)$ mixing boxes from $\lambda$ and $\mu$, controlled by the following lemma whose simple proof we omit.
\begin{lemma}\label{lem:frozenreg}
For every added box $\square_i$ of $\lambda/\mu$ and every frozen box $\square\in\mu$, the factor $\omega'^{-1}(z(\square)/z_i)$ is regular and nonzero at $t=1$.
\end{lemma}

Since all the factors in each term of the sum are regular at $t=1$, every term is regular at $t=1$, and hence so is $\langle\lambda|P_{km,kn}|\mu\rangle/g_\lambda$.

Granted regularity, each $P_{km,kn}$ preserves the lattice $V_R$, and the values $\langle \lambda \mid P_{km,kn} \mid \mu \rangle\big|_{t=1}$ are the matrix coefficients, in the basis $\{\th_\lambda\}$, of the specialized operator $P_{km,kn}|_{t=1}$ on $\Sym_q$. Together with the multiplications by the power sums $p_k$ (the operators $P_{k,0}$), these give the $t=1$ specialization of the action of Theorem \ref{thm:fockrep}.
\end{proof}

At $t=1$ we then finally have 
\begin{proposition}
\label{prop:limitformula}
\begin{equation}
\label{eq:Pmnattequals1}
P_{m,n}|_{t=1}\cdot 1=\sum_{\alpha\vDash n} \wt(\alpha)_{m/n} \th_\alpha =\sum_{\alpha\vDash n}\frac{(-1)^{n-\ell(\alpha)}z_1^{S_1}\cdots z_n^{S_n} }{c_{\alpha-1}(q)} \th_\alpha
\end{equation}
where $S_i=\lceil \frac{im}{n}\rceil -\lceil \frac{(i-1)m}{n}\rceil$. Here $(m,n)=1$.
\end{proposition}
\begin{proof}
From the second part of Proposition \ref{prop:nonvanishingcoeffs} we have that the denominator of $\wt(\alpha)_{m/n}$ is
$$
\prod_{i=2}^{n}\frac{1}{(1-z_i^{-1})(1-qz_{i-1}/z_{i})}
$$
with the convention that zero factors are ignored.

Arranging the $q,t$-contents at $t=1$ into a vector by reading the tableau box by box, we write
$$
z(\alpha)=(1,q,\ldots,q^{\alpha_1-1},1,q,\ldots,q^{\alpha_2-1},\ldots,1,\ldots,q^{\alpha_r-1})
$$ 
and by definition $z(\alpha)_i=z_i$.

It is easy to see that the $1-qz_i/z_{i+1}$ factors are only nonvanishing at the ends of the parts of $\alpha$ so this becomes
$$
\prod_{j=1}^{r}\frac{1}{(1-q^{-1})\cdots (1-q^{-\alpha_j-1})}\prod_{j=1}^{r-1}\frac{1}{1-q^{\alpha_j}}=\frac{(-1)^{n-\ell(\alpha)}z_1\cdots z_n }{(1-q)^{n-1}[\alpha_1]!\cdots [\alpha_{r-1}]![\alpha_{r}-1]!}=\frac{(-1)^{n-\ell(\alpha)}z_1\cdots z_n}{c_{\alpha-1}(q)}
$$
where $[m]!=\prod_{i=1}^{m}(1-q^m)/(1-q)$ and $\alpha-1$ is the composition where we remove $1$ from the last part. 

Since the numerator was $z_1^{S_{m/n}(1)-1}\cdots z_n^{S_{m/n}(n)-1}$ we get the result.
For further reference, we will denote the coefficient for a fixed $\alpha$ the {\em weight of $\alpha$}:
\begin{equation}
    \label{eq:compweightcoprime}
    \wt(\alpha)_{m/n}=\frac{(-1)^{n-\ell(\alpha)}z_1^{S_{m/n}(1)}\cdots z_n^{S_{m/n}(n)}}{c_{\alpha-1}(q)}
\end{equation}
\end{proof}
In order to write down the transition matrix of Shalika germs, we will also need the case when $m,n$ are not coprime.
This is the $t\to 1$ limit of the formula in Theorem \ref{thm:matrixcoeff} at $\mu=\emptyset$. 
\begin{proposition}
\label{prop:limitformulafull}
\begin{equation}
\label{eq:noncoprimeformula}
P_{km,kn}\cdot 1=\sum_{\alpha\vDash kn} \wt(\alpha)_{m/n}\th_\alpha
\end{equation}
where for $\alpha\vDash kn$
\begin{equation}
\label{eq:compweightfull}
\wt(\alpha)_{m/n}=\left(1+\sum_{j=1}^{k-1} q^j \frac{z_{(k-j)n}\cdots z_{(k-1)n}}{z_{(k-j)n+1}\cdots z_{(k-1)n+1}}\right) \frac{(-1)^{kn+k-1-\ell(\alpha)}z_1^{S_{m/n}(1)}\cdots z_{kn}^{S_{m/n}(kn)}}{c_{\alpha-1}(q)}
\end{equation}
\end{proposition}
\begin{proof}
Comparing to Theorem \ref{thm:matrixcoeff} and Proposition \ref{prop:limitformula} this is proved exactly in the same way, but we also have the coefficient 
$$1+\sum_{j=1}^{k-1}q^j \frac{z_{n(k-j)}\cdots z_{n(k-1)}}{z_{n(k-j)+1}\cdots z_{n(k-1)+1}}.$$ 
\end{proof}

\section{Explicit computation of orbital symmetric functions and germs}
\label{sec:combinatorics}

In this section, we state and prove the inductive combinatorial formula for the Shalika germs and the Steinberg germs, as well as the orbital integrals themselves. Our methods on the harmonic analysis side are heavily based on results of \cite{W2}. Currently, this section can be regarded as the most technical part of our computations, but we also hope our results give insight to the techniques in \cite{W2}.

\begin{definition}
\label{def:gamma}
Let $m\geq n\geq 0$, and consider two partitions $\lambda, \mu$ so that $\lambda \vdash m$ has $n$ parts and $\mu\vdash n$. Consider the set $\Upsilon^\mu_\lambda\subset \Z_{\ge 0}^n$ defined by

$$\Upsilon^\mu_\lambda:=\Bigg\{(d_1,\ldots,d_n)\in\Z_{\ge 0}^n\;|\;(\sum_{i=1}^k d_i)(\sum_{i=k+1}^n\lambda_i)-(\sum_{i=1}^k\lambda_i)(\sum_{i=k+1}^nd_i)$$

$$\begin{cases}
=0, & \text{if there is a $j$ so that $k=\sum_{i=1}^j\mu_i$}\\
>0, & \text{otherwise}
\end{cases}, \; \forall k=1,\ldots, n-1 \Bigg\}
$$
\end{definition}

\begin{example}
Let $m=n$, $\mu=(n)$ and $\lambda=(1^n)$. Then 
$$\Upsilon^{(n)}_{(1^n)}=\left\lbrace\vec{d} \;\mid\; (n-k)\cdot \sum_{i=1}^kd_i-k\cdot \sum_{i=k+1}^nd_i>0, \forall k=1,\ldots, n-1\right\rbrace$$
\end{example}
\begin{proposition}
\label{prop:dyckpaths}
 The subset $D_{d,n}\subseteq \Upsilon^{(n)}_{(1^n)}$ given by $\vec{d}$ with $\sum_i d_i=d$ is in bijection with the slope $d/n$ rational Dyck paths $\Dyck_{d,n}^<$ {\em strictly under the diagonal}.
\end{proposition}
\begin{proof}
Giving a Dyck path is the same as giving the sequence of its horizontal steps. If we are looking at $d/n$-Dyck paths, these steps have to sum to $d$ and there are at most $n$ steps, which gives us a sequence $(d_1,\ldots, d_n)$. Since such a Dyck path lies above the line with slope $d/n$, we must have $$(n-k)(\sum_{i=1}^kd_i)-k(d-\sum_{i=1}^k d_i)=
n\sum_{i=1}^k d_i-kd>0
$$ for all $k$ (if we allowed all Dyck paths, for noncoprime $d,n$ equality could also hold). The converse is clear.
\end{proof}
\begin{example}
\label{ex:dyck}
Let $n=4$, $d=3$. The allowed sequences are
$$ (3,0,0,0), (2,1,0,0), (2,0,1,0), (1,2,0,0),(1,1,1,0)$$
and these correspond to the Dyck paths

\begin{center}
\begin{tikzpicture}[scale=0.7]
    \NEpath{0,0}{3}{4}{0,0,0,1,1,1,1};
    \NEpath{4,0}{3}{4}{0,0,1,0,1,1,1};
    \NEpath{8,0}{3}{4}{0,0,1,1,0,1,1};
    \NEpath{12,0}{3}{4}{0,1,0,0,1,1,1};
    \NEpath{16,0}{3}{4}{0,1,0,1,0,1,1};
\end{tikzpicture}
\end{center}
\end{example}

\vspace{1cm}
For any $P=\sum_{\vec{d}}a_{\vec{d}}X^{\vec{d}} \in \Z[[X_1,\ldots,X_n]]$ we define the truncation
\begin{equation} 
\label{eq:restrictedseries}
\Upsilon_\lambda^\mu(P):=\sum_{\vec{d}\in \Upsilon^\mu_\lambda}a_{\vec{d}}X^{\vec{d}}
\end{equation}
\begin{remark}
Our $\Upsilon^\mu_\lambda$ is equal to the ''polynomial" component of the set $\Gamma^\mu_{\lambda,\Z}$ defined in \cite[I 10]{W2}. It appears that this polynomial component constitutes the primary (if not the exclusive) usage of $\Gamma^\mu_{\lambda,\Z}$ in \cite{W2}. We have avoided this notation in order not to get it confused with the Shalika germs $\Gamma_\lambda(-)$. Note also that on pages 856 and 878 of {\em loc. cit.} the confusing notation $\Gamma^{\mu_P}_\lambda$ is used, but this seems to be a misprinted $\Gamma^{\mu}_\lambda P$. 

In general, one should think of the elements of $\Gamma^\mu_{\lambda,\Z}$ and $\Upsilon^\mu_\lambda$ Lie-theoretically as follows. $\Z^n$ is the weight lattice of $GL_n$, and each collection of $n$ integers $\lambda_1,\ldots,\lambda_n$ gives a linear form on $\Z^n$ as defined above. On the other hand, $\mu$ gives a parabolic subgroup of $GL_n$, and the (in)equalities above decide that this linear form should (not) vanish on the relative root subspaces of the corresponding Levi subgroup. This cuts out a cone in the apartment of $T$. Fixing the coordinatewise sum is intersecting this cone with an affine hyperplane. Further restricting to the nonnegative part is intersecting with yet another cone, i.e. the positive orthant. We remark that the definition of $\Upsilon^\mu_\lambda$ is related to the definition of ''Hecke-regular functions" in \cite[Section 4]{AC}. 
\end{remark}

\subsection{Waldspurger's recursion and the elliptic Hall algebra}
\label{sec:w2}
The goal of this subsection is to recall the recursive computation of the Steinberg germs from \cite{W2} and in the inertially elliptic case to compare it to the recursive definition of $\fpq$ from Section \ref{sec:knots}. We also extend the construction to apply to general tamely ramified elements, even when there is no knot in the picture.

To start with, if we suppose $\gamma\in GL_n(F)$ is elliptic, tamely ramified and compact, the recursion of \cite{W2} for determining the Steinberg germs $\Gamma_\lambda^{St}(\gamma)$ proceeds in two steps, as explained in VII 7. of {\em loc. cit.} We now recall this process in our notation. As our eventual goal is to compute the Steinberg and Shalika germs of regular semisimple tamely ramified $\gamma\in \fg\fl_n$ or $\gamma\in GL_n$ by a combination of homogeneity and the computation for compact elements in the group done below, we set up some definitions. Their utility will immediately become clear.

\begin{definition}
\label{def:ellipticmastersymfn}
We now define the {\em orbital symmetric function of $\gamma$}, $\fgamma\in \Sym$ for $$\gamma\in GL_n(\cO), GL_n(F)_c, \fg\fl_n(F)$$ respectively. 
\begin{enumerate}
\item  Let $\gamma\in GL_n(\cO)$ be regular semisimple. By Corollary \ref{cor:shalika} we have the group Shalika germs $\Gamma_\lambda(\gamma)$ for $\gamma$ and the Steinberg germs $\Gamma^{St}_\lambda(\gamma)$. Let the orbital symmetric function of $\gamma$ be 
$$\fgamma:=\sum_\lambda \Gamma_{\lambda^t}(\gamma)\th_\lambda$$ or equivalently by Propositions \ref{prop:wpair1}, \ref{prop:wpair2},
$\fgamma=\sum_\lambda \Gamma^{St}_\lambda(\gamma)e_\lambda$. Note the transposition in the Shalika expansion, forced upon us by Proposition \ref{prop:wpair1}.
\item For general compact $\gamma\in GL_n(F)_c$, not necessarily with $\val(\det(\gamma))\not=0$, we define 
$$\fgamma:=\sum_\lambda \Gamma^{St}_\lambda(\gamma)e_\lambda$$
\item Finally, if $\gamma\in \fg\fl_n(F)$, define $$\fgamma:=\sum_\lambda \Gamma_{\lambda^t}(\gamma)\th_\lambda$$ where $\Gamma_\lambda(\gamma)$ are the Lie algebra Shalika germs from Theorem \ref{thm:shalika-for-Lie-algebra}.
\end{enumerate}
\end{definition}
\begin{remark}
\label{rem:sizeofresfield}
    These symmetric functions are elements of $\Sym$, the ring of symmetric functions over $\C$. As follows from Theorem \ref{thm:mainformula}, they only depend on certain discrete invariants attached to $\gamma$ in the tamely ramified case. Moreover, the coefficients in any of the above expansions are rational functions of $q$, the size of the residue field of $F$, which only have poles at roots of unity. It therefore does no harm to consider $\fgamma$ as an element of $\Sym_q$, the ring of symmetric functions over $\C(q)$. 
\end{remark}
\begin{remark}
If $\gamma\in GL_n(F)$ is topologically unipotent, i.e. $\gamma=1+\gamma'$ where $\gamma'$ is topologically nilpotent, then it follows from Proposition \ref{prop:germcomparison} that the Lie algebra version of $\fgamma$ coincides with the group version $\fgamma$, relating definitions (1) and (3). In this case we also have $\fgamma=\mathbf{f}_{\gamma'}$ on the Lie algebra. Below, we tacitly avoid keeping track of which $\fgamma$ we work with, as it should be clear from the context. Note that in case (2) there is no obvious analog of the Shalika germs. Of course, one may define these via a change of basis a posteriori, but the harmonic analysis meaning is unclear.
\end{remark}
\begin{definition}
\label{def:dyckgerms}
Let $\gamma$ be as in any of the cases (1)--(3) above. The coefficients of the expansion of $\fgamma$ in the complete homogeneous symmetric functions are called the {\em Dyck germs} of $\gamma$. In the group case, one can think of these as the expansions of the corresponding orbital integrals as linear combinations of Fourier transforms of nilpotent orbital integrals truncated to compact elements, although we do not use this fact.  
\end{definition}

Let us recall the setup of \cite[S\'ection VI--VII]{W2}, in slightly simplified form (the simplification being that for us $E=F$ and $r=1$ in the notation of {\em loc. cit.}). Now $F$ is a nonarchimedean local field (of arbitrary characteristic, see Remark \ref{rmk:characteristic}), $\gamma\in GL_n(F)$ is elliptic and tamely ramified so that
$F(\gamma)/F$ is a tamely ramified extension of $F$. Here for any semisimple $\gamma\in GL_n(F)$, we can realize $F(\gamma)$ as the commutative subalgebra of $\fg\fl_n(F)$ given as the center of the centralizer of $\gamma$ in $\fg\fl_n(F)$, and $F(\gamma)^{\times}$ is the corresponding torus in $GL_n(F)$. 
\begin{definition}
\label{def:fcuspidal}
Let $F'/F$ be a tamely ramified extension with ramification index $e=e_{F'/F}$. We call $\delta\in (F')^\times$ cuspidal for $F'/F$ if $(\text{val}_{F'}(\delta),e)=1$ and if the reduction of $t^{-\text{val}_{F'}(\delta)}\delta^e$ in the residue field of $\cO_{F'}$ generates the residue field over that of $\cO_F$ for any uniformizer $t\in F$.
\end{definition}    

We will also need the following \cite[VII 4.]{W2}
\begin{lemma}
\label{lem:factoreta}
Let $E/F$ be a non-trivial tamely ramified extension and $\gamma\in 1+\mathfrak{m}_E$ be such that $F(\gamma)=E$. Then we can always write $\gamma=\eta(1+\delta\gamma')$ with $\eta\in 1+\mathfrak{m}_F$, $\gamma'\in 1+\mathfrak{m}_E$ and $\delta\in E$, such that $\delta$ is cuspidal for $F':=F(\delta)/F$ with $F'/F$ a non-trivial extension.\end{lemma}
Note that in both lemmas, while $\gamma', \delta$ are not uniquely determined by $\gamma$, the integers $e=e_{F'/F}$, $f=[\cO_{F'}/\fm_{F'}:\cO_F/\fm_F]$, $d=\val_{F'}(\delta)$ and the extension $F'$ are all uniquely determined by $\gamma$. In fact, more numerical invariants are preserved. Let us define
\begin{definition}
\label{def:RV}
Let $E/F$ be a tamely ramified extension of degree $n$ and $\gamma\in E$ be such that $E=F(\gamma)$. Take $F^s$ any separable closure of $F$ and let $\gamma_1,...,\gamma_n$ be all conjugates of $\gamma$ in $F^s$. We define $RV_F(\gamma)$ to be the {\em multiset} $RV_F(\gamma):=\{\val_F(\gamma_i-\gamma_j)\;|\;1\le i<j\le n\}$.
We also define  $$\Xi_F(\gamma) = \left(\sum_{r\in RV_F(\gamma)} r\right) - (n - [E^{ur}:F])/2$$ where $E^{ur}$ is the maximal unramified subextension of $E/F$.
\end{definition}
\begin{remark}
\label{rmk:dimspgamma}
The underlying set of our multiset $RV_F(\gamma)$ already appeared as the set $RV(\gamma)$ in Definition \ref{def:puiseux}. In particular, $RV_F(\gamma)$ determines the {\it Puiseux pairs} and {\it Newton pairs} of $\gamma$ as in Definition \ref{def:puiseux}.

We note that when $\gamma$ is tamely ramified, $\Xi(\gamma)$ equals the dimension of the affine Springer fiber $\Sp_\gamma$.
\end{remark}
Given a multiset $S$ of rational numbers and $\alpha\in \Q$, we denote by $\alpha S$ and $S+\alpha$ the multisets of the same cardinality as $S$, given by multiplying each element by $\alpha$, resp. adding $\alpha$ to each element. 

\begin{lemma}\label{lem:reductionnumerics} Let $\gamma$, $\delta$ and $\gamma'$ be as in Lemma \ref{lem:factoreta}. Let $e$ and $f$ be the ramification index and residual degree of $F'/F$ in the corresponding notations. Write $n':=n/ef$.
\begin{enumerate}
    \item We have $RV_F(\gamma)=(\frac{1}{e}RV_{F'}(\gamma')^{ef}\sqcup\{0\}^{n(n-n')/2})+\val_F(\delta)$. 
    \item We have $\Xi_F(\gamma)=f\cdot\Xi_{F'}(\gamma')+\frac{n(n-1)}{2}\val_F(\delta)-\frac{n-n/e}{2}$. 
    \item When $E=F(\gamma)\supsetneq F$ and $E/F$ is totally ramified, the sequence of Newton pairs of $\gamma'$ is given by deleting the last Newton pair $(p_k,q_k)$ of $\gamma$, where $q_k/p_k=\val_F(\delta)$, or more precisely $p_k=e$ and $q_k=e\cdot \val_F(\delta)$.
\end{enumerate}
\end{lemma}

To carry out the recursive computation of $\fgamma$, we first assume our elliptic, tamely ramified group element $\gamma$ is also of the form $\gamma\equiv 1$ mod $\mathfrak{m}_{F(\gamma)}$, following Lemma \ref{lem:factoreta}. 
We also define $\gamma''=\delta\gamma'$, similar to \cite[VII 4.]{W2}.

\begin{example} 
\label{ex:totram}
In the totally ramified case, assuming $\gamma=1+a_d u^{nr_k}+\cdots+a_1 u^{nr_1}$ with $r_k<\cdots < r_1$ (compare to Eq. \eqref{eq:puiseux}), we can write $\gamma=\eta(1+\delta\gamma')$ with $\eta$, $\delta$, $\gamma'$ satisfying the above conditions, in the following way: If $r_d\not\in\Z$, then we simply take $\eta=1$, $\delta=a_du^{nr_k}$, and $\gamma'=1+\frac{a_{k-1}}{a_k}u^{n(r_{k-1}-r_k)}+\cdots$  In general, $\eta$ takes care of possible central translations of $\gamma$, which don't affect the computation by Theorem \ref{thm:shalika-for-Lie-algebra}. For example, if $\gamma=1+u^4+u^6+u^7$ in $GL_4$, we can take 
$\eta=1+u^4$, in which case $$
\gamma\eta^{-1}=1+u^6+u^7-u^{10}-u^{11}+\cdots=1+u^6(1+u-u^4-u^5+\cdots)
$$ so that we can take $\delta=u^6$ and $\gamma'=1+u-u^4-u^5+\cdots$

\end{example}
The first step of the recursion is essentially \cite[Prop. VII 5.]{W2}, recalled in Proposition \ref{prop:wrecursion} below. To state it, we need some notation.
\begin{definition}
\label{def:xm}
Let $X(T)=\sum_{i\geq 0} h_i T^i\in \Sym[[T]]$.
Fix $\lambda'\vdash n'$ and number the boxes in its diagram $1,\ldots, n'$, going from left to right and bottom to top.
Let $$P=\prod_{k=1}^{n'} X(X_k)$$ 
and let $x_n(d\lambda',q)\in\Sym[q^{\pm1}]$ be the coefficient of $T^n$ in the series 
$$\Upsilon_{(1^{n'})}^{\lambda'}(P)$$ evaluated at $X_{\square}=q^{a'(\square)-(\lambda'_{l'(\square)+1}-1)/2}T$, where $\Upsilon_{(1^{n'})}^{\lambda'}$ is as defined in Eq. \eqref{eq:restrictedseries}, $a'(\square), l'(\square)$ are the coarm and coleg lengths as in Eq. \eqref{eq:armslegs} and $\square$ runs over the boxes in $\lambda'$. See also \cite[Section I 10., V 12.]{W2}. 
\end{definition}

\begin{proposition}
\label{prop:wrecursion}
Write $\gamma=1+\delta\gamma'$ and $\gamma'':=\delta\gamma'$ as in Lemma \ref{lem:factoreta}. Let $n'=n/e$, where $e$ is as after Lemma \ref{lem:factoreta}. Evaluating at $q=|k|$ we have the following equality of symmetric functions in $\Sym$:
\begin{equation}
\label{eq:recursion-hside}
\sum_{\lambda\vdash n} \Gamma^{St}_\lambda(\gamma)e_\lambda=(-1)^{n-n'}q^{(n'n d-n)/2}\sum_{\lambda'\vdash n'} (-1)^{n'-\ell(\lambda')}\Gamma^{St}_{\lambda'}(\gamma'')x_n(d\lambda',q)
\end{equation}
where $x_n(d\lambda',q)$ is as defined in Definition \ref{def:xm} and the Steinberg germs are as in Theorem \ref{thm:steinberg}.
\end{proposition}
We refer the reader to \cite[Lemme VII 5., Lemme V 12.]{W2} for details. We warn the interested reader that there are some printing errors in the latter lemma, e.g. the second displayed equation on p. 880 should have a subscripted $X_2$, a sentence after it there seems to be an extra "A" in front of $\lambda$, and on line 7 of p. 881 another subscript seems to have gone astray. Also in the statement the "$q$" should not be a subscript of $t(\lambda)$ but on the same line as $(-1)$. Also in the former Lemma, the exponent of $q$ should have a $-n/2$ instead of $-n/e$.
 
Proposition \ref{prop:wrecursion} does not quite yet cover the induction for general tamely ramified compact $\gamma\in GL_n(F)$ as outlined in \cite[Section VII 7.]{W2}, as the right hand side of \eqref{eq:recursion-hside} has $\gamma''$ which is not congruent to $1$ for an essential reason; $\val_F(\gamma'')$  is in general not an integer.

Write $\gamma'':=\delta\gamma'$. and let $d,f$ be as in Lemma \ref{lem:reductionnumerics} so that $\val_{F}(\det_G(\delta))=nd/e$. By Theorem \ref{thm:steinberg} we have the Steinberg germs $\Gamma_{\lambda}^{St}(\gamma'')$ for each $\lambda\vdash n/e$. Note that even though $\gamma''\in GL_n(F)$, these are partitions of $n/e$. In contrast, we have germs $\Gamma_{\lambda'}^{St}(\gamma')$ for $\lambda'\vdash n/ef$ where $\gamma'$ is now thought of as an element of $G':=Z_{GL_n}(\delta)\cong GL_{n/ef}(F(\delta))$. The associated orbital symmetric functions are
\[\mathbf{f}_{\gamma''}=\sum_{\lambda\vdash n/e}\Gamma_\lambda^{St}(\gamma'')e_\lambda\text{ for }\gamma''\in GL_n(F)\]
and similarly 
\[\mathbf{f}_{\gamma'}=\sum_{\lambda'\vdash n/ef}\Gamma_{\lambda'}^{St}(\gamma')e_{\lambda'}\text{ for }\gamma'\in GL_{n/ef}(F')\]
We have
\begin{proposition}
\label{prop:unramifiedextension}
Let $\gamma''$ and $\gamma'$ be as above. The orbital symmetric functions of $\gamma''$ and $\gamma'$ are related by
$$\mathbf{f}_{\gamma''}=(-1)^{n-n/e}\tau_f(\mathbf{f}_{\gamma'})$$
where $\tau_f$ is the plethysm/Adams operation defined on symmetric functions by $\tau_f: p_r\mapsto p_{rf}$ for all $r\in \Z_{\geq 0}$. Note that when $f=1$ this allows us to compare the Steinberg germs of $\gamma''$ and $\gamma'$ directly up to a sign.
\end{proposition}
\begin{proof}
This is a direct translation of the $r=1$ case of \cite[Proposition VII 2.]{W2}, where the sign is $(-1)^{nd-nd/e}$. But $nd-nd/e\equiv n-n/e $ mod $2$ whenever $\gcd(d,e)=1$. Note that our $\gamma''$ was denoted by $\gamma$ {\it loc. cit.}.
\end{proof}

Let us now compare the construction of $\fpq$ from Section \ref{sec:knots} to $\fgamma$ introduced above, and extend the construction to apply to general tamely ramified extensions of local fields, even when there is no knot in the picture.
We start with a Lemma.
\begin{lemma}
\label{lem:comparedycktoW2} For any $d\in\Z_{\ge 0}$, $n,e\in\Z_{\ge 1}$, with $e|n$ and $\lambda'\vdash n/e$, we have 
$$x_{n}(d\lambda',q)=q^{-(ed\sum_{i=1}^{\ell(\lambda')}(\lambda_i')^2-en')/2}\prod_{i=1}^{\ell(\lambda')}\left(
\sum_{\pi \in \Dyck^<_{e\lambda_i',d\lambda_i'}}q^{\coarea(\pi)}h_\pi\right)$$ where we only sum over Dyck paths strictly under the diagonal. Here $n'=n/e$ as before. (Note that the $n'$ in \cite{W2} is $n/ef$.)
\end{lemma}
\begin{proof}
Since $$x_{n}(d\lambda',q)=\prod_{i=1}^{\ell(\lambda')}x_{e\lambda_i'}(d\lambda_i',q)$$ (see e.g. the second displayed equation of p. 883 in \cite{W2}),
we may restrict to a single factor in the product. 

Proposition \ref{prop:dyckpaths} tells us that the terms in $\Upsilon^{d\lambda'_i}_{1^{n'}}(P)$ which contribute to the coefficient $x_{e\lambda_i'}(d\lambda_i',q)$ of $T^{e\lambda_i'}$ 
of the evaluation of this series at $T_\square=q^{a'(\square)-(d\lambda_i'-1)/2}T$ are in bijection with $(d{\lambda_i'},e\lambda_i')$-Dyck paths strictly under the diagonal.

Collecting these terms, we evaluate $\Upsilon^{\lambda'}_{1^{n'}}(P)$ from Definition \ref{def:xm} at $$T_\square=q^{a'(\square)-(\lambda_{l'(\square)+1}-1)/2}=q^{a'(\square)-(d\lambda_i'-1)/2}T$$ since $l'(\square)+1=1$ for all boxes and $\lambda=d\lambda_i'$. (Note that in the definition of $P$ there is a $+1/2$ instead of a $-1/2$ as in \cite[V 12]{W2}, because our $a'(\square)$ equals $i-1$ from {\em loc. cit.}.)  We can thus factor out $q^{-(d\lambda_i'-1)/2}$. Since we are interested in the coefficient at $e\lambda_i'$, we get an overall factor of $q^{-e\lambda_i'(d\lambda_i'-1)/2}$. Taking the product over $i$ gives 
$q^{-(ed\sum_{i=1}^{\ell(\lambda')}(\lambda_i')^2-en')/2}$ as desired. 

It is not hard to see that the coarea of the Dyck path corresponding to a term $h_\pi T_\pi$ contributing to $x_{e\lambda_i'}(d\lambda_i',q)$ where $T_\pi=T_1^{h(\pi)_1}\cdots T_{d\lambda_i'}^{h(\pi)_{d\lambda_i'}}$ is exactly $\coarea(\pi)=\sum_k (k-1) h(\pi)_k$. Here $h(\pi)_k$ denote the horizontal steps of $\pi$ and $h_\pi$ is the homogeneous symmetric function associated to the partition of $e\lambda_i'$ given by the horizontal steps of $\pi$.

\end{proof}
 Recall the operators  $E_{d,e,\lambda'}$ from Proposition \ref{prop:Emnk}. Similarly, we define 
\begin{equation}
    \label{def:Hlambda}
    H_{d,e,\lambda'}=\prod_i H_{d,e,\lambda_i'}=\prod_{i}\left(\sum_{\pi \in \Dyck_{e\lambda_i',d\lambda_i'}}q^{\coarea(\pi)}h_\pi\right)
\end{equation}
Motivated by Lemma \ref{lem:comparedycktoW2}, we also define $$H^<_{d,e,\lambda'}:=\prod_{i=1}^{\ell(\lambda')}
\left(\sum_{\pi \in \Dyck_{e\lambda_i',d\lambda_i'}^<} q^{\coarea(\pi)}h_\pi\right)$$
We want to relate these functions to the slope $d/e$-plethysms $\varphi_{d/e}$ on symmetric functions, as defined in Section \ref{sec:knots}. Note that these operators involve the area statistic on Dyck paths, rather than the coarea.
\begin{lemma}
\label{lem:coareatoarea}

We have 
$$q^{n'nd/2-n/2}x_n(d\lambda',q)=q^{\frac{(dn'-1)(en'-1)+n'-1}{2}}\prod_{i=1}^{\ell(\lambda_i')}\left(\sum_{\pi\in\Dyck_{e\lambda_i',d\lambda_i'}^<}q^{-\area(\pi)}h_\pi\right)$$

\end{lemma}
\begin{proof}
    We have 
$$H^<_{d,e,\lambda'}=q^{\sum \delta_{\lambda_i'}} \prod_i \left(\sum_{\pi \in \Dyck_{e\lambda_i',d\lambda_i'}^<} q^{-\area(\pi)}h_\pi\right)$$
where $\delta_{\lambda_i'}=\frac{(d\lambda_i'-1)(e\lambda_i'-1)+\lambda_i'-1}{2}$.
Then $$\sum_{i}\delta_{\lambda_i'}=
(ed \sum_i (\lambda_i')^2+\sum_i (-d-e+1)\lambda_i')/2
$$ and therefore by Lemma \ref{lem:comparedycktoW2} $$
x_n(d\lambda',q)=q^{-(ed\sum_{i=1}^{\ell(\lambda')}(\lambda_i')^2-en')/2}H^{<}_{d,e,\lambda'}
=q^{(1-d)n'/2}\prod_i \left(\sum_{\pi \in \Dyck_{e\lambda_i',d\lambda_i'}^<} q^{-\area(\pi)}h_\pi\right)
$$
Further multiplying this by $q^{n(n'd-1)/2}$ we get
$q^{\frac{nn'd-n-dn'+n'}{2}}=q^{\frac{(n-1)(dn'-1)+n'-1}{2}}$ in front.

\end{proof}

For convenience, let us denote $$
H^-_{d,e,\lambda'}=\prod_i\left(\sum_{\pi\in\Dyck_{e\lambda_i',d\lambda_i'}}q^{-\area(\pi)}h_\pi\right)
$$
and
$$H^{-,<}_{d,e,\lambda'}=\prod_i\left(\sum_{\pi\in\Dyck_{e\lambda_i',d\lambda_i'}^<}q^{-\area(\pi)}h_\pi\right)$$

Consider now the operator on symmetric functions which takes each $h_{\lambda'}$ and replaces it by $H^-_{d,e,\lambda'}$. 
We can think of this as the conjugation of $\varphi_{d/e}$ by the involution $\omega$, together with negating all the powers of $q$ that appear in the definition. 

More precisely, take an elliptic tamely ramified $\gamma\in G(\cO)$ as above and let $\gamma',\gamma''$ be as in Proposition \ref{prop:wrecursion}. By the process just described, we compute
\begin{equation}
\label{eq:recursion-eside}
q^{\frac{(dn'-1)(en'-1)+n'-1}{2}}\omega \varphi_{d/e}|_{q\mapsto q^{-1}} \omega^{-1} (\mathbf{f}_{\gamma''})=q^{\frac{(dn'-1)(en'-1)+n'-1}{2}}\sum_{\lambda'\vdash n/e} \sigma_{\lambda'}(\gamma'')
H_{d,e,\lambda'}^-
\end{equation} 
Where $\sigma_{\lambda'}$ are the ''Dyck germs" of Definition \ref{def:dyckgerms} that appear in the $h_{\lambda'}$-expansion of $\mathbf{f}_{\gamma''}$ and the relevance of the $q$-power will become clear soon. 
Our main technical result is
\begin{theorem}
\label{thm:comparerecursion}
The right-hand sides of Eqs. \eqref{eq:recursion-hside} and \eqref{eq:recursion-eside} are equal up to a sign. More precisely,
$$q^{n'n d/2-n/2}(-1)^{n-n'}\sum_{\lambda'\vdash n'} (-1)^{n'-\ell(\lambda')}\Gamma^{St}_{\lambda'}(\gamma'')x_n(d\lambda',q)=(-1)^{n-n'}q^{\delta_{n'}} \sum_{\lambda'\vdash n'} \sigma_{\lambda'}(\gamma'')H^-_{e,d,\lambda'}$$
where $\delta_{n'}:=\frac{(dn'-1)(en'-1)+n'-1}{2}$  In particular, the left-hand sides are also equal up to the same sign.

\end{theorem}
\begin{proof}
Dividing out the sign and using Lemma \ref{lem:coareatoarea}, the LHS reads $$q^{\delta_{n'}}\sum_{\lambda'\vdash n'} (-1)^{\frac{n}{e}-\ell(\lambda')}\Gamma^{St}_{\lambda'}(\gamma'')\prod_{i=1}^{\ell(\lambda')}\left(\sum_{\pi\in\Dyck_{e\lambda_i',d\lambda_i'}^<}q^{-\area(\pi)}h_\pi\right)$$
and we can divide out $q^{\delta_{n'}}$ on both sides. 

On the other hand, by definition we have 
$$H^-_{e,d,\lambda_i'}=\sum_{\alpha\vDash \lambda_i'}\sum_{\substack{\pi \in \Dyck_{e\lambda_i',d\lambda_i'}\\ \touch{\pi}=\alpha}} q^{-\area(\pi)}h_\pi$$
where $\touch(\pi)=\alpha$ specifies that $\pi$ touches the diagonal at $\alpha$.

Given two arbitrary compositions, denote by $\alpha+\beta$ their concatenation. By sorting, this gives a partition of $|\alpha|+|\beta|$ of length $\ell(\alpha)+\ell(\beta)$. Given a collection $\vec{\alpha}$ of compositions $$\alpha^{(1)}\vDash \lambda_1',\ldots,\alpha^{\ell(\lambda')}\vDash \lambda'_{\ell(\lambda')}$$ refining the parts of $\lambda'$,  we write $\vec{\alpha}\leftrightarrow \lambda'$.  

We now further expand the LHS and collect terms as follows. Note that from the equation $\sum_{k=1}^n(-1)^k h_{n-k}e_k=0$ it follows that 
$$e_n=\sum_{\alpha \vDash n}(-1)^{\ell(\alpha)} h_\alpha$$ where we sum over all compositions of $n$.
Writing $\lambda'=(\lambda_1',\ldots,\lambda_{\ell(\lambda')}')$ we then have 
\begin{equation}
\label{eq:altsum}
e_{\lambda'}=\prod_{i=1}^{\ell(\lambda')}(\sum_{\alpha \vDash \lambda_i'} (-1)^{\ell(\alpha)}h_\alpha)
\end{equation}
By definition of the Dyck germs of $\mathbf{f}_{\gamma''}$, $$\sum_{\lambda'\vdash n/e} \Gamma^{St}_{\lambda'}(\gamma'')e_{\lambda'}=\sum_{\lambda'\vdash n/e} \sigma_{\lambda'}(\gamma'')h_{\lambda'}$$
Fix $\mu'\vdash n/e$. Plugging Eq. \eqref{eq:altsum} in the expression
$$\sum_{\lambda'\vdash n/e} \Gamma^{St}_{\lambda'}(\gamma'')e_{\lambda'}$$
and collecting all collections of compositions whose sum has associated partition $\mu'$ we see that
$$\sigma_{\mu'}(\gamma'')=\sum_{\substack{\vec{\alpha}\leftrightarrow \lambda'\vdash n/e\\\text{sort}(\alpha^{(1)}+\cdots+\alpha^{(k)})=\mu'}}
(-1)^{\frac{n}{e}-\ell(\lambda')}\Gamma^{St}_{\lambda'}(\gamma'')$$ where the sum runs over all $\lambda'\vdash n/e$ and all collections of compositions $\vec{\alpha}$ refining the parts of $\lambda'$, such that the partition given by adding the compositions and sorting is exactly $\mu'$.
It now remains to replace $H^-_{d,e,\lambda'}$ by a similar expansion.
By definition we have 
$$
H^-_{d,e,\lambda'}=\prod_{i=1}^{\ell(\lambda')}
(\sum_{\alpha\vDash \lambda_i'}\sum_{\substack{\pi \in \Dyck_{e\lambda_i',d\lambda_i'}\\ \touch{\pi}=\alpha}} q^{-\area(\pi)}h_\pi)
$$
Picking one of the summands over $\alpha$, we notice that  $$\sum_{\substack{\pi \in \Dyck_{e\lambda_i',d\lambda_i'}\\ \touch{\pi}=\alpha}} q^{-\area(\pi)}h_\pi=
\prod_{i=1}^{\ell(\alpha)}H^{-,<}_{d,e,\alpha_i}$$ 
where we use the fact that the area statistic is additive on concatenation of Dyck paths (but note the coarea is not). Here $H^{-,<}_{d,e,\alpha_i}$ is defined as above.
Again collecting all $\vec{\alpha}\leftrightarrow \mu'$ we see that 
the terms contributing to $H_{d,e,\mu'}^{-}$ on the right are of the form
$$(-1)^{\frac{n}{e}-\ell(\lambda')}\Gamma^{St}_{\lambda'}(\gamma'') \prod_{k=1}^{\ell(\lambda')}\prod_{i=1}^{\ell(\alpha^{(k)})}H^{-,<}_{d,e,\alpha^{(k)}_i}$$
where $\vec{\alpha}$ refines parts of $\lambda'$ and sums and sorts to $\mu'$. Summing over all such collections
we get the desired result.
\end{proof}

We have the following corollary, which is a form of Theorem \ref{thm:mastersymfn} from the introduction. 
\begin{corollary}
\label{cor:inertiallyelliptic}
Let $\gamma\in GL_n(\cO)$ be inertially elliptic and tamely ramified, with Newton pairs $(\vec{p},\vec{q})$. 
Consider the orbital symmetric function $\fgamma$ for $\gamma$ as defined in Definition \ref{def:ellipticmastersymfn} and the annular symmetric function $\fpq$ introduced in Definition \ref{def:degmaster}. Then $\fgamma=q^{\Xi(\gamma)} \omega \fpq|_{q\mapsto q^{-1}}$. Here $\Xi(\gamma)$ is as defined in Definition \ref{def:RV}.
\end{corollary}
\begin{proof}
    By assumption $\gamma$ lies in the ring of integers of the totally and tamely ramified extension $F(\gamma)/F$. Dividing $\gamma$ by an element in $\cO_F^{\times}$, we may assume $\gamma$ is topologically unipotent, i.e. $\gamma\equiv1(\operatorname{mod} \mathfrak{m}_{F(\gamma)})$. The resulting $\gamma$ is of the type Lemma \ref{lem:factoreta} and Proposition \ref{prop:wrecursion} apply to.
    
    Write $\gamma=\eta(1+\delta\gamma')$ as in Lemma \ref{lem:factoreta}.  By Theorem \ref{thm:comparerecursion} combined with the $f=1$ case of Proposition \ref{prop:unramifiedextension} we see that the orbital symmetric functions of $\gamma\in G(F)$ and $\gamma'\in Z_{G(F)}(\delta)$ are related by the equation $$\fgamma=q^{\frac{(dn'-1)(en'-1)+n'-1}{2}}(-1)^{n-n'}\omega\varphi_{d/e}|_{q\mapsto q^{-1}}\omega^{-1}((-1)^{n-n'}\tau_1(\mathbf{f}_{\gamma'}))$$
    where $d, e$ are as defined after Lemma \ref{lem:factoreta}.

    First, we note the overall sign cancels. 
    Second, to see that the power of $q$ adds up to $\Xi(\gamma)$, note that the exponent for each step is $\frac{1}{2}((dn'-1)(n-1)+n'-1)=\frac{1}{2}(dn'(n-1)+n'-n)$. On the other hand, as both $\gamma$ and $\gamma'$ are inertially elliptic, by Lemma \ref{lem:reductionnumerics}(2), $\Xi(\gamma)-\Xi(\gamma')$ is equal to $\frac{1}{2}(\frac{d}{e}n(n-1)-(n-n'))=\frac{1}{2}(dn'(n-1)+n'-n)$.
    By Lemma \ref{lem:reductionnumerics}(3), this proves the corollary inductively.
\end{proof}

    We note that in Theorem \ref{thm:mastersymfn} we worked with topologically nilpotent elements, which can be taken to be just $\gamma-1$ where $\gamma$ is in Corollary \ref{cor:inertiallyelliptic}. The Newton pairs associated to $\gamma$ (see Definition \ref{def:puiseux}) is obviously equal to the Newton pairs associated to $\gamma-1$.

\begin{remark}
    Lemma \ref{lem:factoreta} could have worked without $\eta$ but with the price that $F'/F$ might be a trivial extension. In fact, one can run the recursion without choosing an $\eta$ at each step
    and this still gives us a sequence of pairs of integers $(\vec{p},\vec{q})$ and hence an associated annular symmetric function $\frs$. See e.g. Example \ref{ex:reductioncase} below for an example. That the result is independent of which way we proceed is clear from that $\Gamma_{\lambda}(\gamma)=\Gamma_{\lambda}(\gamma/\eta)$ for $\eta\in\cO_F^{\times}$. One may ask whether this is clear from the combinatorial setting of Section \ref{sec:knots} where one starts simply with a sequence $(\vec{p},\vec{q})$. Indeed, the relevant symmetry can be deduced for the refined annular symmetric functions $\fqrs$ and hence for the $\fpq$ as well, using Cherednik-Danilenko's ''reduction cases" in \cite[(4.4)--(4.5)]{CD2}. 
\end{remark}
\begin{example}
\label{ex:reductioncase}
As a continuation of Example \ref{ex:totram} and the above Remark, we write the recursion of Corollary \ref{cor:inertiallyelliptic} in two different ways. Let $\gamma=1+u^4+u^6+u^7$. Proceeding as in the Corollary, we get 
$$
\gamma=1+u^4+u^6+u^7=(1+u^4)(1+u^6(1+u-u^4-u^5+\cdots))
$$
and further $\gamma'=1+u-u^4-u^5+\cdots=1+u(1-u^3-u^4+\cdots)$. This shows that the orbital symmetric function is 
$$\fgamma=q^8\omega\varphi_{3/2}(\varphi_{1/2}(e_1))|_{q\mapsto q^{-1}}$$ which is also written out in Example \ref{ex:twothirteen2}.
If we were to use Lemma \ref{lem:factoreta} but without $\eta$, we would get $\gamma=1+u^4(1+u^2+u^3)$ and 
$\gamma'=1+u^2+u^3=1+u^2(1+u)$. This gives 
$$\fgamma=q^8\omega\varphi_{1/1}(\varphi_{1/2}(\varphi_{1/2}(e_1)))|_{q\mapsto q^{-1}}$$ We leave it to the reader to verify these two are the same symmetric function.
\end{example}

When $\gamma$ is not inertially elliptic, we still have a finite algorithm to compute the orbital symmetric function $\fgamma$. From Lemma \ref{lem:factoreta} we have the following.
\begin{proposition}
\label{prop:discreteinvariants}
    Take a topologically unipotent elliptic tamely ramified $\gamma\in GL_n(F)$. There exist
    \begin{enumerate}
        \item A sequence of triples of integers $(f_i, q_i, p_i), i=1,\ldots, k$ with $(q_i,p_i)=1$,
        \item A sequence of fields $F=F_k\subsetneq F_{k-1} \subsetneq \cdots\subsetneq F_1\subsetneq F_0=F(\gamma)$, and
        \item A sequence of elements $\gamma_1,\ldots, \gamma_{k}=\gamma\in F_0^{\times}$
    \end{enumerate}
    such that for $i=1,...,k$ we have
    \begin{enumerate}
\renewcommand{\labelenumi}{(\alph{enumi})}
        \item Each $\gamma_i$ is topologically unipotent, i.e. $\gamma_{i}\in 1+\mathfrak{m}_{F(\gamma_i)}$.
        \item There exists $\eta_i\in 1+\mathfrak{m}_{F_i}$ and $\delta_i\in F_{i-1}$ such that $\gamma_i=\eta_i(1+\delta_i\gamma_{i-1})$ and that $\delta_i$ is $F_{i-1}/F_i$-cuspidal. In particular $F_{i-1}=F_i(\delta_i)$.
        \item  $f_i$ (resp. $p_i$) is the residue degree (resp. ramification degree) of $F_{i-1}/F_i$, and $q_i=\val_{F_{i-1}}(\delta_i)$.
    \end{enumerate}
    Moreover, the sequence of triples $(f_i, q_i, p_i)$ is uniquely determined by $\gamma$.
\end{proposition}
We also get the following generalization of Lemma \ref{lem:reductionnumerics}(3).
\begin{lemma}\label{lem:reductionnumerics2}
Let $\gamma$, $\delta$ and $\gamma'$ be as in Lemma \ref{lem:factoreta}. Let $e$ and $f$ be the ramification index and residual degree of $F'/F$ in the corresponding notations. Write $n':=n/ef$. Then the sequence of triples in Proposition \ref{prop:discreteinvariants} associated to $\gamma'$ is given by deleting the last triple $(f_k,q_k,p_k)$ is the sequence of triples associated to $\gamma$. We have $q_k/p_k=\val_F(\delta)$, $e=p_k$ and $f=f_k$.
\end{lemma}

Denote the residue field of $F_i$ in Proposition \ref{prop:discreteinvariants} by $k_i$. Combining the Proposition \ref{prop:wrecursion}, Proposition \ref{prop:unramifiedextension}, Proposition \ref{prop:discreteinvariants} and Lemma \ref{lem:reductionnumerics2}  we get the following Theorem, which is the most general expression for the orbital symmetric function of a compact, elliptic and tamely ramified $\gamma \in \GL_n(F)_c$. 

\begin{theorem}
\label{thm:mainformula}
Let $\gamma\in GL_n(F)$ be topologically unipotent, elliptic and tamely ramified with discrete invariants  $(\vec{f},\vec{q},\vec{p})$. Then
$$\mathbf{f}_{\gamma_{i}}=|k_i|^{\frac{(dn/e-1)(n-1)+n/e-1}{2}}\omega \varphi_{q_i/p_i}|_{q\mapsto |k_i|^{-1}}\omega^{-1}(\tau_{f_i}(\mathbf{f}_{\gamma_{i-1}}))$$ 
Here, at the $i$-th step, $d=q_i$ and $e=p_i$ are the valuation datum and ramification index, and $n=n_i:=\prod_{j=1}^{i}p_jf_j=[F(\gamma_i):F]$ (so $n/e=n_i/p_i$); the exponent and the base $|k_i|$ refer to these step-$i$ quantities.

For $b\geq 1$ let $\tau'_b:\Sym_q\to\Sym_q$ be the operator sending $p_k\mapsto p_{bk}$ for all $k\geq 1$ and $q\mapsto q^b$.
In particular, we may write 
$$\fgamma=q^{\Xi(\gamma)}\omega \varphi_{q_k/p_k}(\tau'_{f_k}(\varphi_{q_{k-1}/p_{k-1}}(\tau'_{f_{k-1}}(\cdots \varphi_{q_1/p_1}(\tau'_{f_1}(e_1))\cdots)))|_{q\mapsto q^{-1}}$$
\end{theorem}
\begin{proof}
    The proof is essentially the same as for Corollary \ref{cor:inertiallyelliptic}, with the difference that now one also applies the plethysm $\tau_b$ with $b>1$. As the cardinality of the residue field in Proposition \ref{prop:wrecursion} is that of the base field, and we are applying the proposition recursively, we need to raise the variable $q$ to the residue degree at each step. For any symmetric function $g$ of degree $a/b$ with $a,b\in \Z_{\geq 1}$ and $b|a$, we have
    $$\langle \tau_b(g), e_\lambda \rangle =\begin{cases}
        0, & \lambda\neq b \lambda' \text{ for any } \lambda'\vdash a/b\\
        (-1)^{a-a/b} \langle g, e_{\lambda'}\rangle, & \lambda=b\lambda' \text{ for some } \lambda' \vdash a/b
    \end{cases}$$
    and $$\langle \tau_b(g), h_\lambda \rangle =\begin{cases}
        0, & \lambda\neq b \lambda' \text{ for any } \lambda'\vdash a/b\\
        \langle g, h_{\lambda'}\rangle, & \lambda=b\lambda' \text{ for some } \lambda' \vdash a/b
    \end{cases}$$
    In particular, any sequence of operators of the form $\omega \tau'_b \omega$ in the given expression may be replaced by $(-1)^{a-a/b}\tau'_b$ where $a/b$ is the degree of the symmetric function these operators are being applied to.
    This introduces the overall sign $(-1)^{\sum_{i=1}^{k} (\prod_{j=1}^i p_jf_j)-f_i^{-1}(\prod_{j=1}^{i}p_jf_j)}$. 
    However, applying Proposition \ref{prop:unramifiedextension} repeatedly gives the same overall sign, and together they cancel.
    Lastly, the power $q^{\Xi(\gamma)}$ is computed in the same way as in Corollary \ref{cor:inertiallyelliptic}, with Lemma \ref{lem:reductionnumerics2} replacing Lemma \ref{lem:reductionnumerics}(c).
\end{proof}

Just as in the discussion following Corollary \ref{cor:inertiallyelliptic}, note that the first case considered in Theorem \ref{thm:mainformula} applies verbatim to $\gamma+1$ where $\gamma\in \fg\fl_n(F)$  is topologically nilpotent, elliptic, and tamely ramified. In particular, the theorem gives a formula for the Lie algebra $\fgamma$ in Definition \ref{def:ellipticmastersymfn} (3). Note also that more generally than in the Theorem, if $\gamma\in \GL_n(F)$ is compact of valuation $\neq 0$, it is still possible to compute $\fgamma$ using Proposition \ref{prop:unramifiedextension} combined with the Theorem. These cases cover all of the three possibilities in Definition \ref{def:ellipticmastersymfn}.

Motivated by the above and Definition \ref{def:degmaster}, we also define the combinatorial counterpart of $\fgamma$ as above. Namely, given {\em any} sequence of triples of positive integers $(\vec{f},\vec{q},\vec{p})$ we define
\begin{equation}
\label{eq:degenerate_combcounterpart}
    \mathbf{f}_{(\vec{f},\vec{q},\vec{p})}:=\varphi_{q_k/p_k}(\tau_{f_k}'(\varphi_{q_{k-1}/p_{k-1}}(\tau_{f_{k-1}}'(\cdots \varphi_{q_1/p_1}(\tau_{f_1}'(e_1))\cdots)))
\end{equation}

\begin{proposition}[Homogeneity of the orbital symmetric function]\label{prop:homog-orbital}
Let $\nabla$ denote the Bergeron--Garsia operator from Definition \ref{def:nablaoperator}. For any tamely ramified compact $\gamma \in \fg\fl_n(F)$, the following equation holds in $\Sym_q$:
\begin{equation}
\label{eq:homogeneity}
\nabla|_{t=1}\fgamma=\mathbf{f}_{\varpi\gamma},
\end{equation}
where the notation $\nabla|_{t=1}$ refers to specializing the parameter $t$ in $\Sym_{q,t}$ to $1$. Equivalently, $\nabla|_{t=1}$ scales each $\th_\lambda$ in the expansion $\fgamma=\sum_\lambda \Gamma_{\lambda^t}(\gamma)\th_\lambda$ by $q^{n(\lambda)}$, recovering the classical homogeneity property of Shalika germs (Proposition \ref{prop:homogeneity}; see also \cite[Lemme 1.2.]{W1}).
\end{proposition}
\begin{proof}
At $t=1$, the action of $\nabla$ on $\Sym_q$ coincides with the slope $1/1$-plethysm $\varphi_{1/1}$. Combined with Theorem \ref{thm:mainformula}, this corresponds to changing the last Newton pair from $(p_k, q_k)$ to $(p_k, q_k+p_k)$, i.e., to multiplying $\gamma$ by $\varpi$.
\end{proof}

As a logical consequence of Proposition \ref{prop:homog-orbital} together with the formula for unramified extensions, we also get a new proof of \cite[Th\'eor\`eme 1.3.]{W1}, which addresses the following situation. Suppose $F'/F$ is unramified of degree $f$, $X\in \cO_{F'}$ generates the residue field $k'$ of $F'$ and $Y\in\mathfrak{gl}_{n/f}(\cO_{F'})$ is such that $F'(Y)$ is an unramified extension of degree $n/f$. Let $0\le a<b$ be integers and $\gamma:=1+\varpi^aX+\varpi^bY$, $\gamma':=1+\varpi^{b-a}Y$. We have

\begin{theorem}
\label{thm:unramified}
$$\fgamma=|k|^{a(n^2-n)/2}\omega (\nabla|^{a}_{t=1,q=1/|k|}\omega\tau_f(\mathbf{f}_{\gamma'}))$$
where $\nabla$ is the operator from Definition \ref{def:nablaoperator}.
\end{theorem}
\begin{proof}
This follows from the fact that $\nabla|_{t=1}^a=\varphi_{a/1}$ as an operator on symmetric functions and Theorem \ref{thm:mainformula} with $f_1=f, q_1=a, p_1=1$. 
\end{proof}
\begin{remark}
In order to compare the Shalika germs denoted by $s_\lambda(\gamma)$ in \cite{W1} to ours, we notice there is a factor of $c^{Wal}_\lambda(q)$ (Eq. \eqref{eq:ccoeffwal}) and another of 
$c_{\lambda'}^{Wal}(q^f)$ inside the plethysm used in {\em loc. cit.}. This is explained by the fact that there is a mismatch between \cite{W1,W2}, namely our orbital symmetric function $\fgamma$ 
is defined to cohere with the latter paper \cite{W2}, 
whereas in \cite{W1} paper the Shalika germs are collected into a generating function 
$$\sum_\lambda s_{\lambda^t}(\gamma,q)c^{Wal}_\lambda(q) h_\lambda$$ instead of $\fgamma=\sum_\lambda \Gamma_{\lambda^t}(\gamma) \th_\lambda$ which has an additional plethysm $X\mapsto X/(q-1)$. Composing this with $\tau_f$ explains the power $q\mapsto q^f$.

Finally, on the LHS of \cite[Th\'eor\`eme 1.3.]{W1} we have factors of the form $q^{-an(\lambda^t)}$ which are exactly the ones coming from homogeneity of Shalika germs as in Proposition \ref{prop:homog-orbital}. This matches the appearance of $\nabla|_{t=1}^a$ above.
\end{remark}
\begin{remark}
In \cite{W2} an unramified character and some ''twisted" Steinberg germs appear. While these are not studied in the present paper, they also have nice expressions and combinatorics in terms of symmetric functions. For example, the fundamental lemma proved in \cite{W2} can be given meaning in this language. We will do this elsewhere.
\end{remark}

\subsection{A canonical $t$-deformation}
\label{sec:tdeformation}
We now discuss a canonical $t$-deformation of $\fgamma$ as defined above for 
$\gamma\in G(F)$ or $\gamma\in \fg(F)$ that are tamely ramified and elliptic.

Note that by induction, as explained in Theorem \ref{thm:mainformula}, $\fgamma$ is constructed using the steps in Theorem \ref{thm:comparerecursion} as well as Proposition \ref{prop:unramifiedextension} (or Theorem \ref{thm:unramified}), which are operations on symmetric functions, namely compositions of slope $m/n$ plethysms $\varphi_{m/n}: \Sym_q\to \Sym_q$, the specialized nabla operator $\nabla|_{t=1}$, scalar multiplication, and the Adams operations $\tau_f$. 
Promoting $\nabla|_{t=1}$ to $\nabla: \Sym_{q,t}\to \Sym_{q,t}$, the slope $m/n$-plethysms $\varphi_{m/n}$ to a family of operators coming from the elliptic Hall algebra via $\widehat{\varphi}_{m/n}:\Sym_{q,t}\to \cE^{m/n}$, and keeping the $\tau_f$ as they are, we may run the similar recursion which only depends on the datum of $\gamma$.
In particular, we define 
\begin{definition}
\label{def:fullfgamma}
Let $\gamma\in \fg\fl_n$ be tamely ramified, topologically nilpotent and elliptic. Let $(\vec{f},\vec{p},\vec{q})$ be the discrete invariants attached to $\gamma$ by Proposition \ref{prop:discreteinvariants}. 

Eq. \eqref{eq:degenerate_combcounterpart} deforms to involve a $t$ as explained above, and with this motivation we define the {\em refined annular symmetric function} of $(\vec{f},\vec{p},\vec{q})$  as $$
\widehat{\mathbf{f}}_{(\vec{f},\vec{p},\vec{q})}:=
\widehat{\varphi}_{q_k/p_k}(\tau_{f_k}'(\cdots\tau_{f_3}'(\widehat{\varphi}_{q_{2}/p_{2}}(\tau_{f_2}'(\widehat{\varphi}_{q_1/p_1}(\tau_{f_1}'(e_1))\cdot 1))\cdot 1))\cdots)\cdot 1$$
Here $k$ is the number of steps in Proposition \ref{prop:discreteinvariants}, each $\widehat{\varphi}_{q_j/p_j}(g)\in\cE^{q_j/p_j}$ is an operator on $\Sym_{q,t}$, and the symbol ``$\cdot 1$'' denotes evaluating that operator on $1\in\Sym_{q,t}$ before the next level is applied.
\end{definition}

Note that in the totally ramified case we recover the refined annular symmetric function $\fqpq$ from Section \ref{sec:knots}.
This symmetric function, while carrying all the information and behaving nicely with respect to combinatorics, is again not the direct generalization of the functions $\fgamma$, as we need to twist by $\omega$.
Motivated by this, we define 
\begin{equation}
\label{eq:deffgammahat}
    \widehat{\mathbf{f}}_{\gamma}:=(qt)^{\Xi(\gamma)}\omega \widehat{\mathbf{f}}_{(\vec{f},\vec{p},\vec{q})}|_{q\mapsto q^{-1},t\mapsto t^{-1}}
\end{equation}
\begin{example}
    With these conventions, we for example have 
    $$\widehat{\mathbf{f}}_{\gamma}=qt\omega\nabla e_2$$
    in the case $\gamma=u^3\in \fg\fl_2(F)$.
\end{example}
In general we have proved the following.
\begin{theorem}
\label{thm:tdeformation}
Let $\gamma\in \fg(F)$ be compact, tamely ramified and elliptic. Then 
the orbital symmetric function admits a canonical $t$-deformation, namely 
$$\fqgamma=\sum_\lambda \widetilde{\Gamma}_{\lambda^t}(\gamma)\tH_\lambda$$
where $\tH_{\lambda}$ are the modified Macdonald polynomials. In particular, the Shalika germs $\Gamma_\lambda(\gamma)$ admit a canonical $t$-deformation.
\end{theorem}

\subsection{Transition matrices for Shalika germs}
\label{sec:transitionmatrices}
In this section we state and prove the main combinatorial formula for the transition matrices of Shalika germs. Throughout, $\gamma\in\fgl_n(F)$ is elliptic, topologically nilpotent and tamely ramified, and
$$\fgamma=\sum_{\lambda\vdash n} \Gamma_{\lambda^t}(\gamma)\,\th_\lambda$$
is the Shalika expansion of its orbital symmetric function, so that the germ $\Gamma_{\lambda^t}(\gamma)$ is the coefficient of $\th_\lambda$.

Recall from Theorem \ref{thm:mainformula} that $\fgamma$ is assembled along a tower of extensions. The ``cabling'' step passing from $\mathbf{f}_{\gamma_{i-1}}$ to $\mathbf{f}_{\gamma_{i}}$, with new Newton pair $(p_i,q_i)$, expands $\mathbf{f}_{\gamma_{i-1}}$ in the $\{h_\lambda\}$ and replaces each $h_\lambda$ by $H^-_{q_i,p_i,\lambda}$; when $f_i>1$ it is preceded by the unramified twist $\tau_{f_i}$, which is the ``slope $0/f_i$'' case of the same process. On the level of Shalika expansions we treat the ramified and unramified steps separately, so that each step of the recursion of Theorem \ref{thm:mainformula} becomes a single change of basis.

By Definition \ref{def:ellipticmastersymfn}, the transition matrix along one such step is precisely the matrix of the operator $\omega\, \varphi_{q_i/p_i}|_{q\mapsto q^{-1}}\, \omega$ (in the ramified case) or $\tau_{f_i}$ (in the unramified case), written in the bases $\{\th_\lambda\}_{\lambda\vdash n}$ and $\{\th_{\lambda'}\}_{\lambda'\vdash n/e}$. In either case we write
\begin{equation}
    M^{d/e}=\{M^{d/e}_{\lambda,\lambda'}\},
\end{equation}
with the convention that $M^{0/e}$ is the matrix of $\tau_e$. Concretely, writing $\gamma'\in\fgl_{n/e}(F')$ for the element preceding $\gamma$ in the recursion --- so that $F'/F$ is totally tamely ramified of degree $e$ when $d\ge1$, and unramified of degree $e$ when $d=0$ --- the matrix $M^{d/e}$ relates their germs by
\begin{equation}
\label{eq:germtransition}
\Gamma_{\lambda^t}(\gamma)=\sum_{\lambda'\vdash n/e} M^{d/e}_{\lambda,\lambda'}\,\Gamma_{{\lambda'}^t}(\gamma').
\end{equation}
\begin{theorem}
\label{thm:transitionmatrix}
Keep the notation above: $\gamma$ is elliptic, topologically nilpotent and tamely ramified, $\gamma'$ is the preceding element in the recursion, and $M^{d/e}$ is the transition matrix of \eqref{eq:germtransition}, with $e\geq 1$ and either $d=0$ or $d\geq 1$ with $(e,d)=1$. Considering the symmetric group of permutations on $n/e$ letters, denote by $|S_{\lambda'}\cap C_\mu|$ the number of permutations simultaneously lying in the Young (parabolic) subgroup $S_{\lambda'}:=S_{\lambda'_1}\times\cdots\times S_{\lambda'_{\ell}}$ and in the conjugacy class $C_\mu$ of permutations with cycle type $\mu$.

When $d\geq 1$, the transition matrix $M^{d/e}$ of Shalika germs
has a combinatorial description as follows:
\begin{equation}
\label{eq:transitionmatrix}
M_{\lambda,\lambda'}^{d/e}=\left(c_{\lambda'}|_{q\mapsto q^{-1}}\; q^{s} \sum_{\mu\vdash n/e} \frac{|S_{\lambda'}\cap C_\mu|}{b_\mu|_{q\mapsto q^{-1}}\, \lambda'!}\prod_{i=1}^{\ell(\mu)}(\sum_{\alpha\vDash e\mu_i} \wt_{d/e}(\alpha)_{q\mapsto q^{-1}} q^{-n(\alpha^t)}\th_\alpha)\right)\Bigg|_{\th_\lambda}
\end{equation}
where $c_{\lambda'}$ is as in Eq. \eqref{eq:ccoeff}, $b_\mu:=\prod_i (1-q^{\mu_i})$, $s:=n({\lambda'}^t)+\frac{(dn'-1)(en'-1)+n'-1}{2}$ and $\wt(\alpha)_{d/e}$ is defined in Eq. \eqref{eq:compweightfull}. When $d=0$ the transition matrix is given by
\begin{equation}
\label{eq:transitionmatrixunr}
M_{\lambda,\lambda'}^{0/e}=\left(c_{\lambda'}\sum_{\mu\vdash n/e} \frac{|S_{\lambda'}\cap C_\mu|}{b_\mu \lambda'!}\prod_{i=1}^{\ell(\mu)}(\sum_{\alpha\vDash e\mu_i} \wt_{0/e}(\alpha)\th_\alpha)\right)\Bigg|_{\th_\lambda}
\end{equation}
where $\wt_{0/e}(\alpha)$ is still defined by Eq. \eqref{eq:compweightfull} but with $S_{0/e}(i):=0$ for all $i$ when $d=0$. 
\end{theorem}
\begin{proof}
We need to compute the slope ''$d/e$ plethysm" of the functions $\th_{\lambda'}$, i.e. expand them in the $p_{\mu}$ and replace each $p_k$ by $\varphi_{d/e}(p_k)=P_{dk,ek}$ and bring the result back to the basis $\th_\lambda$.  In order to conform to the recipe in Theorem \ref{thm:mainformula} we also need to conjugate the plethysm by $\omega$ and invert $q$, as well as multiply the result by $q^{\frac{(dn'-1)(en'-1)+n'-1}{2}}$. In the case $d=0$ there is no conjugation by $\omega$ and we simply have $\varphi_{0/e}(p_k):=\tau_e(p_k)=p_{ek}$.

We first notice that by \cite{ReEc} the {\em untransformed} complete homogeneous symmetric functions satisfy 
$$h_{\lambda'}=\sum_\mu \frac{|S_{\lambda'}\cap C_\mu|}{\lambda'!}p_\mu$$ where $|S_{\lambda'}\cap C_\mu|$ is as defined above.
Since $$\th_{\lambda'}=c_{\lambda'} h_{\lambda'}\left[\frac{X}{1-q}\right],$$
we get
$$\th_{\lambda'}=c_{\lambda'} \sum_{\mu\vdash n'} \frac{|S_{\lambda'}\cap C_\mu|}{\lambda'! b_\mu} p_\mu$$
where $b_\mu=\prod_i (1-q^{\mu_i})$ is the reciprocal of the principal specialization of $p_\mu$ and $c_{\lambda'}$ is as before.
Let now $d>0$. 

We first note that the matrix elements of $\omega \varphi_{d/e}|_{q\mapsto q^{-1}} \omega$ in the $\th_\lambda$-basis are the same as those of $\varphi_{d/e}$ in the $\omega \th_\lambda$-basis. On the other hand, $\omega\th_\lambda=q^{n(\lambda^t)}\th_\lambda[X;q^{-1}]$ by Eq. \eqref{eq:invertq}. Applying $\varphi_{d/e}|_{q\mapsto q^{-1}}$ to $\th_\lambda[X;q^{-1}]$ is clearly the same as applying $\varphi_{d/e}$ to $\th_\lambda$ and then inverting $q$. By Proposition \ref{prop:limitformulafull} we get $$\varphi_{d/e}(\th_{\lambda'})=\left(c_{\lambda'}\sum_{\mu\vdash n/e} \frac{|S_{\lambda'}\cap C_\mu|}{b_\mu \lambda'!}\prod_{i=1}^{\ell(\mu)}(\sum_{\alpha\vDash e\mu_i} \wt_{d/e}(\alpha)\th_\alpha)\right)$$ and by the above argument applying $\omega$ and $\varphi_{d/e}|_{q\mapsto q^{-1}}$ gives
 $$\varphi_{d/e}|_{q\mapsto q^{-1}}(\omega \th_{\lambda'})=\left(c_{\lambda'}|_{q\mapsto q^{-1}}\,q^{n({\lambda'}^t)}\sum_{\mu\vdash n/e} \frac{|S_{\lambda'}\cap C_\mu|}{b_\mu|_{q\mapsto q^{-1}}\, \lambda'!}\prod_{i=1}^{\ell(\mu)}(\sum_{\alpha\vDash e\mu_i} \wt_{d/e}(\alpha)|_{q\mapsto q^{-1}}\th_\alpha[X;q^{-1}])\right)$$
Applying $\omega$ once more gives
$$\omega \varphi_{d/e}|_{q\mapsto q^{-1}}(\omega \th_{\lambda'})=\left(c_{\lambda'}|_{q\mapsto q^{-1}}\,q^{n({\lambda'}^t)}\sum_{\mu\vdash n/e} \frac{|S_{\lambda'}\cap C_\mu|}{b_\mu|_{q\mapsto q^{-1}}\, \lambda'!}\prod_{i=1}^{\ell(\mu)}(\sum_{\alpha\vDash e\mu_i} \wt_{d/e}(\alpha)|_{q\mapsto q^{-1}}q^{-n(\alpha^t)}\th_\alpha)\right)$$
and multiplying by the prefactor $q^{\frac{(dn'-1)(en'-1)+n'-1}{2}}$ gives the result.

The proof for the $d=0$ case is similar, except we do not need to apply $\omega$ nor invert $q$. Note that Eq. \eqref{eq:compweightfull} still holds in this case by \cite{NegutPieri} and is closely related to the classical Pieri rule for Macdonald polynomials.
\end{proof}
\begin{remark}
One may view Theorem \ref{thm:transitionmatrix} as giving a combinatorial expression for the ''$\lambda'$-colored" orbital symmetric functions of torus knots at $a=0, t=1$.
\end{remark}

\subsubsection{Integrality properties}
In this subsection we describe a combinatorial reorganization of the Shalika germs which, after a suitable renormalization, reveals striking integrality properties. These are not visible from the transition-matrix formulas above, and rest on the following change-of-basis fact.
\begin{proposition}
\label{prop:changebasis}
The symmetric functions $h_\lambda$ expand with $\Z[q]$-coefficients in the basis $h_\lambda\left[\frac{X}{1-q}\right]$, and the symmetric functions $\th_\lambda$ expand with $\Z[q]$-coefficients in the basis $h_\lambda$.
\end{proposition}
\begin{proof}
The bases $\{m_\lambda\},\{h_\lambda\}$ are dual with respect to the Hall inner product, so that $m_\lambda[X(1-q)]$ is the basis dual to $h_\lambda\left[\frac{X}{1-q}\right]$ by standard properties of plethysm.
Therefore, we need to check that 
$$\langle h_\lambda, m_\mu[X(1-q)]\rangle=\langle h_\lambda[X(1-q)],m_\mu\rangle$$ is in $\Z[q]$. The inner product is nonzero only when there exists an integer matrix with row sums $\lambda$ and column sums $\mu$, with only a single nonzero entry in each row (see e.g. \cite[3.17.]{Zelevinsky}). It is also integral by e.g. \cite[p. 52]{HaimanCDM}. The second statement follows from a similar argument.
\end{proof}

Combining this with the integrality of the transition matrix of Steinberg germs --- equivalently, of $\fgamma$ expanded in the $h_\lambda$-basis, which holds by Theorem \ref{thm:comparerecursion} (we write $M(h)$ for the conjugation of $M$ to the appropriate basis, as in \cite{ReEc}) --- shows that the coefficients of $\fgamma$ in the basis $h_\lambda\left[\frac{X}{1-q}\right]$ are integral. Since $\th_\lambda=c_\lambda h_\lambda\left[\frac{X}{1-q}\right]$, these coefficients are exactly $c_\lambda \Gamma_{\lambda^t}(\gamma)$. This motivates the following.

\begin{defproposition}\label{defprop:renormalized}
The {\em renormalized} Shalika germs $$\Gamma^{ren}_{\lambda}(\gamma):=c^{Wal}_{\lambda^t} \Gamma_\lambda(\gamma)=(-1)^nc_{\lambda^t}\Gamma_\lambda(\gamma)$$ are integral, i.e. $\Gamma^{ren}_{\lambda}(\gamma)\in \Z[q]$.
\end{defproposition}
\begin{proof}
    We have $\th_1=h_1$. By Theorem \ref{thm:mainformula} we can get $\fgamma$ up to a sign by multiplying $h_1$ by various $M^{d/e}$ together. This results in an obviously integral expression for $\fgamma$ in the basis $h_\lambda$. By Proposition \ref{prop:changebasis} we get that the expansion in the $h_\lambda\left[\frac{X}{1-q}\right]$-basis is also integral.
\end{proof}
\begin{remark}
The renormalized Shalika germs are not in $\N[q]$ in general, even up to an overall sign. In particular it is easy to find examples for which $c_\lambda \Gamma_\lambda(\gamma)$ has both positive and negative integer coefficients. We use the overall sign $(-1)^n$ to avoid additional bookkeeping below.
\end{remark}

\begin{remark}[Comparison to Shelstad's normalization]\label{rmk:shelstad}
One can also renormalize the Shalika germs in a representation-theoretic way, further dividing by $q^f-1$ and $q^{\Xi(\gamma)+n-1}$, where $f$ is the degree of the maximal unramified subextension of $F(\gamma)/F$. More precisely, by 
\cite[Eq. (3.6)]{Tsai20} one has
$$(q^f-1)^{-1}q^{-\Xi(\gamma)-n+1}(q-1)^n\,\Gamma_{(n)}(\gamma)=1.$$ For the regular unipotent orbit, this normalization coincides with the one used by Shelstad in \cite{Shelstad}. In general, we set 
$$\Gamma^{dW}_{\lambda}(\gamma):=(q^f-1)^{-1}q^{-\Xi(\gamma)-n+1}c^{Wal}_{\lambda^t}\Gamma_{\lambda}(\gamma)=(q^f-1)^{-1}q^{-\Xi(\gamma)-n+1}\,\Gamma^{ren}_\lambda(\gamma).$$
We call this the {\em degenerate Whittaker normalization}, following \cite{MW87,Shelstad}. This normalization is natural from the point of view of automorphic forms and endoscopy, but does not directly yield the positivity properties of our renormalization.
\end{remark}

\begin{example}
    Let $\gamma=u^3$ in $\fg\fl_2$. Then 
    
    \vspace{0.2cm}
    \begin{center}
    \begin{tabular}{|c|c|c|}
    \hline
 & (11) & (2) \\ \hline
                Ordinary  & $\frac{-1}{q - 1}$ & $\frac{q^2}{q - 1}$\\\hline
        Renormalized & $1-q^2$ & $q^3-q^2$   \\ \hline
        Degenerate Whittaker&  $-q^{-1}-q^{-2}$ & 1\\ \hline
    \end{tabular}
\end{center}
\vspace{0.2cm}

Similarly, if 
$\gamma=u^4\in \fg\fl_3$ we have

\vspace{0.2cm}
    \begin{center}
    \begin{tabular}{|c|c|c|c|}
    \hline
 & (111) & (21) & (3) \\ \hline
                Ordinary  & $\frac{-1}{-q^{3} + q^{2} + q - 1}$ & $\frac{2 q^{4} + q^{3}}{-q^{3} + q^{2} + q - 1}$ &  $\frac{-q^{5}}{-q^{2} + 2 q - 1}$ \\\hline
        Renormalized & $q^{3} - 1$  & $-2 q^{5} + q^{4} + q^{3}$
 & $q^{6} - q^{5}$  \\ \hline
        Degenerate Whittaker& $q^{-3}+q^{-4}+q^{-5}$  & $-2q^{-1}-q^{-2}$ & $1$ \\ \hline
    \end{tabular}
\end{center}
\vspace{0.2cm}

In the latter example, one can directly compare the subregular germ to the formula in \cite[(10.2)]{Repka}: the factor $2+q^{-1}=-q\,\Gamma^{dW}_{(21)}$ is exactly the $r\equiv 1\ (\mathrm{mod}\ 3)$ case there, and the remaining discriminant factors stem from a normalization difference, just as in the comparison between the regular germs in \cite{Shelstad} and \cite{Repka1}. More generally, specializing Theorem~\ref{thm:mainformula} to the regular and subregular orbits reproves the germ computations of \cite{Repka1,Repka} for tamely ramified elliptic $\gamma$. The regular germ is $\Gamma^{dW}_{(n)}\equiv 1$, which is the leading term of \cite{Repka1} in the normalization of \cite{Shelstad}. For the subregular orbit one finds $\Gamma^{dW}_{(n-1,1)}=-\sum_{j=1}^{n}q^{-S_{m/n}(j)}$ when $\gamma=u^m$ is totally ramified of slope $m/n$, and $\Gamma^{dW}_{(2,1)}=-3q^{-a}$ for the unramified cubic of depth $a$ (obtained through the operator $\tau_f$); each agrees with \cite[\S10]{Repka} after multiplication by $q^{-\Xi(\gamma)}=q^{-\dim\Sp_\gamma}$, the factor passing from the geometric to the degenerate-Whittaker normalization. As this derivation runs through the recursion of \cite{W2} and the elliptic Hall algebra action rather than a direct computation of orbital integrals, it constitutes a new proof of Repka's results in the tamely ramified case.
\end{example}
We now turn to the finer behavior of the germs $\Gamma^{dW}_\lambda$ and $\Gamma^{ren}_\lambda$. Computations suggest the following.
\begin{conjecture}
\label{conj:dWintegrality}
    The $\Gamma^{dW}_{\lambda}(\gamma)$ are integral polynomials in $q^{-1}$, i.e. lie in $\Z[q^{-1}]$.
\end{conjecture}
We know of a combinatorial proof of this conjecture for inertially elliptic $\gamma$, based on the graph formula below, which will appear in future work.

We now give a combinatorial description of the degenerate-Whittaker germs $\Gamma^{dW}_\lambda$ of Remark \ref{rmk:shelstad}, conjecturally integral in $q^{-1}$. Let $\lambda,\mu \vdash n\geq 1$, and let $\cG(\lambda)$ be the set of directed graphs (loops allowed) with vertex set the boxes of the Ferrers diagram of $\lambda$, labeled by $\{1,\ldots, n\}$, with edges only between boxes in the same row, and with every vertex of in- and out-degree $1$. Let $\cG(\lambda,\mu)\subset \cG(\lambda)$ be the subset of graphs whose connected components sort to the partition $\mu$; note that $\mu$ is necessarily a {\em refinement} of $\lambda$.
\begin{lemma} 
There is a natural bijection
$S_\lambda\cap C_\mu\leftrightarrow \cG(\lambda,\mu)$.
\end{lemma}
\begin{proof}
Writing a cycle decomposition for elements on the left gives rise to a graph by drawing the boxes labeled $1,\ldots, n$ and adding edges $a_i\to a_{i+1}$ for each cycle $(a_1\cdots a_k)$. The converse is clear.
\end{proof}
Next, we note that by deleting at least one edge from each cycle of a graph $\mathtt{G}\in \cG(\lambda)$, we get a composition of $n$, by remembering the ordering on the original boxes of $\lambda$. This composition naturally refines $\lambda$. Accordingly for $\mathtt{G}\in \cG(\lambda)$, we say $\alpha\vDash n$ refines $\mathtt{G}$ if we can obtain the composition $\alpha$ by deleting edges from $\mathtt{G}$. Finally, for $e\geq 1$, let $e\mathtt{G}$ be the graph obtained by $e$-dilating each cycle in $\mathtt{G}$.

We may further associate to each $\mathtt{G}'\in\cG(\lambda',\mu)$ and a composition $\alpha \leq e \mathtt{G}'$ exactly $\prod_{i=1}^{\ell(\mu)}\mu_i$ different graphs by cyclic permutation of vertices in $\mathtt{G}'$. It is easy to see these graphs $\mathtt{G}\leq e \mathtt{G}'$ are the ones coming exactly from $e\mathtt{G}'$ by removal of one or more edges so that the resulting composition is $\alpha$. 

Next, define the weight of a graph. For each vertex $v$, let
\begin{itemize}
\item[(a)] $\delta_{e\mathtt{G}'}(v):=\lceil j\,d/e\rceil-\lceil (j-1)\,d/e\rceil=S_{d/e}(j)$, where $j$ is the position of $v$ in its circuit of the $e$-dilated graph $e\mathtt{G}'$ (for any fixed origin of the circuit; the final sum is independent of this choice),
\item[(b)] $\coarm_{\mathtt{G}}(v):=\chi(v)-1$, where $\chi(v)$ is the position of $v$ in the chain of $\mathtt{G}$ it belongs to after edges are deleted.
\end{itemize}
The weight is then $$\wt(\mathtt{G})_{d/e}=q^{-\sum_{v}\delta_{e\mathtt{G}'}(v)\,\coarm_{\mathtt{G}}(v)}.$$
\begin{theorem}\label{thm:graphformula}
The transition matrix ${}^{dW}M^{d/e}_{\lambda,\lambda'}$ of the degenerate-Whittaker Shalika germs $\Gamma^{dW}$ is given by
\begin{equation}
\label{eq:transitionmatrix2}
{}^{dW}M_{\lambda,\lambda'}^{d/e}=(-1)^{n-\ell(\lambda)}\frac{1}{\lambda'!}\sum_{\mathtt{G}'\in\cG(\lambda')} (-1)^{n/e-\ell(\alpha(\mathtt{G'}))}\sum_{\substack{\mathtt{G}\leq e\mathtt{G}'\\ \mathrm{sort}(\mathtt{G})=\lambda}}\wt(\mathtt{G})_{d/e}
\end{equation}
\end{theorem}
This is a purely combinatorial statement: it identifies the transition matrix of Eq.~\eqref{eq:transitionmatrix}, converted to the degenerate-Whittaker normalization via Remark \ref{rmk:shelstad}, with a weighted sum over dilated graphs. The proof is a direct, if somewhat lengthy, comparison of the two sides, and will appear separately.

\begin{remark}[The renormalized version]\label{rmk:rengraph}
By Remark \ref{rmk:shelstad}, the two normalizations differ by the factor $(q^f-1)$ as well as by a monomial:
$$\Gamma^{ren}_\lambda(\gamma)=(q^f-1)\,q^{\,\Xi(\gamma)+n-1}\,\Gamma^{dW}_\lambda(\gamma).$$
On the level of transition matrices,
$${}^{ren}M^{d/e}=q^{\,(\Xi(\gamma)-\Xi(\gamma'))+(n-n')}\;{}^{dW}\!M^{d/e}=q^{\,n-n'+\Delta}\;{}^{dW}\!M^{d/e}.$$
\end{remark}

\subsection{The formulas for orbital integrals}
In this section, we give a combinatorial reformulation of the orbital integrals of characteristic functions of standard parahorics, in particular we prove Theorem \ref{thm:weightpolynomials} from the introduction.

\begin{theorem}
\label{thm:orbintformula}
Let $\gamma\in \fg(F)$ be compact, tamely ramified and regular semisimple, and let $\1_\lambda$ be the characteristic function of the standard parahoric $\bP_\lambda$ associated to $\lambda\vdash n$, divided out by its measure (with the normalization of measures as before). Then
$$I_\gamma(\1_\lambda)=\langle \fgamma,h_\lambda\rangle$$ where we pair using the Hall inner product and $\fgamma$ is as above.
\end{theorem}
\begin{proof}
From Theorem \ref{thm:steinberg}, 
we have $$I_\gamma(\1_\lambda)=\sum_\mu \Gamma^{St}_\mu(\gamma) \St_{\mu,c}(\1_\lambda)$$ and by Theorem \ref{thm:comparerecursion} plus Definition \ref{def:ellipticmastersymfn} we have 
$$\sum_\mu \Gamma^{St}_\mu(\gamma) e_\mu=\sum_\mu \Gamma_\mu(\gamma) \th_\mu=\sum_{\mu} \sigma_\mu(\gamma)h_\mu=\fgamma$$
The result then follows from Propositions \ref{prop:wpair1}, \ref{prop:wpair2} and Proposition \ref{prop:germcomparison}.
\end{proof}

\begin{corollary}
\label{cor:sphericalorbital}
Let $\gamma$ be inertially elliptic and tamely ramified, with Newton pairs $(\vec p,\vec q)$, and let $\1_{(n)}$ be the (normalized) characteristic function of the hyperspecial maximal parahoric $GL_n(\cO)$, i.e.\ the spherical test function. Then the spherical orbital integral is the reversed area generating function of the $k$-nested Dyck paths of Definition \ref{def:nesteddyck}:
$$I_\gamma(\1_{(n)})=\langle\fgamma,h_n\rangle=q^{\Xi(\gamma)}\sum_{\Pi\in\nDyck_{(\vec p,\vec q)}}q^{-\area(\Pi)}.$$
Equivalently, it equals $q^{\dim\Sp_\gamma}\langle\fpq,e_n\rangle|_{q\mapsto q^{-1}}$, the weight polynomial of the spherical affine Springer fiber.
\end{corollary}
\begin{proof}
Since $h_n=\sum_{\mu\vdash n}m_\mu$, we have $\langle h_\mu,h_n\rangle=1$ for every $\mu\vdash n$; hence by Theorem \ref{thm:orbintformula}, $I_\gamma(\1_{(n)})=\langle\fgamma,h_n\rangle=\sum_{\mu}\sigma_\mu(\gamma)$ is the sum of all Dyck germs. By Corollary \ref{cor:inertiallyelliptic}, Proposition \ref{prop:annularformula} and $\omega e_\mu=h_\mu$,
$$\fgamma=q^{\Xi(\gamma)}\,\omega\,\fpq|_{q\mapsto q^{-1}}=q^{\Xi(\gamma)}\sum_{\Pi\in\nDyck_{(\vec p,\vec q)}}q^{-\area(\Pi)}\,h_{\lambda(\Pi)},$$
and summing the coefficients gives the formula. The equivalent form follows since $\omega$ is a Hall isometry with $\omega h_n=e_n$, and $\Xi(\gamma)=\dim\Sp_\gamma$ by Remark \ref{rmk:dimspgamma}.
\end{proof}

The same computation evaluates every parahoric orbital integral, not just the spherical one, once one records how much each nested Dyck path contributes to the coefficient of $h_\lambda$.
\begin{corollary}
\label{cor:orbitalnesteddyck}
Let $\gamma$ be inertially elliptic and tamely ramified, with Newton pairs $(\vec p,\vec q)$ and $n=p_1\cdots p_k$, and let $\lambda\vdash n$. Then
\begin{equation}
\label{eq:orbitalnesteddyck}
I_\gamma(\1_\lambda)=\langle\fgamma,h_\lambda\rangle=q^{\Xi(\gamma)}\sum_{\Pi\in\nDyck_{(\vec p,\vec q)}}q^{-\area(\Pi)}\,\bigl|\MM(\lambda,\lambda(\Pi))\bigr|,
\end{equation}
where $\lambda(\Pi)\vdash n$ is the partition given by "top-level" horizontal runs of the nested Dyck path $\Pi$ (Definition \ref{def:nesteddyck}) and $\MM(\lambda,\lambda(\Pi))$ is the set of non-negative integer matrices with row sums $\lambda$ and column sums $\lambda(\Pi)$ (Definition \ref{def:rowcolmatrices}). In particular $I_\gamma(\1_\lambda)$ is a polynomial in $q$ with non-negative integer coefficients.
\end{corollary}
\begin{proof}
By Corollary \ref{cor:inertiallyelliptic}, Proposition \ref{prop:annularformula} and $\omega e_\mu=h_\mu$, as in the previous proof,
$$\fgamma=q^{\Xi(\gamma)}\sum_{\Pi\in\nDyck_{(\vec p,\vec q)}}q^{-\area(\Pi)}\,h_{\lambda(\Pi)}.$$
Pairing with $h_\lambda$ using Theorem \ref{thm:weightpolynomials} (equivalently Theorem \ref{thm:orbintformula}) and the identity $\langle h_\lambda,h_\mu\rangle=|\MM(\lambda,\mu)|$ gives \eqref{eq:orbitalnesteddyck}. The pairing identity is classical: $\langle h_\lambda,h_\mu\rangle$ is the number of ways to fill a matrix with non-negative integer entries and the prescribed row and column sums, since $h_\mu=\sum_\nu\langle h_\mu,h_\nu\rangle m_\nu$ and $\langle h_\lambda,m_\nu\rangle=\delta_{\lambda\nu}$ expand $\prod_j h_{\mu_j}=\prod_j\sum_{|\alpha^{(j)}|=\mu_j}m_{\alpha^{(j)}}$ into monomials. Non-negativity is then manifest, as $q^{\Xi(\gamma)-\area(\Pi)}$ has non-negative exponent (the area of a nested Dyck path is at most $\Xi(\gamma)$, with equality for the unique path of maximal area) and $|\MM(\lambda,\lambda(\Pi))|\in\Z_{\ge0}$.
\end{proof}

The spherical case $\lambda=(n)$ of Corollary \ref{cor:orbitalnesteddyck} is Corollary \ref{cor:sphericalorbital}, since $|\MM((n),\mu)|=1$ for all $\mu\vdash n$; the Iwahori case $\lambda=(1^n)$ has $|\MM((1^n),\mu)|=\binom{n}{\mu}=n!/\prod_i\mu_i!$, the number of set compositions of type $\mu$. Worked instances of both appear in Example \ref{ex:twothirteen2}.

\begin{remark}[Nested rational parking functions]
\label{rmk:nestedparking}
In parallel with nested Dyck paths, it makes sense to speak of $k$-nested parking functions, which are unsurprisingly related to the parahoric versions in \eqref{eq:orbitalnesteddyck}. For a nested Dyck path $\Pi$ with top-level run partition $\lambda(\Pi)$, the number $|\MM(\lambda,\lambda(\Pi))|=\langle h_\lambda, h_{\lambda(\Pi)}\rangle$ is exactly the number of ways to label the $n$ top-level boxes of $\Pi$ with "cars'' of content $\lambda$, using the label $i$ with multiplicity $\lambda_i$, where the labels within a single top run are weakly increading. We call such a labeled path $(\Pi,L)$ a {\em nested rational parking function of content $\lambda$}, and Corollary \ref{cor:orbitalnesteddyck} becomes
$$I_\gamma(\1_\lambda)=q^{\Xi(\gamma)}\sum_{(\Pi,L)}q^{-\area(\Pi)},$$
the sum ranging over nested rational parking functions of content $\lambda$ for $(\vec p,\vec q)$. The two extremes are the familiar ones: content $(n)$ gives a single labeling per path (all cars equal), recovering the unlabeled count of Corollary \ref{cor:sphericalorbital}, while content $(1^n)$ gives the {\em standard} nested rational parking functions, in which the cars $1,\ldots,n$ are distinct. In particular the Iwahori orbital integral $I_\gamma(\1_{(1^n)})$ is the reversed-area enumerator of standard nested rational parking functions.

\end{remark}

As a corollary of the proof, we get 
\begin{corollary}
\label{cor:mastersymfnmanycomponents}
For $\gamma \in \fm(F)\subset \fg(F)$ as in Theorem \ref{thm:orbintformula}, where $\fm$ is the Lie algebra of a Levi subgroup conjugate to $L(\mu)$, let $\gamma_1,\ldots,\gamma_{\ell(\mu)}$ be the ''diagonal blocks" of $\gamma$. Then we have the following identity of orbital symmetric functions:
$$\fgamma=\left|\det(\ad(\gamma)_{\fg/\fm})\right|^{-1/2} \prod_{i=1}^{\ell(\mu)} \mathbf{f}_{\gamma_i}$$
where $\mathbf{f}_{\gamma_i}$ is the orbital symmetric function of $\gamma_i$ as defined in Definition \ref{def:ellipticmastersymfn}.
\end{corollary}
\begin{proof}
As in the statement, suppose $\gamma$ belongs to a Levi of type $\mu$, WLOG to the standard one and has blocks $\gamma_1,\ldots, \gamma_{\ell(\mu)}$. Then by Proposition \ref{prop:reductiontolevi-alg}
$$I^G_\gamma(\1_\lambda)=\left|\det(\ad(\gamma)_{\fg/\fm})\right|^{-1/2} I_\gamma^M(\Res_\fm^\fg(\1_\lambda))$$
Let us write $\Res_\mu=\Res_\fm^\fg$. By \cite[Lemme IV 3.]{W2} and Lemma \ref{lem:definitionsagree}, we get 
$$\Res_\mu(\1_\lambda)=\sum_{m\in \MM(\lambda,\mu)} \otimes_{j=1}^{\ell(\mu)} \1_{m_{\cdot,j}}$$ where $\1_{m_{\cdot,j}}$ is the characteristic function of the corresponding standard parahoric normalized by its measure and $\MM(\lambda,\mu)$ is as in Definition \ref{def:rowcolmatrices}.
It is clear that this implies 
$$I_\gamma^M(\Res_\mu(\1_\lambda))=\sum_{m\in \MM(\lambda,\mu)}\prod_{j=1}^{\ell(\mu)}\langle h_{m_{\cdot,j}},\mathbf{f}_{\gamma_j}\rangle$$
On the other hand, the first displayed equation on \cite[p. 883]{W2} implies that we may write 
the RHS of the above equation as 
$$\langle h_\lambda, \prod \mathbf{f}_{\gamma_i}\rangle$$
Multiplying by $\left|\det(\ad(\gamma)_{\fg/\fm})\right|^{-1/2}$, we are done.
\end{proof}

\section{Examples}
\label{sec:examples}
\begin{example}
\label{ex:twothirteen2}
Let $\gamma=u^{7}+u^6\in \fg\fl_4(F)$ and $\cchar(k)\neq 2$, following Examples \ref{ex:totram} and \ref{ex:reductioncase}. 
Let us write down the annular symmetric function:
$$\mathbf{f}_{(1,2),(3,2)}=\varphi_{3/2}(\varphi_{1/2}(e_1))=\varphi_{3/2}(e_2)=$$
$$=\left(q^{2} + q + 1\right)e_{1,1,1,1} + \left(q^{5} + 2 q^{4} + 4 q^{3} + 2 q^{2} + 2 q\right)e_{2,1,1} + $$ $$\left(q^{6} + q^{4} + q^{2}\right)e_{2,2} + \left(q^{7} + q^{6} + 2 q^{5} + q^{4}\right)e_{3,1} + q^{8}e_{4}$$ Indeed, $2$-nested Dyck paths in this case are simply the 23 Dyck paths in a $6\times 4$ rectangle with these horizontal steps and area statistics. 
We recover $\fgamma$ from Corollary \ref{cor:inertiallyelliptic}, namely 
$$\fgamma=q^8\omega \mathbf{f}_{(1,2),(3,2)}|_{q\mapsto q^{-1}}=\left(q^{8} + q^{7} + q^{6}\right)h_{1,1,1,1} + \left(2 q^{7} + 2 q^{6} + 4 q^{5} + 2 q^{4} + q^{3}\right)h_{2,1,1} $$$$ + \left(q^{6} + q^{4} + q^{2}\right)h_{2,2} + \left(q^{4} + 2 q^{3} + q^{2} + q\right)h_{3,1} + h_{4}$$
The weight polynomial of the spherical affine Springer fiber is (equivalently $\langle\fgamma,h_4\rangle$, by the $\omega$-isometry: $\fgamma=q^{\Xi(\gamma)}\omega\fpq|_{q\mapsto q^{-1}}$ and $\omega h_4=e_4$; the Iwahori line below is likewise $\langle\fgamma,h_{1111}\rangle$)
$$q^{\dim \Sp_\gamma}\langle \mathbf{f}_{(1,2),(3,2)},e_4\rangle|_{q\mapsto q^{-1}}=
1+q+2q^2+3q^3+4q^4+4q^5+4q^6+3q^7+q^8
$$
and that of the Iwahori affine Springer fiber is 
$$q^{\dim \Sp_\gamma}\langle \mathbf{f}_{(1,2),(3,2)},e_{1111}\rangle|_{q\mapsto q^{-1}}=
1+4q+10q^2+20q^3+34q^4+48q^5+54q^6+48q^7+24q^8
$$ 
Note that the first one is just the sum of the coefficients of the various $h_\lambda$ in $\fgamma$. It agrees up to $q\mapsto q^{-1}$ with the computation in \cite[Eq. (3.1)]{CD} -- it seems that there is a typo in that paper, repeating one from Piontkowski's work \cite{Piontkowski}.

Finally, we illustrate Theorem \ref{thm:transitionforSha}, i.e. the combinatorial formula for Shalika germs. On the second step of our induction $n=4, n'=2$. Suppose we want to compute the entry 
$M^{3/2}_{211,2}$ of our transition matrix in the $\th_{\lambda}$-basis. 
We have
$$c_{2}=(1-q)(1-q^2), \; b_{2}=(1-q^2),\; b_{11}=(1-q)^2, \lambda'!=2$$
and therefore
$\th_{2}=\frac{q+1}{2}p_{11}+\frac{1-q}{2}p_{2}$. Now since $(p_2,q_2)=(3,2)$ we must apply the slope $3/2$ plethysm and replace $p_{11}\mapsto P_{3,2}^2, p_{(2)}\mapsto P_{6,4},$ up to conjugation by $\omega$.

Now by formula \eqref{eq:Pmnattequals1} we have
$P_{3,2}=\frac{1}{1-q}\th_{11}-\frac{q}{1-q}\th_{2}
$ so $$P_{3,2}^2=\frac{1}{(1-q)^2}\th_{1111}-\frac{2q}{(1-q)^2}\th_{211}+\frac{q^2}{(1-q)^2}\th_{22}$$
There are $8$ compositions of $4$, and we compute 
$$S_{3/2}(1)=2,S_{3/2}(2)=1,S_{3/2}(3)=2,S_{3/2}(4)=1$$
Plugging this in to Eq. \eqref{eq:compweightfull} gives
$$
\wt(2+1+1)_{3/2}=\frac{q(1+q^2)}{(1-q)^2(1-q^2)},\ \wt(1+2+1)_{3/2}=\frac{2q^2}{(1-q)^2(1-q^2)},\ $$$$ \wt(1+1+2)_{3/2}=\frac{q(1+q)}{(1-q)^3},
$$
so that the coefficient of $\th_{211}$
in $P_{6,4}$ is $$\frac{q(1+q^2)}{(1-q)^2(1-q^2)}+\frac{2q^2}{(1-q)^2(1-q^2)}+\frac{q(1+q)}{(1-q)^3}=\frac{-2q^2 - 2q}{(q- 1)^3}$$
Taken together, we get 
$$\frac{-2q(q+1)}{2(1-q)^2}+\frac{(1-q)(-2q^2 - 2q)}{2(q- 1)^3}=0$$

One can verify that the slope $3/2$ plethysm of $\th_2$ has vanishing coefficient for $\th_{211}$. More generally, one can compute that
$$M^{3/2}=\left(\begin{array}{rrrrr}
\frac{1}{q^{5} - q^{3} - q^{2} + 1} & \frac{-q^{6} - q^{5}}{q^{4} - q^{3} - q + 1} & \frac{q^{7}}{q^{3} - q^{2} - q + 1} & 0 & 0 \\
0 & 0 & \frac{q^{6}}{q^{2} - 2 q + 1} & \frac{-2 q^{8}}{q^{2} - 2 q + 1} & \frac{q^{10}}{q^{2} - 2 q + 1}
\end{array}\right)$$
where the rows are indexed by the partitions $(2),(11)$ and the columns are indexed by $(4),(31),(22),(211),(1111)$.

\end{example}
\begin{example}
The simplest elliptic case with three Puiseux pairs appears in \cite[Eq. (3.8.)]{CD} as well as in \cite{Piontkowski} as an example where previous methods fail. 
This example corresponds to the plane curve singularity $\C[[t^8, t^{12}+t^{14}+t^{15}]]$, so we have
$(p_1,q_1)=(p_2,q_2)=(2,1), (p_3,q_3)=(2,3)$. The dimension of the ASF is $42$ in this case. Using Sage, we compute 
\begin{eqnarray*}
&\langle \fgamma,h_8\rangle=\\
&q^{42} + 7 q^{41} + 24 q^{40} + 56 q^{39} + 104 q^{38} + 166 q^{37} + 236 q^{36} + 306 q^{35} + 370 q^{34} +\\
&424 q^{33} + 465 q^{32} + 492 q^{31} + 507 q^{30} + 510 q^{29} + 504 q^{28} + 488 q^{27} + 466 q^{26} + \\
&437 q^{25} + 406 q^{24} + 370 q^{23} + 335 q^{22} + 298 q^{21} + 264 q^{20} + 230 q^{19} + 199 q^{18} + \\
&168 q^{17} + 143 q^{16} + 118 q^{15} + 97 q^{14} + 78 q^{13} + 63 q^{12} + 48 q^{11} + 38 q^{10} + \\
&28 q^{9} + 21 q^{8} + 15 q^{7} + 11 q^{6} + 7 q^{5} + 5 q^{4} + 3 q^{3} + 2 q^{2} + q + 1
\end{eqnarray*}
which by Theorem \ref{thm:weightpolynomials} is the weight polynomial of the compactified Jacobian in this case. We refer the reader to the attached computer program for computing the Shalika germs and other data in this case.
\end{example}
\begin{example}
Let $G=GL_4$ and $\gamma=u^3\oplus(-u^3)\in\fg\fl_2(F)\times\fg\fl_2(F)\subset\fg\fl_4(F)$, with $u$ the matrix \eqref{eq:umatrix} for $\fg\fl_2$. This is an element whose characteristic polynomial is $x^4-\varpi^6$, so that the link is a $(6,4)$-torus link. (Note that $u^6\in\fg\fl_4(F)$ itself is not regular semisimple, which is why we pass to the isoclinic sum of two slope-$3/2$ blocks with the sign twist.) The element $\gamma$ already lies in a Levi isomorphic to $GL_2\times GL_2$, and on each of the blocks we have an equivalued element of valuation $3/2$. We compute the orbital symmetric function to be the product of $\left|\det(\ad(\gamma)_{\fg/\fm})\right|^{-1/2}=q^6$ and the two factors in this case, namely 
$\fgamma=q^6(qh_{11}+h_{2})^2$.The Shalika expansion of $\fgamma$ reads
 $$\fgamma=\left(\frac{q^{10}}{q^{2} - 2 q + 1}\right)\th_{1111} + \left(\frac{-2 q^8}{q^{2} - 2 q + 1}\right)\th_{211} + \left(\frac{q^{6}}{q^{2} - 2 q + 1}\right)\th_{22}$$
Theorem \ref{thm:orbintformula} gives that $I_\gamma(\1_{(4)})=q^8+2q^7+q^6$ and 
$I_\gamma(\1_{(1^4)})=24q^8+24q^7+6q^6$. Note that up to $q\leftrightarrow t$, the first result agrees with the numerator of \cite[Example 1.3.]{HM} at $a=0, q=1$.
\end{example}
\begin{example}
Let us work out an unramified example. 
Suppose $k=\F_q, p\neq 2$, and $a\in \F_q^\times-(\F_q^\times)^2$. Let 
$$\gamma=\begin{pmatrix}
0 & at \\ t & 0
\end{pmatrix}$$
Then $\gamma$ splits over a degree two unramified extension of $F$. By Hilbert's Theorem 90, stable conjugacy in $GL_n$ is rational conjugacy, so by \cite[(3.5.4)]{ZhiweiLects} we should have 
$$I_\gamma^{GL_n}(\1_{\fg(\cO)})=\mathrm{SO}_\gamma(\1_{\fs\fl_n(\cO)})=q+2$$ where $\mathrm{SO}_\gamma$ is the stable orbital integral in $SL_n$, defined as in \cite[Section 3.5.3.]{ZhiweiLects}.

Indeed, we are in a situation where $1+\gamma\sim 1+tX$ with $X$ generating $\F_{q^2}$ over $\F_q$, and by Theorem \ref{thm:unramified} and Proposition \ref{prop:germcomparison} 
$$\fgamma=q\omega \nabla_{t=1}\omega\tau_2(h_1)|_{q\mapsto q^{-1}}=qh_{11} + 2 h_{2}
$$
By Theorem \ref{thm:orbintformula} we get 
$$I_\gamma(\1_{(2)})=\langle \fgamma, h_2\rangle=q+2$$ as desired.

\end{example}
\begin{example}
Next, we discuss the simplest ''mixed" example. Let $a$ be as above, and $$\gamma=\begin{pmatrix}
    0 & 0 & 0 & at^2 \\
    1 & 0 & 0 & 0\\
    0 & 1 & 0 & 0 \\
    0 & 0 & 1 & 0
\end{pmatrix}$$
This element is elliptic and splits over a degree $4$ extension, with a maximal 
unramified subextension of degree $2$.
The discrete invariants from Proposition \ref{prop:discreteinvariants} are $(f,q,p)=(2,1,2)$, so 
$$\fgamma=q^2\omega \varphi_{1/2}|_{q\mapsto q^{-1}}\omega(\tau_2(e_1))=q^{2}h_{22} + 2 qh_{31} + 2h_{4}$$
and for example $I_\gamma(\1_{(4)})=\langle\fgamma,h_{4}\rangle =q^2+2q+2$. We note that $I_\gamma(\1_{(11111)})=6q^2+8q+2$. Note that the Iwahori affine Springer fiber in this case has six components (cf. Theorem \ref{thm:components}).
\end{example}
\begin{example}
Finally, we consider an example of the form considered in \cite[1.5.]{W1}. Let $f_1=f_2=2$ and $F=F_2\subset F_1\subset F_0$ be a tower of unramified extensions, both of degree 2. Suppose $X_i\in \cO_{F_i}$ are so that their reduction in the residue field generates the residue field over that of $F_{i+1}$. Consider $\gamma=1+tX_1+t^2X_0$. The discrete invariants from Theorem \ref{thm:mainformula} are $(2,1,1),(2,1,1)$ and the orbital symmetric function is 
$$
\fgamma=(q^{8} + 2 q^{6}) h_{1111} + (4 q^{7} + 8 q^{5} + 4 q^{4} + 4 q^{3}) h_{211} + (4 q^{6} + 6 q^{4} + 4 q^{2}) h_{22} + (4 q^{3} + 4 q^{2} + 4 q) h_{31} + 4 h_{4}
$$
Note that on the second step of the recursion we apply $\tau_2$ to a degree two symmetric function. To compare with the transition matrix in Example \ref{ex:twothirteen2} we compute

$$M^{0/2}=\left(\begin{array}{rrrrr}
\frac{2}{q^{5} - q^{3} - q^{2} + 1} & \frac{-2 q^{2} - 2}{q^{4} - q^{3} - q + 1} & \frac{q^{2} + 4 q + 1}{q^{3} - q^{2} - q + 1} & \frac{-4 q}{q^{2} - 2 q + 1} & \frac{q^{2} + q}{q^{2} - 2 q + 1} \\
0 & 0 & \frac{4}{q^{2} - 2 q + 1} & \frac{-4 q - 4}{q^{2} - 2 q + 1} & \frac{q^{2} + 2 q + 1}{q^{2} - 2 q + 1}
\end{array}\right)$$
Note that the $(2),(211)$-entry is {\em not} zero in this example.
This is because the coefficient of $\th_{211}$
in $P_{0,4}$ is $\frac{-4 q^{2} - 4}{q^{3} - 3 q^{2} + 3 q - 1}$ and in $P_{0,2}^2$ it is $\frac{-4 q - 4}{q^{2} - 2 q + 1}$, which gives
$$\frac{1}{2}\left(\frac{-4 q^{2} - 4}{q^{3} - 3 q^{2} + 3 q - 1}(1-q)+\frac{-4 q - 4}{q^{2} - 2 q + 1}(1+q)\right)=\frac{-4 q}{q^{2} - 2 q + 1}$$
\end{example}

\section{Applications}
\label{sec:applications}
\subsection{Affine Springer fibers}
Let $G=GL_n/F$ where $F=k((t))$ with $k=\F_q$. Appropriately modifying the definition of $\Sp_\gamma$ below to account for mixed characteristic $F$, we get similar results but leave these for the interested reader. Suppose $\bP\subset G(F)$ is a parahoric subgroup. Let $\Fl_\bP=G(F)/\bP$ be the corresponding partial affine flag variety.

\begin{definition} Let $\gamma\in\fg(F)^{rs}$.
The affine Springer fiber is the reduced ind-subscheme of $\Fl_\bP$ defined by
$$\Sp_\gamma^\bP(k)=\{g\bP|\Ad(g^{-1})\gamma\in \Lie(\bP)\}$$ Let $T_\gamma$ be the centralizer of $\gamma$. Then it acts naturally on $\Sp^\bP_\gamma$. Let $S_\gamma$ be the maximal unramified subtorus of $T_\gamma$, and $X_*(S_\gamma)=:\Lambda_\gamma$ its cocharacter lattice. As in \cite[Section 15]{GKM}, the centralizer action gives rise to an action of $\Lambda_\gamma$ on $\Sp_\gamma^\bP$.
\end{definition}
Recall that $1_\bP$ is the characteristic function of $\bP$ divided by the measure of $\bP$. Unraveling the definitions, it is not hard to prove (see for example \cite[Theorem 15.8.]{GKM})
\begin{proposition}
\label{prop:numberofpoints} 
We have $$|\Sp_\gamma^\bP(k)/\Lambda_\gamma(k)|=I_\gamma(\1_\bP).$$ See Definition \ref{def:haarmeasures} for the normalization in $I_\gamma$.
\end{proposition}
Recall that Theorem \ref{thm:orbintformula} says that if $\bP$ is of type $\lambda$ and $\gamma$ is tamely ramified and regular semisimple, then
$$I_\gamma(\1_\bP)=\langle \fgamma,h_\lambda\rangle$$ which is a polynomial with nonnegative coefficients. This result, combined with Proposition \ref{prop:numberofpoints} implies
\begin{corollary}
\label{cor:asfcount}
When $\gamma$ is tamely ramified, the number of points $|\Sp_\gamma^\bP(k)/\Lambda_\gamma(k)|$ is a polynomial in $q$ with nonnegative coefficients.
\end{corollary}
Let us also note the following application to components of affine Springer fibers.
\begin{theorem}
\label{thm:components}
Let $\gamma\in \fg(F)$ be compact and regular semisimple. 
Then the number of geometric components of $\Sp_\gamma^{\bI}/\Lambda_\gamma$ stable under $Gal(\overline{k}/k)$ is always a divisor of  $|W|=n!$.
\end{theorem}
\begin{proof}
By \cite[Eq. (4.5)]{Tsai20}, the number of Frobenius-stable geometric components is the coefficient of the leading term in $q$ of the integral of $\1_{11\cdots 1}$ along the orbit of $\gamma$. By Theorem \ref{thm:orbintformula}, this orbital integral can be computed using $\fgamma$ by pairing it with $h_{11\cdots 1}$. On the other hand, $\fgamma$ is formed by multiplying the orbital symmetric functions for the blocks of $\gamma$. Suppose for a moment $\gamma$ is totally ramified. By Lemma \ref{lem:smallestpower} the smallest power in the Dyck germs is $1$, it appears with coefficient one, and it appears for the least dominant partition. Since pairing with $h_\lambda$ does not introduce powers of $q$, the highest power of $q$ appearing in 
$$\langle \fgamma,h_{11\cdots 1}\rangle$$ is $q^{\dim \Sp_\gamma}$ and it appears with coefficient $\langle h_\lambda,h_{11\cdots 1}\rangle=\frac{|S_n|}{|S_{\lambda}|}$
where $\lambda$ is the smallest partition in the dominance order appearing in the Dyck expansion of $\fgamma$. This proves the claim in the totally ramified case. 

For the general case where the construction $\fgamma$ involves the operator $\tau_f$, cf. Proposition \ref{prop:unramifiedextension}, note that the plethysm $\tau_f$ is designed so that for any homogeneous symmetric function $\mathbf{f}$ of degree $n'=n/f$, $\lambda\vdash n$, we have
$$\langle h_\lambda,\tau_f(\mathbf{f})\rangle=\begin{cases}
\langle h_{\lambda/f}, \mathbf{f}\rangle, & \text{ if $\lambda$ is divisible by $f$}\\
0, & \text{ if $\lambda$ is not divisible by $f$}
\end{cases}$$
In particular, we may reduce these cases to the computation above. 
\end{proof}

\begin{lemma}
\label{lem:smallestpower}
Let $\gamma$ be inertially elliptic, with Newton pairs $(\vec p,\vec q)$. Among the $k$-nested Dyck paths of Definition \ref{def:nesteddyck} there is a unique path $\Pi_0\in\nDyck_{(\vec p,\vec q)}$ of area $0$, obtained by assigning to each incoming run the unique area-zero Dyck path in its box. Consequently the smallest power of $q$ in the coefficients of $\fpq=\sum_\lambda\sigma_\lambda(\gamma)e_\lambda$ is $q^0$; it occurs with coefficient $1$, only in front of $e_{\lambda(\Pi_0)}$, and $\lambda(\Pi_0)$ is the smallest in dominance order among all $\lambda$ with $\sigma_\lambda(\gamma)\neq0$.
\end{lemma}
\begin{proof}
By Proposition \ref{prop:annularformula}, $\fpq=\sum_{\Pi\in\nDyck_{(\vec p,\vec q)}}q^{\area(\Pi)}e_{\lambda(\Pi)}$ with $\area(\Pi)=\sum_{\pi\in\Pi}\area(\pi)\ge0$, and $\area(\Pi)=0$ iff every constituent path has area $0$. In each box $\Dyck_{hq_j,hp_j}$ there is exactly one area-$0$ Dyck path, namely the maximal staircase lying just under the diagonal; choosing it at every level yields the unique area-$0$ nested path $\Pi_0$. Hence $q^0$ occurs in $\fpq$ exactly once, with coefficient $1$, in front of $e_{\lambda(\Pi_0)}$. Finally $\lambda(\Pi_0)$ is minimal in the dominance order: at each level the area-$0$ path has the most evenly spread horizontal runs, and dominance order is preserved under the nesting.
\end{proof}

\begin{remark}
Let $m/n$ be the minimal root valuation of $\gamma$. When $\gamma$ is inertially elliptic, the above shows that the minimal partition appearing in the $e_\lambda$-expansion is formed from the horizontal steps of the maximal staircase partition fitting under a line of slope $m/n$, i.e. the one with parts $\lfloor \frac{(m-k)n}{m}\rfloor, k=1,\ldots m$. For example, when $m/n=3/7$ this gives the partition $4+2(+0)$, and the corresponding Dyck path in the $3\times 7$ rectangle has horizontal steps $3, 2, 2$. In particular when $m/n\geq 1$ the horizontal steps give the one-column partition.
It is easy to extend this to $\gamma$ non-elliptic by multiplying the corresponding $e_\lambda$ together.
This gives another (slightly more general) proof of a Theorem of Z. Yun in type A, which states that the {\em minimal reduction type} of $\gamma$ determines the number of components in the Iwahori affine Springer fiber.
\end{remark}
\begin{remark}
Theorem \ref{thm:components} proves \cite[Conjecture 8.7.]{Tsai20} in type A.
From the main result of \cite{Tsai20}, there are always exactly $n!$ components when the depth is $>1$. In fact the last statement is true for depth $\ge 1$ because any depth-$1$ element either differs from a depth $>1$ element by a central element, or is contained in a Levi subalgebra in which case we can reduce the assertion to the Levi case as in the proof of Theorem \ref{thm:orbintformula}.
\end{remark}
Theorem \ref{thm:components} has the following interesting corollary about the $W$-representation given by $H^*(\Sp_\gamma^{\bI}/\Lambda_\gamma)$. Let us assume that $\gamma$ is tamely ramified, and note that the top degree part of the cohomology is always pure. A well-known argument using finite Springer theory tells us there is a graded isomorphism of vector spaces:
$$H^*(\Sp_\gamma^{\bP_\lambda}/\Lambda_\gamma)\cong H^*(\Sp_\gamma^\bI/\Lambda_\gamma)^{W_\lambda}$$
In particular, knowing the dimensions of the top degree cohomologies of each $\Sp_\gamma^{\bP_\lambda}/\Lambda$ tells us exactly all the dimensions of the $W_\lambda$-invariants of the representation on top degree cohomology of $\Sp^\bI_\gamma/\Lambda_\gamma$. Recall that using the Hall inner product and Frobenius reciprocity, this is the same as knowing the inner products of the Frobenius character with $h_\lambda$. Since the $h_\lambda$ are a basis of the ring of symmetric functions, this uniquely determines the representation. A similar argument shows that assuming purity, $\fgamma$ actually determines the Frobenius character of $H^*(\Sp^\bI_\gamma/\Lambda_\gamma)$ in the elliptic case. In fact, the Frobenius character will simply be $\fgamma$ if these assumptions are satisfied.

More precisely, we get 
\begin{theorem}
\label{thm:topdegreerep}
Let $\gamma$ be tamely ramified. The $W=S_n$-representation on $H^{top}(\Sp^\bI_\gamma/\Lambda_\gamma)$ has Frobenius character $h_\nu$, where $\nu$ is the smallest partition in dominance order appearing in the $h_\lambda-$expansion of $\fgamma$.
In particular, when $\gamma$ has depth $\geq 1$, this is the regular representation by above.

If we assume the purity conjecture and that $\gamma$ is elliptic, $\fgamma$ is the Frobenius character of $H^*(\Sp^\bI_\gamma/\Lambda_\gamma)$.
\end{theorem}
\begin{proof}
For the first statement, let $C_\gamma$ be the number of components in the Iwahori affine Springer fiber of $\gamma$, defined as in Theorem \ref{thm:components}. As in the proof of Theorem \ref{thm:components}, the orbital symmetric function has the form 
$$\fgamma=C_\gamma q^{\Xi(\gamma)} h_\nu+O(q^{\Xi(\gamma)-1})$$ Pairing this with $h_\lambda$ for varying $\lambda\vdash n$ and taking the leading term in $q$ gives the trace of Frobenius on $H^{top}(\Sp_\gamma^{\bP_\lambda}/\Lambda_\gamma)$ for varying $\lambda$. But since the top cohomology is always pure, the coefficient of $q^{\Xi(\gamma)}$ in $\langle \fgamma, h_\lambda\rangle$ in fact equals the top Betti number of $\Sp_\gamma^{\bP_\lambda}/\Lambda_\gamma$. Since an $S_n$-representation is uniquely determined by the dimensions of its invariants under Young subgroups and these are given exactly by $\langle h_\nu, h_\lambda\rangle$ in this case, we are done.

For the second statement, if purity holds we see $\langle \fgamma, h_\lambda\rangle$ is the ordinary Poincar\'e polynomial of $\Sp_\gamma^{\bP_\lambda}/\Lambda_\gamma$, in particular the graded dimension of the space of $S_\lambda$-invariants of $H^*(\Sp_\gamma^\bI/\Lambda_\gamma)$. Since these pairings determine $\fgamma$ uniquely, the graded Frobenius character of $H^*(\Sp_\gamma^\bI/\Lambda_\gamma)$ equals $\sum_\lambda \langle \fgamma,h_\lambda\rangle m_\lambda=\fgamma$.
\end{proof}
\begin{remark}
Note that this proves \cite[Conjecture 7.17.]{GKO} in type A.
\end{remark}

\begin{remark}
\label{rmk:purity}
The purity conjecture assumed in the second part of Theorem \ref{thm:topdegreerep} is known for several classes of inertially elliptic elements, for example the ''homogeneous" ones. Combined with our results, this gives an unconditional identification of $\fgamma$ with the Frobenius character of $H^*(\Sp^\bI_\gamma/\Lambda_\gamma)$ in those cases. 
\end{remark}


\subsection{Compactified Jacobians}
\label{sec:cptjac}
In this section, we apply Theorem \ref{thm:weightpolynomials} to show that the point-counts of compactified Jacobians of rational, unibranch plane curves are polynomials in $q$.

Let us recall some relevant material from \cite{Laumon}. Let $C$ be a reduced, projective and geometrically connected curve over the residue field $k$, with only planar singularities. Suppose for simplicity that the normalization of $C$ is rational. Let $\bPic(C)$ be the compactified Picard scheme of $C$. It is the moduli space whose closed points parametrize torsion-free rank one sheaves on $C$. For each $c\in \text{Sing}(C)$ fix an isomorphism $\widehat{\cO}_{C,c}\cong k[[x,y]]/f$ and let $\Sp_c$ be the affine Springer fiber associated to $\gamma_c:=\gamma_f\in \fg\fl_{\deg_x f}$ where $\gamma_f$ is the companion matrix of $f$. Let $\Lambda_c$ be the lattice part of the centralizer of $\gamma_f$ and $\Lambda=\Pic(C)/\Jac(C)$. Fix a section $\Lambda\to \Pic(C)$ of the quotient map.

From \cite[Proposition 2.3.1.]{Laumon} we have 
\begin{proposition}
\label{prop:gkmlaumon}
There is a universal homeomorphism $$\prod_{c\in \text{Sing}(C)} \Sp_c/\Lambda_c \to \bPic(C)/\Lambda$$
\end{proposition}
If $k$ is a finite field, we have
\begin{corollary}
\label{cor:comparecounts}
Let $k'/k$ be a finite extension. Then 
$$\left|\prod_{c\in\text{Sing}(C)}\Sp_c(k')/\Lambda_c\right|=\left|\bPic(C)(k')/\Lambda\right|$$
\end{corollary}
Combined with Corollary \ref{cor:asfcount}, we have 
\begin{theorem}
\label{thm:polynomialcount}\label{cor:nonnegative}
The number of points on $\bPic(C)$ is a polynomial in $q=|k|$. In addition, it is a polynomial with nonnegative integer coefficients.
\end{theorem}
A standard spreading-out argument, combined with \cite[Theorem 1]{Katz} and the previous Theorem gives
\begin{corollary}
Let $k=\C$. Then $X=\bPic(C)$ is strongly polynomial-count in the sense of Katz \cite{Katz}, and the
$E$-polynomial $$E_X(x,y):=\sum_{p,q} e_{p,q} x^py^q$$ is given by the 
weight polynomial of $\bPic(C)$ as $E_X(x,y)=P_X(xy)$, defined by
$$P_X(q)=\sum_{i,j}(-1)^iq^j\dim \gr_{W}^jH^i(\bPic(C))$$
\end{corollary}
Yet another corollary of Corollary \ref{cor:asfcount} together with Corollary \ref{cor:comparecounts} and Definition \ref{def:superpolynomial} is a virtual version of 
\cite[Conjecture 2.4.(iii)]{CD}, which compares Betti numbers of Jacobian factors with superpolynomials at $q=1$ (i.e. $t=1$ in our notation). Using Definition \ref{def:superpolynomial}, the more precise statement is that
\begin{proposition}
\label{prop:superpolycptjac}
For unibranch $C=\{f(x,y)=0\}\subset \C^2$, the weight polynomial of $\overline{\mathrm{Jac}}(C)$ is given by the superpolynomial at $a=0, t=1$, with $q$ replaced by $1/q$, up to multiplying by $q^{\dim \overline{\mathrm{Jac}}(C)}$. That is,
$$P_{\overline{\mathrm{Jac}}(C)}(q)=q^{\dim \overline{\mathrm{Jac}}(C)}\mathbf{P}_{\text{Link}(C)}(a=0,q=1/q,t=1)$$
\end{proposition}

\subsection{Bounds on orbital integrals}\label{sec:boundsOI}
Our explicit formulas have consequences for the size of regular semisimple orbital integrals on $GL_n$, a question of independent interest in analytic number theory.

Shin and Templier \cite{ST} prove a Sato--Tate equidistribution theorem for families of automorphic representations, whose proof requires a uniform bound on orbital integrals of unramified Hecke operators \cite[Proposition 7.1]{ST}. That bound has the shape
$$0\le I_\gamma(\tau_\lambda)\le q^{a+b\Vert\lambda\Vert}\,D_G(\gamma)^{-e/2},$$
where $D_G$ is the Weyl discriminant and the constants $a,b,e\ge 0$ are shown to exist but are not made explicit.

For tamely ramified regular semisimple $\gamma$, we evaluate the integrals $I_\gamma(\1_\lambda)$ explicitly, rather than bounding them. For more general Hecke functions $\tau_\lambda$, knowledge of the germs reduces computing $I_\gamma(\tau_\lambda)$ to understanding unipotent orbital integrals or truncated Steinberg evaluations of $\tau_\lambda$.  

In a different direction, Langlands' Beyond Endoscopy proposal \cite{Beyond} requires explicit evaluation of orbital integrals. For $GL_2$, this is carried out using elementary methods by Langlands and extended by Altug \cite{Altug} and Espinosa Lara \cite{MEEL}. For $GL_n$ with $n>2$, Arthur \cite{Arthur} suggests that explicit local orbital integrals could yield analogous results; our formulas provide the necessary local input in the tamely ramified case.

\end{document}